\documentclass[a4 paper, 11pt]{smfart}

\usepackage{smfthm}
\usepackage{smfenum}
\usepackage{amssymb}
\usepackage{amscd}

\usepackage{textcomp}  

\usepackage{amssymb,amsmath,latexsym,epsfig,euscript}

\catcode`\@=11

\newdimen\f@ntdim
\def\ovg{\f@ntdim=\fontdimen6\font \font\f@nt=lasy10 at \f@ntdim
\ouvreg\f@nt\everypar={\ouvreg\f@nt}}
\def\ogg{\f@ntdim=\fontdimen6\font \font\f@nt=lasyb10 at \f@ntdim
\ouvreg\f@nt\everypar={\ouvreg\f@nt}}
\def\ouvreg#1{\leavevmode\hbox{#1\unskip\kern+0.20em(\kern-0.20em(\kern+0.20em}%
\ignorespaces\nobreak}
\def\fmg{\fermeg\f@nt\everypar={}} 
\def\fermeg#1{\relax \ifhmode \unskip\kern+0.20em\else \leavevmode \fi%
\hbox{#1)\kern-0.20em)\kern+0.20em\ignorespaces}}

\catcode`\@=12 

\def\enoncesans#1#2 #3{\begin{thm}#1#3\end{thm}}
\def\enoncenormal#1#2 #3{\begin{prop}#1#3\end{prop}}
\def\cal#1{\mathcal{#1}}

\font\ninebf=cmbxti10 at 9pt

\def\^#1{\if#1i{\accent"5E\i}\else{\accent"5E #1}\fi} 
\def\"#1{\if#1i{\accent"7F\i}\else{\accent"7F #1}\fi} 
\catcode`\Ž=\active \def Ž{\'e}
\catcode`\ˆ=\active \def ˆ{\`a}
\catcode`\=\active \def {\`e}
\catcode`\=\active \def {\`u}
\catcode`\‰=\active \def ‰{\^a}
\catcode`\=\active \def {\^e}
\catcode`\"=\active \def "{\^\i}
\catcode`\™=\active \def ™{\^o}
\catcode`\ž=\active \def ž{\^u}
\catcode`\Š=\active \def Š{\"a}
\catcode`\'=\active \def '{\"e}
\catcode`\•=\active \def •{\"\i}
\catcode`\š=\active \def š{\"o}
\catcode`\Ÿ=\active \def Ÿ{\"u}
\catcode`\=\active \def {\c c}

\catcode`\é=\active \def Ž{\'e}
\catcode`\à=\active \def ˆ{\`a}
\catcode`\è=\active \def {\`e}
\catcode`\ù=\active \def {\`u}
\catcode`\â=\active \def ‰{\^a}
\catcode`\ê=\active \def {\^e}
\catcode`\î=\active \def "{\^\i}
\catcode`\ô=\active \def ™{\^o}
\catcode`\û=\active \def ž{\^u}
\catcode`\ä=\active \def Š{\"a}
\catcode`\ë=\active \def '{\"e}
\catcode`\ï=\active \def •{\"\i}
\catcode`\ö=\active \def š{\"o}
\catcode`\ü=\active \def Ÿ{\"u}
\catcode`\ç=\active \def {\c c}

\def\enmarged#1{\goodbreak\ifodd\count0\marginpar{\hbox{\ninebf#1\kern18mm}}
\else\marginpar{\hbox{\kern18mm\ninebf#1}}\fi}

\catcode`\é=\active \def é{\'e}
\catcode`\à=\active \def à{\`a}
\catcode`\À=\active \def À{\`A}
\catcode`\è=\active \def è{\`e}
\catcode`\ù=\active \def ù{\`u}
\catcode`\â=\active \def â{\^a}
\catcode`\ê=\active \def ê{\^e}
\catcode`\î=\active \def î{\^\i}
\catcode`\ô=\active \def ô{\^o}
\catcode`\û=\active \def û{\^u}
\catcode`\ä=\active \def ä{\"a}
\catcode`\ë=\active \def ë{\"e}
\catcode`\ï=\active \def ï{\"\i}
\catcode`\ö=\active \def ö{\"o}
\catcode`\ü=\active \def ü{\"u}
\catcode`\ç=\active \def ç{\c c}

\def\codim {{\rm \mbox{codim}}}

\def\Conv {{ \mbox{\rm Conv}}}

\def\Vol {{\mbox{\rm Vol}}}

\def\Ord {\mbox{\rm ord}}

\def\dist {{\rm \mbox{dist}}}
\def\init {{\rm \mbox{init}}}
\def\initial {{\mbox{\scriptsize \rm init}}}

\def\rang {\mbox{\rm rang\,}}
\def\ess {{ \mbox{\rm \scriptsize ess}}}
\def\abs {{ \mbox{\rm \scriptsize abs}}}

\def\Div {\mbox{\rm Div}}
\def\div {\mbox{\rm div}}
\def\Tors {\mbox{\rm Tors}}

\def\irr {\mbox{\rm \scriptsize Irr}}

\def\geom {{\mbox{\rm \scriptsize g{\'e}om}}}

\def\arith {{\mbox{\rm \scriptsize arith}}}

\def\Spec   {\mbox{\rm Spec}\,}
\def\MV {\mbox{\rm MV}}
\def\MI {\mbox{\rm MI}}

\def\Card {\mbox{\rm Card}}
\def\card {{\mbox{\rm \scriptsize Card}}}

\def\Expo {\mbox{\rm Supp}}
\def\expo {\mbox{\rm \scriptsize Supp}}

\def\Spec {\mbox{\rm Spec}}
\def\DPC  {\mbox{\rm DPC}}
\def\dpc  {\mbox{\rm \scriptsize DPC}}
\def\Pan  {\mbox{\rm Pan}}
\def\pan  {\mbox{\rm \scriptsize Pan}}
\def\longueur  {\mbox{\rm long}}

\def\hnorm {{\widehat{h}}}
\def\munorm {{\widehat{\mu}}}

\def\geom {{\mathcal A}}
\def\arith {{{\mathcal A}, \alpha}}

\def\Faces {{\mbox{\rm F}}}
\def\faces {{\mbox{\rm \scriptsize F}}}
\def\res {\mbox{\rm r{\'e}s}}

\def\ov#1{{\overline{#1}}}
\def\un#1{{\underline{#1}}}
\def\wh#1{{\widehat{#1}}}
\def\wt#1{{\widetilde{#1}}}

\def \A {{\bf A}}
\def \C {{\bf C}}
\def \G {{\bf G}}
\def \K {{\bf K}}
\def \N {{\bf N}}
\def \P {{\bf P}}
\def \Q {{\bf Q}}
\def \R {{\bf R}}
\def \Z {{\bf Z}}

\def \cA {\mathcal{A}}
\def \cB {\mathcal{B}}

\def \cM {\mathcal{M}}

\def \cO {\mathcal{O}}

\def \cX {\mathcal{X}}

\def\cad {c'est-{\`a}-dire } 
\def\spdg {sans perte de généralité }

\def\Qbar {{\overline{\Q}}}

\title[Quelques aspects diophantiens des vari{\'e}t{\'e}s toriques projectives]{
Quelques aspects diophantiens des\\ vari{\'e}t{\'e}s toriques projectives}

\author{Patrice Philippon} 

\address{Institut de Math{\'e}matiques de Jussieu (U.M.R. 7586), Projet 
G{\'e}om{\'e}trie et Dynamique.
Case~7012, 
2 place Jussieu, 
75251 Paris Cedex 05,
France.}
\email{pph@math.jussieu.fr}
\urladdr{http://www.math.jussieu.fr/\~{}pph/}
\thanks{P. Philippon a {\'e}t{\'e} partiellement financ{\'e} par une allocation de recherche de la {\em Fondation Alexander von Humboldt} pendant la r{\'e}alisation de ce travail.} 

\author{Mart{\'\i}n Sombra}

\address{Universitat de Barcelona, 
Departament d'{\`A}lgebra i Geometria. 
Gran Via 585,
08007 Barcelona, Espagne.}
\email{sombra@ub.edu}
\urladdr{http://atlas.mat.ub.es/personals/sombra/} 
\thanks{ M. Sombra a \'et\'e financ\'e par le {\it Programme Ram\'on y Cajal} du  
Ministerio de Educaci\'on y Ciencia, Espagne.}

\subjclass{Primaire: 11G50; Secondaire: 14M25, 14G40.} 

\keywords{Vari{\'e}t{\'e} torique, hauteur normalis{\'e}e, multihauteurs, 
fonction de Hilbert arithm{\'e}tique, poids de Chow, volume mixte, indice
d'obstruction, minimums successifs.}

\begin{document}
\maketitle

\vspace{-0.5cm}
\hbox to\hsize{\hfill\vbox{\hsize=5.5cm\parskip=0pt{\small\it 
An Wolfgang Schmidt\par
Der oft und schick\par
'nen Edelsatz\par
Geschmiedet hat}}}
\vspace{1cm}


\begin{abstract}
\setlength{\baselineskip}{12pt}
On prsente plusieurs facettes des vari{\'e}t{\'e}s toriques projectives, intressantes du point de vue 
de la gomtrie diophantienne. On montre comment la thorie s'explicite sur un certain nombre d'exemples 
significatifs et on tablit galement un thorme de type Bzout pour les poids de Chow 
des varits projectives. 
\end{abstract}


\begin{altabstract}
\setlength{\baselineskip}{12pt}
{\bf Some diophantine aspects of projective toric varieties.} 
We present several facets of projective toric varieties, of interest from the point of view of Diophantine geometry. 
We make explicit the theory in a number of meaningful examples and we also prove a Bzout type theorem for Chow weight of projective varieties.

\end{altabstract}




\setlength{\baselineskip}{13.55pt}
 
%
%
 

\typeout{Introduction}

\section*{Introduction et résultats}

Les vari{\'e}t{\'e}s toriques jouent un r{\^o}le important au carrefour de l'alg\`ebre, la g{\'e}om{\'e}trie et la combinatoire. 
Elles constituent une classe de variétés suffisamment rigide pour que
beaucoup des invariants s'explicitent en termes combinatoires, 
et en même temps suffisamment riche pour permettre de tester et illustrer 
diverses conjectures et thŽories abstraites. 
Elle trouve application 
dans de nombreuses branches des mathŽmatiques~: gŽomŽtrie algŽbrique bien sžr, algbre commutative, combinatoire, calcul formel, gŽomŽtries symplectique et kŠhlerienne, topologie et physique mathŽmatique,
{\it voir} par exemple~\cite{Ful93}, \cite{GKZ94}, \cite{Stu96}, \cite{Cox01},
\cite{Audin}, 
\cite{Don02}.

Par définition, les variŽtŽs toriques projectives sont 
les com\-pac\-ti\-fi\-ca\-tions Žquivariantes des translatŽs de sous-tores 
des tores multiplicatifs $\G_m^N$. 
Du point de vue de la gŽomŽtrie diophantienne, ces vari{\'e}t{\'e}s se 
trouvent ˆ la croisŽe des problmes de Lehmer et de Bogomolov g\'en\'eralis\'es sur les tores. 
En effet, on sait que lorsqu'elles ne sont pas de torsion,  
les minorations pour la hauteur normalisŽe de ces variŽtŽs 
sont de nature fondamentalement arithmŽtique, 
dŽpendant essen\-tiel\-le\-ment du corps de dŽfinition de la variŽtŽ. 
Au contraire, pour les sous-variŽtŽs de $\G_m^N$ qui ne sont pas 
des translatŽes de sous-tores on dispose de minorations ne dépendant que de leur géométrie, {\it voir}~\cite{AD}, \cite{AD04}, \cite{Dav02}. 
D'un autre c™tŽ, R.~Ferretti~\cite{Fer03} a exploit\'e les variŽtŽs toriques 
afin de trouver des exemples concrets de l'extension 
du thŽorme du sous-espace 
aux variŽtŽs projectives, qu'il a obtenue avec 
J.-H.~Evertse~\cite{EF02}.

Dans~\cite{PS04} ({\it voir} aussi~\cite{PS05})
nous avons ŽtudiŽ l'un des invariants 
arith\-mŽ\-ti\-ques les plus significatifs des variŽtŽs dans le cas torique, à savoir 
leur hauteur normalisŽe. 
Cet invariant est l'analogue arithmétique du degré, il mesure la complexité binaire d'une 
représentation de la variété et contrôle aussi la distribution des points algébriques 
de petite hauteur sur la vari\'et\'e. 
Dans \cite{PS04}, on a donn\'e une expression explicite pour la hauteur 
normalisée d'une variété torique et 
plus généralement pour la multihauteur d'un tore par rapport à plusieurs plongements 
monomiaux. 

Ces résultats sont en parfait parallle avec la  thŽorie gŽomŽtrique. 
En fait, on construit un objet adélique $\Theta_X$ associé à une variété torique $X$ 
(constitué par une famille finie de fonctions concaves et affines par morceaux)
qui est le pendant arithmétique
du polytope classiquement associé à l'action du tore et dont l'intégrale donne la hauteur. 
Gr‰ce ˆ cette approche, il est possible de calculer explicitement
cette quantitŽ 
pour n'importe qu'elle variŽtŽ torique particulire et de tester utilement des conjectures et rŽsultats.

Le prŽsent texte a le double propos d'introduire le lecteur 
ˆ l'Žtude des 
variŽtŽs toriques 
dŽbutŽe dans~\cite{PS04} (\S~\ref{Chow}, \S~\ref{hauteurs}), 
et de présenter de nouvelles applications des variétés toriques ˆ 
des problmes diophantiens ou 
d'origine diophantienne (\S~\ref{obstruction}, \S~\ref{voltang}, \S~\ref{minimums}, \S~\ref{diviseurs}).

\medskip

Dans le~\S~\ref{obstruction} on s'intŽresse aux indices d'obstruction successifs 
des variŽtŽs toriques d\'efinies sur un corps alg\'ebriquement clos $\K$. 
Il s'agit des plus petits degrŽs de formes d'une suite sŽcante (soit globalement, soit dans un ouvert) dŽcoupant un ensemble algŽbrique ayant la variŽtŽ comme composante. 
Différentes variantes de ces indices jouent un rôle important dans les généralisations des problèmes de Lehmer 
et de Bogomolov par exemple, {\it voir}~\cite{AD}.

Les sous-vari{\'e}t{\'e}s toriques de $\P^N(\K)$ correspondent à des id{\'e}aux {binomiaux} premiers et homog{\`e}nes de l'anneau $\K[x_0, \dots,x_N]$. 
De plus, les bin{\^o}mes engendrant l'id{\'e}al d'une vari{\'e}t{\'e} torique $X$ 
s'explicitent en termes d'un certain $\Z$-module $\Gamma_X \subset \Z^{N+1}$ 
naturellement associ{\'e} à $X$, {\it voir}~\cite{ES96} ou~\S~\ref{obstruction}.

On montre que le premier indice d'obstruction d'une variété torique $X$ 
est Žgal au premier minimum de $\Gamma_X$ par rapport {\`a} une m{\'e}trique convenable; et 
plus g{\'e}n{\'e}ralement, que les indices d'obstruction successifs  
$\omega_i(X; (\P^N)^\circ)$ de cette vari{\'e}t{\'e} 
relatifs ˆ l'ouvert $(\P^N)^\circ$ 
coïncident avec les minimums successifs de $\Gamma_X$ et qu'ils 
se r{\'e}alisent par des Žquations binomiales (Proposition~\ref{framboise}). 
{\it Via} ce r{\'e}sultat, 
le deuxi{\`e}me th{\'e}or{\`e}me de Minkowski se traduit en des estimations pour
le produit des indices d'obstruction successifs, qui pr{\'e}cisent dans
le cas torique les estimations de M.~Chardin~\cite{Cha89} et de Chardin et 
P.~Philippon~\cite{CP99}:

\begin{prop} \label{intro2}
Soit $X\subset \P^N$ une variété torique de dimension $n$, alors 
$$
\deg(X) \le \omega_1(X; (\P^N)^\circ)
\cdots  
\omega_{N-n}(X; (\P^N)^\circ)
\le (N+1)^{N-n} \, \deg(X).
$$
\end{prop}

En outre, le rŽseau $\Gamma_X$ s'identifie au rŽseau des pŽriodes de l'application exponentielle 
restreinte ˆ l'espace tangent en l'origine de $X^\circ$. On retrouve ainsi au~\S~\ref{voltang} certains rŽsultats de~\cite{BP} reliant degrŽ et multi-degrŽs d'un sous-groupe algŽbrique d'un tore multiplicatif au volume de son rŽseau des pŽriodes et ˆ la hauteur de son espace tangent.

Dans l'article~\cite{PS06} on poursuit une étude approfondie des indices d'obstruction des variétés toriques et de son application à la minoration de la hauteur des points dans ces variétés.

\medskip

Soit maintenant $X\subset{\bf P}^N$ une vari{\'e}t{\'e} quelconque de dimension $n$ 
et $\tau=(\tau_0,\dots,\tau_N)\in{\bf R}^{N+1}$ un vecteur {\it poids}.
Soient 
$n+1$ groupes $U_0, \dots, U_n$ de $N+1$ variables chacun et 
considérons 
la {\it forme de Chow de $X$}
$$
Ch_X =\sum_{a\in \N^{(n+1)(N+1)}} c_a U_0^{a_0}\cdots U_n^{a_n} 
\in  \K[U_0, \dots, U_n]
$$ 
Le {\em poids de Chow relatif ˆ $\tau$} (ou {\em $\tau$-poids de Chow}) de $X$ est
défini comme le poids de sa forme de Chow par rapport au vecteur $(\tau,\dots, \tau)\in \R^{(n+1)(N+1)}$, 
\cad
$$
e_\tau(X) := \max \{ \langle a_0,\tau\rangle + \cdots +\langle a_n,\tau\rangle :
a\in \N^{(n+1)(N+1)} \mbox{ tel que } c_a\ne 0 \}\enspace,$$
où $\langle\cdot,\cdot\rangle$ désigne le produit scalaire ordinaire. 

Le poids de Chow a ŽtŽ introduit par D. Mumford~\cite{Mum} en
relation avec la stabilitŽ des variŽtŽs projectives. On le
retrouve (et en particulier l'ŽnoncŽ~\ref{poidsChow} du
prŽsent texte) dans un travail de S.K. Donaldson~\cite{Don02} montrant la relation 
entre stabilitŽ des variŽtŽs toriques et existence de mŽtrique 
kŠhlerienne ˆ courbure constante. 
Il apparaît Žgalement en gŽomŽtrie diophantienne au travers des travaux de Ferretti~\cite{Fer03}, Evertse et Ferretti~\cite{EF02} et dans notre formule pour  la hauteur d'une variété torique~\cite{PS04}.

\smallskip

Pour $\tau \in \Z^{N+1}$, consid{\'e}rons l'action du sous-groupe {\`a} un param{\`e}tre   
$$
*_\tau : \G_m \times \P^N \to \P^N \enspace, \quad \quad 
(t,(x_0: \cdots : x_N)) \mapsto (t^{\tau_0}\, x_0 : \cdots :
t^{\tau_N}\, x_N)
$$
et la {\it dŽformation torique} $ X_\tau$ de $X$ associ{\'e}e, d{\'e}finie comme l'adhŽrence de Zariski de l'ensemble 
$$
\{((1:t),  t *_\tau x) \, : \ t\in{\bf G}_m,\ x\in X\}
\subset{\bf P}^1\times{\bf P}^{N} \enspace. 
$$
La {\it variété initiale de $X$ relative au poids  $\tau \in \Z^{N+1}$}
est par définition
$$
\init_\tau(X) := \iota^*(X_\tau \cdot (\{(0:1)\}\times \P^{N})) \in
Z_n(\P^N) \enspace, 
$$
o{\`u} $ \iota: \P^N \to \P^1\times \P^N $ d{\'e}signe l'inclusion $(x_0:\cdots:x_N) \mapsto ((0:1), (x_0:\cdots:x_N))$. Autrement-dit, $\init_\tau(X)$ est le {\it cycle limite} $\lim_{t\to \infty} t*_\tau X$ de $X$ sous l'action $*_\tau$. C'est un cycle de m{\^e}me dimension et degr{\'e} que $X$.

On montre au~\S~\ref{Chow} que lorsque $\tau\in\N^{N+1}$, le poids de Chow relatif ˆ $\tau$ s'interprte comme un bi-degrŽ d'une variante de la dŽformation torique ci-dessus et se comporte donc comme une hauteur. 
Comme consŽquence de cette interprŽtation, on dŽmontre un thŽorme de BŽzout pour le poids de Chow
du cycle intersection $X\cdot H$, 
 qui prŽcise la majoration obtenue par Ferretti~\cite[Prop.~4.3]{Fer03}.

\begin{thm}\label{BezoutpoidsChow}
Soit $X\subset\P^N$ une variŽtŽ projective et $ H \in\Div(\P^N)$ un diviseur ne contenant pas $X$, alors pour $\tau\in\Z^{N+1}$ on a 
$$\begin{array}{rcl}
e_\tau(X\cdot H) &= &e_\tau(X)\,\deg(H) + e_\tau(H)\,\deg(X) - (\tau_0+\dots+\tau_N)\,\deg(H)\,\deg(X)\\[2mm] 
&&\kern5.3cm\displaystyle - \sum_{Y\in\irr(\initial_\tau(X))} m(X_\tau\cdot H_\tau;\iota(Y)) \, \deg(Y)
\end{array}$$
o{\`u} la somme porte sur les composantes irr{\'e}ductibles de ${\rm init}_\tau(X)$.

En particulier, si $H$ est effectif on a pour tout $\tau\in\R^{N+1}$
$$
e_\tau(X\cdot H) \leq e_\tau(X) \, \deg(H) + e_\tau(H)\,\deg(X) - (\tau_0+\dots+\tau_N)\,\deg(H)\,\deg(X)
\enspace,$$
avec ŽgalitŽ si et seulement si les variŽtŽs initiales de $X$ et $H$ s'intersectent proprement. 
\end{thm}

\begin{rem}
Lorsque $H$ est effectif, le terme correctif $\sum_{Y} m(X_\tau\cdot H_\tau;\iota(Y)) \, \deg(Y)$ est le degr\'e la partie du cycle intersection $X_\tau\cdot H_\tau$,
contenue dans la vari\'et\'e initiale de $X$.  
\end{rem}

Suivant l'attitude gŽnŽrale adoptŽe dans ce texte nous Žcrivons Žgalement au~\S~\ref{diviseurs} ce thŽorme en termes combinatoires pour l'intersection d'une variŽtŽ torique avec un diviseur monomial. 
Comme autre application de la formule de la hauteur d'une variété torique, on 
en dŽduit dans ce cas un thŽorme de BŽzout {\it exact} pour la hauteur normalisŽe du cycle intersection 
({\it voir} Corollaire~\ref{capre}). 

\medskip

L'étude des points algébriques de petite hauteur (ou {\it petits points}) 
a reçu une attention considérable au cours des dernières années. 
La taille des petits points dans une variété est quantifié par 
ses minimums algébriques successifs. 

Soit $X \subset \P^N(\Qbar)$  une vari{\'e}t{\'e} quasi-projective quelconque, 
de dimension $n$.
Pour $ \theta\ge 0$ on pose $X(\theta)$ pour l'ensemble des points de $X$ de hauteur normalisée $\hnorm$ majorée par $\theta$. 
Pour $i=1, \dots, n+1$, le {\it $i$-ème minimum algŽbrique de $X$ par rapport 
\`a la hauteur normalisŽe} est 
$$
\wh{\mu}_i(X) := \inf \big\{\theta : \dim (\ov{X(\theta)}) \ge n-i+1\big\}\enspace, 
$$
o\`u $\ov{X(\theta)}$ d\'esigne l'adh\'erence de Zariski de $X(\theta)$.
On {\'e}crit \,$
\wh{\mu}^\ess(X):= \wh{\mu}_{1}(X) $ \ 
et \ $ \wh{\mu}^\abs(X):= \wh{\mu}_{n+1} (X)
$\, 
pour les minimums {\em essentiel}  et {\em absolu} respectivement;
on a $\wh{\mu}_1(X)\geq\dots \geq \wh{\mu}_{n+1}(X) \newline \geq 0$.

La r{\'e}partition de la  hauteur des points alg{\'e}briques d'une variŽtŽ projective {\it fermŽe} $X$ est en relation  avec sa  hauteur, le lien est donn{\'e} par le {\it th{\'e}or{\`e}me des minimums successifs}~\cite[Thm. 5.2 et Lem. 6.5]{Zha95}: 
\begin{equation} \label{zhang} 
\wh{\mu}_1 (X)+ \cdots + \wh{\mu}_{n+1} (X) 
\, \le \,  
\frac{\wh{h}(X) }{\deg (X)}
\, \le \,  (n+1) \, \wh{\mu}_1(X)  
\enspace.\end{equation}

Comme application de la formule pour la hauteur d'une vari\'et\'e torique 
de \cite{PS04}, 
on construit au \S~6 des exemples montrant 
que toute configuration 
possible des minimums successifs 
est arbitrairement approchable et que le quotient ${\wh{h}(X) }/{\deg (X)}$
peut atteindre n'importe quelle valeur dans l'intervalle autorisé par 
l'encadrement ci-dessus:

\begin{thm} \label{densite}
Soient $n, N\in \N$ tels que $N \ge 3\,n+1$ 
et 
$\mu_1, \dots, \mu_{n+1}, \nu  \in \R$ tels que 
$$\mu_1\geq \dots\geq \mu_{n+1}\geq 0  \quad \mbox{ et } \quad  
\mu_1+\dots +\mu_{n+1} \le \nu < (n+1)\,\mu_1 \enspace.
$$ 
Alors pour $0<\varepsilon_1\leq (n+1)\mu_1-\nu$, 
$\varepsilon_2>0$ arbitraires, 
il existe 
une variŽtŽ torique $X\subset\P^N$ de dimension $n$
telle que 
$$
0 < \mu_i - \wh{\mu}_i(X) \le \varepsilon_1 \enspace \mbox{ pour } i=1,\dots,n+1 
\mbox{ et } \enspace
\left|\frac{\wh{h}(X)}{\deg(X)} - \nu\right| < \varepsilon_2\mu_1 \enspace.
$$
De plus, la variété $X$ peut être choisie de degré 
$\leq (4n^2 \, \varepsilon_2^{-1})^n$ et 
définie 
sur une extension kummerienne $K=\Q(2^{1/\ell})$ 
de degré $\leq\lfloor\log(2)\,\varepsilon_1^{-1}\rfloor +1$.
\end{thm}

Les exemples construits prŽsentent une codimension minimale $N-n=2n+1$ 
de l'ordre de la dimension de la variŽtŽ produite. 
La question se pose donc de savoir ce qu'il en est pour les 
variŽtŽs de petite codimension ({\it cf.} 
Proposition~\ref{cashyper}).

\bigskip
Pour faciliter la lecture, nous avons tressŽ dans le texte plusieurs paragraphes d'introduction aux variŽtŽs toriques (\S~\ref{vartor}), aux poids de Chow (\S~\ref{Chow}) et ˆ la thŽorie de l'intersection multi-projective (\S~\ref{vartor}). 
Les paragraphes~\S~\ref{obstruction} à~\S~\ref{minimums} doivent pouvoir tre lus de manire essentiellement indŽpendante.
Dans les paragraphes \ref{vartor} \`a~\ref{Chow} nous prenons un corps de base $\K$ alg\'ebriquement clos quelconque ou  
$\C$ pour le~\S~\ref{voltang}.  
L'arithm\'etique appara\^\i t \`a partir de la fin du \S~\ref{Chow}, o\`u nous posons nos conventions sur les places et valeurs absolues des corps de nombres (Convention~\ref{convention}).


\bigskip 

\setcounter{tocdepth}{2}

\vbox{\typeout{Matieres}

\tableofcontents}


\renewcommand{\thesection}{\Roman{section}}

%

\typeout{Generalites sur les varietes toriques projectives}

\section{G{\'e}n{\'e}ralit{\'e}s sur les variŽtŽs toriques projectives} \label{vartor}

\setcounter{equation}{0}
\renewcommand{\theequation}{\thesection.\arabic{equation}}

On note $\G_m^n:= (\K^\times)^n$ le tore alg{\'e}brique et  $\P^N$ l'espace 
projectif sur $\K$, de dimension $n$ et $N$ respectivement. Une vari{\'e}t{\'e} 
est toujours suppos{\'e}e r{\'e}duite et irr{\'e}ductible. 
Pour une famille de polyn{\^o}mes homog{\`e}nes $f_1,\dots,f_s \in 
\K[x_0,\dots,x_N]$ on pose $Z(f_1,\dots,f_s) \subset\P^N$ l'ensemble de ses zéros communs. Réciproquement, pour un ensemble algŽbrique $Z\subset\P^N$ on pose $I(Z)$ son idŽal de dŽfinition dans $\K[x_0,\dots,x_N]$.

On note $\R_+$ et $\R_+^\times$ les ensembles des nombres r{\'e}els non-n{\'e}gatifs et positifs, respectivement. On note encore $\N$ et $\N^\times$ les entiers naturels avec et sans 0, respectivement. Pour $N, D \in \N$ 
on pose $\N^{N+1}_D:= \Big\{a\in\N^{N+1} \, : \ a_0+\dots+a_N=D \Big\}$.

\bigskip

Dans ce paragraphe on donne un bref aper{\c c}u des propri{\'e}t{\'e}s géométriques 
des vari{\'e}t{\'e}s toriques. Pour plus de détails, on renvoie le lecteur 
à~\cite{GKZ94},~\cite{Ful93},~\cite{Ewa96},~\cite{Cox01}.

\smallskip

Soit $\cA= ( a_0,  \dots, a_N) \in (\Z^n)^{N+1}$ une suite de $N+1$ vecteurs de $\Z^n$, on consid{\`e}re l'action diagonale de $\G_m^n$ sur $\P^N$
$$
*_\geom : \G_m^n \times \P^N \to \P^N
\quad \quad , \quad \quad  
(s, x) \mapsto (s^{a_0} \, x_0 : \cdots : s^{a_N} \, x_N)
\enspace.
$$ 
On s'int{\'e}ressera {\`a} l'adh{\'e}rence de Zariski des orbites de cette action; pour un point $\alpha= (\alpha_0 : \cdots : \alpha_N) \in \P^N$ on pose
$$
X_\arith := \ov{\G_m^n *_\cA \alpha} \ \subset \P^N
$$
la {\it vari{\'e}t{\'e} torique projective} associ{\'e}e au couple $(\arith)$. 
Autrement-dit, $X_\arith$ est l'adh{\'e}rence de Zariski de l'image de
l'application monomiale
$$
\varphi_{\cA,\alpha} := *_\geom|_\alpha:  \G_m^n \to \P^N 
\enspace , \quad \quad  
s \mapsto (\alpha_{0} \, s^{a_0} : \cdots :   
\alpha_{N} \, s^{a_N} )
\enspace .
$$ 
C'est une vari{\'e}t{\'e} torique projective au sens de~\cite{GKZ94}, \cad
une sous-vari{\'e}t{\'e} de $\P^N$ stable sous l'action d'un tore
$\G_m^n$, avec une orbite dense $X_\arith^\circ:=\G_m^n *_\cA \alpha$. 

\smallskip 
Lorsque le point $\alpha$ est contenu dans un des sous-espaces standard de $\P^N$, la sous-vari{\'e}t{\'e} $X_\arith$ toute enti{\`e}re reste dans cet espace, puisque l'action $*_\cA$ est diagonale. Quitte {\`a} se restreindre au sous-espace standard
minimal contenant $\alpha$, on peut supposer \spdg $\alpha \in (\P^N)^\circ:=\P^N \setminus \{ x_0\cdots x_N =0\} $ et on fixe d{\'e}sormais pour $\alpha$ des 
coordonn{\'e}es  $(\alpha_0,\dots,\alpha_N) \in (\K^\times)^{N+1}$.

On notera $X_\geom$ la vari{\'e}t{\'e} torique associ{\'e}e {\`a}  $\geom \in
(\Z^n)^{N+1}$ et $( 1, \dots, 1)  \in (\K^\times)^{N+1}$. Dans ce cas, l'orbite principale 
$X_\geom^\circ$ est un sous-tore du tore $\G_m^N\cong(\P^N)^\circ$ et en fait 
tous les sous-groupes connexes de $(\P^N)^\circ$ sont de cette forme. 
Dans le cas g{\'e}n{\'e}ral
$$
X_\arith^\circ = \alpha \cdot X_\cA^\circ
$$
o{\`u} $\cdot$ d{\'e}signe la multiplication dans $(\P^N)^\circ$, \cad que 
l'orbite principale de $X_\arith$ est le translat{\'e} d'un sous-tore.

\medskip 

Le couple $(\arith)$ peut s'interpr{\'e}ter comme une suite finie de mon{\^o}mes 
$\alpha_{0}  \,  s^{a_0}, \dots, \alpha_{N} \, s^{a_N}$ de l'anneau des polyn{\^o}mes de Laurent $\K[s_1^{\pm 1} , \dots, s_n^{\pm 1}]$, on dit que les $a_i$ sont les {\em exposants} et les $\alpha_i$ les {\em coefficients}. Le {\it support} 
$$\Expo(\cA) := \{ a_{i_0}, \dots, a_{i_M} \}  \ \subset \Z^n$$ 
est l'ensemble des exposants {distincts} dans $\cA$. 
Les variétés $X_\arith$ et $X_{\expo(\cA)}$ sont lin{\'e}airement isomorphes par l'application
$$
\P^N \to \P^M  
\enspace, \quad \quad 
(x_0 : \cdots : x_N ) \mapsto
({\alpha_{i_0}^{-1}} \, x_{i_0} : \cdots : 
{\alpha_{i_M}^{-1}} \, x_{i_M} ) \enspace,
$$
et donc leurs propri{\'e}t{\'e}s {\it g{\'e}om{\'e}triques} sont les mêmes. 
Pour ces propriŽtŽs, on peut donc se ramener au cas habituel o{\`u} les
$a_i$ sont tous distincts et $\alpha_i =1$ pour tout $i$. Soit
$L_\geom \subset \Z^n $ le sous-module engendr{\'e} par les diff{\'e}rences
des vecteurs $a_0, \dots, a_N$, on a par exemple 
$$
\dim(X_\arith) = \rang_\Z (L_{\expo(\cA)})= \rang_\Z(L_\cA).$$ 
On introduit encore le polytope $Q_\cA \subset \R^n$ enveloppe convexe des vecteurs $a_0,\dots, a_N$.

\begin{lem} \label{truite} 
Soient $\cA =\{ a_0,\dots,a_N\}$ et $\cB=\{b_0,\dots,b_M\} \subset \Z^n$ des ensembles de cardinal $N+1$ et $M+1$ respectivement, avec $N\le M$. Posons $\pi:\P^M\to \P^N$ la projection lin{\'e}aire standard $(y_0:\dots:y_M) \mapsto (y_0:\dots:y_N)$, on a   

{\noindent (a)} \ 
si $\cA \subset \cB$ alors $X_\cA= \ov{\pi(X_\cB)}$; 
\smallskip 

{\noindent (b)} \ 
si $Q_\cA= Q_\cB$ alors $\pi$ est un morphisme rŽgulier et fini (au sens des fibres), de degr{\'e} $\deg(\pi)=[L_\cB:L_\cA]$. Si de plus  $L_\cA=L_\cB$ alors $\pi$ est un isomorphisme entre $X_\cB^\circ$ et $X_\cA^\circ$;
\smallskip 

{\noindent (c)} \ 
si $Q_\cA= Q_\cB=Q$ et si pour toute face $F$ de $Q$ on a $L_{\cA\cap F} = L_{\cB\cap F}$, alors $\pi:X_\cB\to X_\cA$ est un isomorphisme.  \end{lem}

\begin{demo} 
La partie~{\it (a)} est cons{\'e}quence directe des d{\'e}finitions. Pour les parties~{\it (b)} et~{\it (c)}, l'hypoth{\`e}se $Q_\cA=Q_\cB=Q$ entra{\^\i}ne que la projection $\pi$ est compatible avec la d{\'e}composition en orbites~(\ref{decomposition}) ci-dessous. Pour chaque face $F$ de $Q$, $\pi: X_{\cB,F}^\circ \to X_{\cA,F}^\circ$ est une application monomiale de degr{\'e} $[L_{\cB \cap F}: \pi^*(L_{\cA\cap F})] = [L_{\cB \cap F}: L_{\cA \cap F}]$, en particulier on voit que $\pi$ est r{\'e}guli{\`e}re {\`a} fibres finies (et de degrŽ $1$ entre $X_\cA^\circ$ et $X_\cB^\circ$ si $L_\cA=L_\cB$). 

Dans~{\it (c)}, l'hypoth{\`e}se $L_{\cB \cap F} =  L_{\cA \cap F}$ entra{\^\i}ne d{\'e}j{\`a} que $\pi$ est une bijection. Pour chaque exposant $b_i\in \cB$ consid{\'e}rons la face $F$ de $Q$ de dimension minimale contenant $b_i$. On {\'e}crit alors $b_i= \sum_{j: \, a_j \in \cA \cap F} \lambda_{i,j} \, a_j $ avec $\lambda_{i,j} \in \Z$, $\sum_j \lambda_{i,j} = 0$. L'inverse est alors $X_\cA \to X_\cB$, $x \mapsto (x^{\lambda_0} : \cdots  : x^{\lambda_M})$. 
\end{demo}

Suivant la philosophie g{\'e}n{\'e}rale, les 
propri{\'e}t{\'e}s g{\'e}om{\'e}triques des variétés toriques se tra\-dui\-sent en des {\'e}nonc{\'e}s combinatoires 
sur  les vecteurs $a_0, \dots, a_N \in \Z^n$ d{\'e}finissant l'action. En ce qui concerne la th{\'e}orie de l'intersection g{\'e}om{\'e}trique de ces vari{\'e}t{\'e}s, le r{\'e}sultat le plus fondamental est que le degr{\'e} s'identifie au 
volume du polytope $Q_\cA$. 

Le sous-module $L_\cA$ est un r{\'e}seau de l'espace lin{\'e}aire $L_\cA \otimes_\Z \R \subset \R^n$; on consid{\`e}re la forme volume $\mu_\cA$ sur 
cet espace lin{\'e}aire, invariante par translations et normalis{\'e}e de sorte que 
$$
\mu_\cA(L_\cA \otimes_\Z \R/ L_\cA )=1 
\enspace;
$$ 
autrement-dit, de sorte que le volume d'un domaine fondamental soit 1. 
Soit $r:=\rang_\Z(L_\cA)$, le degr{\'e} de $X_\arith $ s'explicite comme $r!$ fois 
le volume normalisé du polytope associ{\'e}~\cite[\S~I.6, Thm.~2.3]{GKZ94}, \cite[p.~111]{Ful93}
$$
\deg(X_\arith) = r! \mu_\cA(Q_\cA)
\enspace. 
$$
Soit maintenant $\eta : \Z^r \hookrightarrow \Z^n$ une application lin{\'e}aire
injective telle que $\eta(\Z^r) = L_\cA$ et posons $b_i:= \eta^{-1}(a_i) \in \Z^r$ puis $\cB:= (b_0, \dots, b_N) \in (\Z^r)^{N+1} $, alors $X_{\cB, \alpha} = X_\arith$~\cite[Ch.~5,~Prop.~1.2]{GKZ94}, et ainsi on peut toujours supposer \spdg 
$L_\cA = \Z^n$. Dans cette situation, la forme volume  $\mu_\cA$ 
co{\"\i}ncide avec la forme volume euclidienne $\Vol_{n}$ sur $\R^n$ et donc
\begin{equation}\label{degvol}
\dim(X_{\cA,\alpha}) = n
\enspace , \quad \quad
\deg(X_{\cA,\alpha}) = n!\, \Vol_{n}(Q_\cA)\enspace.
\end{equation}

\medskip

Plus g{\'e}n{\'e}ralement, les multidegr{\'e}s du tore $\G_m^n$ plongŽ dans un produit d'espaces projectifs {\it via} plusieurs applications monomiales s'expriment comme 
volumes mixtes des polytopes associ{\'e}s ˆ ces applications. Nous explicitons cela maintenant.

\smallskip 
 
Soit $Z\subset\P^{N_0}\times\dots\times\P^{N_m}$ une sous-variété de dimension $n$ et $c=(c_0,\dots,c_m)\in\N^{m+1}_n$ avec $0\le c_i \le N_i$, le {\it multidegr{\'e} de $Z$ d'indice $c$} est d{\'e}fini par
$$
\deg_{c}(Z) := \Card\Big( X \cap \pi_0^{-1}(E_0) \cap \cdots \cap
\pi_m^{-1}(E_m) \Big)
$$ 
o{\`u} $\pi_i$ d\'esigne la
projection $\P^{N_0}\times\dots\times\P^{N_m} \to \P^{N_i}$ et 
$E_i$ est un sous-espace linéaire g{\'e}n{\'e}rique de $\P^{N_i}$ de codimension $c_i$.

Pour $i=0, \dots, m$ on fixe
un groupe  $x_i= \{ x_{i , 0}, \dots, x_{i, N_i}\} $
de $N_i+1$
variables chacun 
et on consid{\`e}re
l'anneau $\K[x_0, \dots, x_m]$, multigradu{\'e} par
$\deg(x_{i, j}) := e_i$
o{\`u} $e_0, \dots, e_m \in \Z^{m+1}$
d{\'e}signent les vecteurs de la base standard. 
Ainsi $I(Z) \subset \K[x_0,\dots,x_m]$ est un idéal multi-homogène de 
rang $N_0+\cdots+N_m-n$. 
Soient
$ d_0, \dots, d_n \in \N^{m+1} $;
pour chaque $ d_i \in \N^{m+1}$
on introduit un groupe de variables $U_i= \{U_{i, 0},
\dots, U_{i, M_i}\}$, 
$M_i+1= \prod_{j=0}^m {d_j+N_j\choose N_j}$, 
et 
on consid{\`e}re la {\it forme r{\'e}sultante d'indice} $d_0, \dots, d_n$ de $I(Z)$: 
$$
\res_{d_0, \dots, d_n}(I(Z)) \in \K[U_0, \dots, U_n] \enspace.
$$ 
On renvoie le lecteur {\`a}~\cite[\S~3]{Rem01a} ou encore à~\cite[\S~I.2]{PS04} pour la d{\'e}finition et propri{\'e}t{\'e}s de base des formes r{\'e}sultantes. 
Maintenant, {\`a} un vecteur $c=(c_0,\dots, c_m)\in\N^{m+1}_n$ comme ci-dessus on associe l'indice partiel
$$
d(c):=(\underbrace{e_0,\dots,e_0}_{c_0\ \scriptstyle \rm fois}, \dots, \underbrace{e_m,\dots,e_m}_{c_m\ \scriptstyle \rm  fois})\in (\N^{m+1})^n 
\enspace,$$
o{\`u} $e_0,\dots,e_m$ d{\'e}signent les vecteurs de la base standard de $\R^{m+1}$. Pour toute choix de $d_0\in\N^{m+1}\setminus\{\mathbf 0\}$ on a~\cite[Prop.~3.4 et 2.11]{Rem01a}
$$\deg_c(Z) = \deg_{U_0}(\res_{d_0,d(c)}(I(Z)))\enspace.$$

Soit encore $D=(D_0,\dots,D_m)\in(\N^\times)^{m+1}$ et considérons
$$
\begin{array}{rll}
\Psi_D : & \P^{N_0} \times \cdots \times \P^{N_m} & \longrightarrow  
\P^{ {D_0+N_0 \choose N_0} \cdots {D_m+N_m \choose N_m} -1} \\[2mm]
&(x_0, \dots, x_m) & \longmapsto  \Big( x_0^{b_0} \cdots x_m^{b_m}
\, : \ b_0 \in \N^{N_0+1}_{D_0}, \dots, \ b_m \in \N^{N_m+1}_{D_m} \Big)
\end{array}
$$
le {\it plongement mixte} associŽ (composition des plongements de Veronese et Segre), on sait que~\cite[\S~2.3]{Rem01b}
\begin{equation}\label{corbeau}
\deg(\Psi_D(Z)) = \sum_{c\in\N^{m+1}_n}{n\choose c}\deg_c(Z)D^c
\enspace.\end{equation}

\smallskip

Consid{\'e}rons maintenant des ensembles convexes $Q_1, \dots, Q_n \subset \R^n$, leur {\it volume mixte} (ou {\it multi-volume}) est d{\'e}fini par la
formule de type inclusion-exclusion
\begin{equation*} 
\MV(Q_1, \dots, Q_n) := 
\sum_{j=1}^n (-1)^{n - j} 
\sum_{1 \le i_1 < \cdots < i_j \le n} 
\Vol_n(Q_{i_1} + \cdots + Q_{i_j}) \enspace. 
\end{equation*} 
Cette notion g{\'e}n{\'e}ralise le volume d'un ensemble convexe car 
$ \MV(Q, \dots, Q) = n! \, \Vol_n(Q)$. Le volume mixte est positif ou nul, sym{\'e}trique et lin{\'e}aire en chaque variable $Q_i$ par rapport {\`a} la somme de Minkowski. On renvoie {\`a}~\cite[\S~7.4]{CLO98} ou~\cite{Ewa96} pour ses propri{\'e}t{\'e}s de base. 

\smallskip

InterprŽtons ces notions dans le cas torique: soient $\cA_0\in(\Z^n)^{N_0+1}$,\dots, $\cA_m\in(\Z^n)^{N_m+1}$ tels que $L_{\cA_0}+\dots+L_{\cA_m}=\Z^n$. Posons $\underline{\cA}:=(\cA_0,\dots,\cA_m)$ et  considŽrons l'action associée $*_{\underline{\cA}}$ de $\G_m^n$ sur le produit d'espaces projectifs $\P^{N_0}\times\dots\times\P^{N_m}$:
$$\begin{array}{rccl}
*_{\underline{\cA}}: &\G_m^n\times\P^{N_0}\times\dots\times\P^{N_m} &\longrightarrow &\P^{N_0}\times\dots\times\P^{N_m}\\
&(s,x_0,\dots,x_m) &\longmapsto & (s*_{\cA_0} x_0,\dots, s*_{\cA_m} x_m)
\end{array}\enspace.$$ 
On considère la variété torique multi-projective $X_{\underline{\cA}} \subset \P^{N_0}\times\dots\times\P^{N_m}$, 
adhŽrence de Zariski de l'orbite du point $((1:\dots:1),\dots,(1:\dots:1))$, 
et plus gŽnŽralement on pose $X_{\underline{\cA},\underline{\alpha}}$ pour  l'adhérence de l'orbite de $\underline{\alpha}=(\alpha_0,\dots,\alpha_m)\in\P^{N_0}\times\dots \times\P^{N_m}$. 

\begin{prop}\label{multidegvol}
Soit $c\in\N^{m+1}_n$, dans la situation ci-dessus on a 
$$\deg_{c}(X_{\underline{\cA}}) =\MV_n(\underbrace{Q_{\cA_0},\dots, Q_{\cA_0}}_{c_0\ \scriptstyle \rm fois},\dots, \underbrace{Q_{\cA_m},\dots,Q_{\cA_m}}_{c_m\ \scriptstyle \rm fois}) \enspace.$$
\end{prop}

\begin{demo}
Pour $D=(D_0,\dots,D_m)\in(\N^\times)^{m+1}$ on pose $D\cdot\underline{\cA} := D_0\cA_0+\dots+D_m\cA_m$ l'ensemble des sommes de $D_i$ ŽlŽments pris dans $\cA_i$ pour $i=0,\dots,m$ et $Q_{D\cdot\underline{\cA}} = D_0Q_{\cA_0}+\dots+D_mQ_{\cA_m}$ l'enveloppe convexe de $D\cdot\underline{\cA}$. On vŽrifie $\Psi_D(X_{\underline{\cA}}) = X_{D\cdot\underline{\cA}}$, et par la multilinŽaritŽ du volume mixte
$$n!\Vol_n(Q_{D\cdot\underline{\cA}}) = \sum_{c\in\N^{m+1}_n} {n\choose c} \MV_n(\underbrace{Q_{\cA_0},\dots, Q_{\cA_0}}_{c_0\ \scriptstyle \rm fois},\dots, \underbrace{Q_{\cA_m},\dots,Q_{\cA_m}}_{c_m\ \scriptstyle \rm fois}) D^c\enspace.$$
On a $\deg(\Psi_D(X_{\underline{\cA}})) = n!\Vol_n(Q_{D\cdot\underline{\cA}})$ pour tout $D$ 
et la comparaison de l'identitŽ ci-dessus avec~(\ref{corbeau}) permet de conclure.
\end{demo}

\medskip 
Les orbites de l'action $*_\cA$ sur $X_\arith$ sont en correspondance avec l'ensemble $\Faces(Q_\cA)$  des faces du polytope $Q_\cA$. Pour chaque face $ P$ on associe un point $\alpha_P := (\alpha_{P, \, 0} : \cdots : \alpha_{P, \, N} ) \in \P^N$ d{\'e}fini par $  \alpha_{P, \, j}  :=\alpha_j$ si \,$ a_j \in P $  et  $  \alpha_{P, \, j}:= 0$ sinon; la bijection est  donn{\'e}e par~\cite[Ch.~5, Prop.~1.9]{GKZ94}, \cite[\S~3.1]{Ful93} 
$$
P \mapsto X_{\arith, P}^\circ := \G_m^n *_\cA \alpha_P \ \subset \P^N 
\enspace. 
$$
On a la d{\'e}composition 
\begin{equation} \label{decomposition}
X_\arith  = \bigsqcup_{P \in \faces(Q_\cA) } X_{\arith, P}^\circ 
\enspace. 
\end{equation} 
Posons  $N(P):= {\rm Card}\{ i \, : \ a_i \in P\} -1$ et
$$
\cA(P) = (a_i \, : \ a_i \in P) \in(\Z^n)^{N(P)+1}
\enspace,\quad \quad 
\alpha(P) := ( \alpha_i \, : \  a_i \in P) \in ({\K}^\times)^{N(P)+1} \enspace,
$$
on v{\'e}rifie que $X_{\arith, P}^\circ \subset \P^N$ 
est l'orbite principale d'une vari{\'e}t{\'e} torique 
contenue
dans un sous-espace standard de dimension ${N(P) }$. 
Restreinte {\`a} ce sous-espace, elle
s'identifie {\`a} la sous-vari{\'e}t{\'e} 
$ X_{\cA(P),\alpha(P)}^\circ \subset \P^{N(P) }$, 
de dimension {\'e}gale {\`a} la dimension (r{\'e}elle) de la face  $P$.

Plus g{\'e}n{\'e}ralement, soit  $X_{\un{\cA}} \subset \P^{N_0} \times \cdots
\times \P^{N_m}$ une vari{\'e}t{\'e} torique multi-projective  et posons 
$ Q_{\un{\cA}}:=(Q_{\cA_0},\dots,Q_{\cA_m})$ la famille de $m+1$
polytopes associ{\'e}e. 

Pour un polytope $Q \subset \R^n$ et un vecteur $b \in \R^n$ on note $Q(b)$ la face de $Q$ d\'efinie comme l'ensemble 
des points $u\in Q$  maximisant la fonctionnelle lin{\'e}aire $u \mapsto \langle b,u\rangle $. 
Consid{\'e}rons  les combinaisons des faces des $Q_{\cA_i}$ obtenues par cette m{\'e}thode:
$$
\Faces(Q_{\un{\cA}}) := \{ (Q_{\cA_0}(b),\dots,Q_{\cA_m}(b)) \, : \ b \in \R^n\} 
\subset \Faces(Q_{\cA_0}) \times \cdots
\times \Faces(Q_{\cA_m}) \enspace. 
$$ 
On peut vŽrifier que l'ensemble des sommes de Minkowski $Q_{\cA_0}(b) + \dots + Q_{\cA_m}(b)$ pour $b\in\R^n$, co{\"\i}ncide avec l'ensemble des faces de la somme de Minkowski $Q_{\cA_0} + \dots + Q_{\cA_m}$. 

Soit $\underline{\alpha}=(\alpha_0,\dots,\alpha_m)\in\P^{N_0} \times \dots \times \P^{N_m}$, les orbites de l'action $*_{\un \cA}$ sur $X_{\un \cA,\un\alpha}$ sont en correspondance avec l'ensemble $\Faces(Q_{\un{\cA}})$. Pour chaque $\un P=(P_0,\dots,P_m) \in \Faces(Q_{\un{\cA}})$ on consid{\`e}re le point 
$\alpha_{\un P} := (\alpha_{P_0},\dots,\alpha_{P_m}) \in \P^{N_0} \times \cdots
\times \P^{N_m}$; 
la bijection est  donn{\'e}e par
$$
\un P \mapsto X_{\un \cA, \un\alpha, \un P}^\circ := \G_m^n *_{\un \cA}
\alpha_{\un P} \ \subset \P^{N_0} \times \cdots
\times \P^{N_m}
\enspace, 
$$
et cette correspondance préserve la dimension, car $\dim(X_{\un \cA, \un\alpha, \un P}) = \dim(P_0+\dots+P_m)$. 
La d{\'e}monstration est une l\'eg{\`e}re variante de~\cite[Ch.~5, Prop.~1.9]{GKZ94}.

\bigskip 

Classiquement les vari{\'e}t{\'e}s toriques sont construites {\`a} partir d'{\it {\'e}ventails}, {\it voir } par exemple \cite{Ful93}. Rappelons qu'un Žventail est un ensemble de c™nes 
simpliciaux de sommet l'origine engendrŽs par un nombre fini de points de $\Z^n$, tels que toute face et toute intersection de c™nes de l'Žventail appartienne encore ˆ l'Žventail. 
Le lien avec notre présentation des variétés toriques {\it via} des polytopes r\'esulte de ce 
qu'un polytope $Q$  à sommets dans $\Z^n$
d{\'e}termine un Žventail et donc aussi une vari{\'e}t{\'e} torique au {\it sens classique}, qui de plus est munie d'un fibr{\'e} en droites ample. Dans les alinéa suivants (tir{\'e}s essentiellement de~\cite{Cox01} et \cite{Ful93}) 
on explicite cette construction.

\smallskip 

Supposons \spdg $\dim(Q)=n$ et posons $\Faces_i(Q)$ l'ensemble des faces de $Q$ de dimension $i$, de sorte que 
$\displaystyle \Faces(Q):=\sqcup_{i=0}^n\Faces_i(Q)$. 
Pour chaque hyperface (\cad face de codimension 1) $F$ de $Q$ on prends le vecteur normal int{\'e}rieur primitif $v_F \in \Z^n$ et on note 
$$m_F:=-\min\{\langle x,v_F\rangle:x\in Q\}\in\Z\enspace,$$
o\`u $\langle\cdot,\cdot \rangle$ d\'esigne le produit scalaire ordinaire sur $\R^n$.  
On peut donc d{\'e}crire $Q$ comme intersection de demi-espaces
\begin{equation} \label{betterave} 
Q= \bigcap_{F \in \faces_{n-1}(Q)} \{ x\in \R^n\, : \ \langle x, v_F\rangle\ge
-m_F\} \enspace.  
\end{equation} 
{\`A} chaque face $P$ de $Q$ on  associe le c{\^o}ne $\sigma_P$ engendrŽ par les vecteurs $v_F$ correspondant aux hyperfaces $F \supset P$; on a $\dim(P)+\dim(\sigma_P)=n$. L'ensemble de ces c{\^o}nes forme un {\'e}ventail {\it complet} $\Sigma$ de $\R^n$, \cad un Žventail recouvrant $\R^n$. 

La figure suivante montre la correspondance $P \mapsto \sigma_P$ entre faces de $Q$ et c™nes de  $\Sigma$, pour $Q$ le simplexe standard de $\R^2$. 


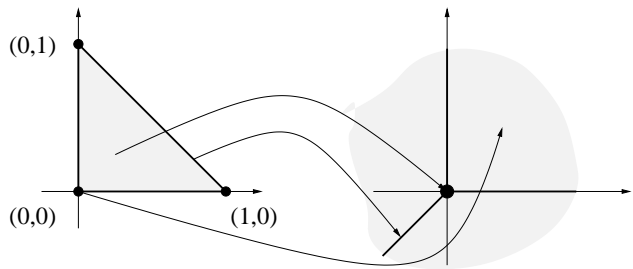
\bigbreak

\begin{figure}[htbp]

$$\epsfig{file=figures/eventail.eps, height= 35mm}$$

\end{figure}

\bigbreak


Notons $\cX_\Sigma$ la vari{\'e}t{\'e} torique (non plong{\'e}e dans un espace projectif) d{\'e}finie par cet {\'e}ventail. 
Pour chaque hyperface $F$ de $ Q$ on note $D_F$ le diviseur correspondant {\`a} la sous-vari{\'e}t{\'e} torique de $\cX_\Sigma$ associ{\'e}e ˆ l'arte $\sigma_F$ de 
$\Sigma$ et on introduit un diviseur de Weil sur $\cX_\Sigma$ 
invariant sous l'action de $\G_m^n$
$$
D_Q:= \sum_{F\in \faces_{n-1}(Q)} m_F \, D_F \enspace, 
$$
qui est en fait un diviseur de Cartier ample~\cite[\S~3.4]{Ful93}. 
Pour $a \in \Z^n$, l'application monomiale $\chi^a:\G_m^n \to \G_m, s \mapsto s^a$
peut se voir comme une fonction rationnelle sur $\cX_\Sigma$. Cette fonction rationnelle s'étend en une section globale de $\cO(D_Q)$ si et seulement si $a \in Q$, et en fait
$$
\Gamma(\cX_\Sigma, \cO(D_Q)) = \bigoplus_{a \in Q \cap \Z^n} \K\cdot \chi^a
\enspace. 
$$
En posant finalement $\cA_Q:=Q\cap\Z^n= \{ a_0,\dots,a_N\} $, on vŽrifie que l'application $\G_m^n \to \P^N, s \mapsto (s^{a_0}, \dots, s^{a_N})$ s'{\'e}tend en un morphisme rŽgulier $\cX_\Sigma \to \P^N$, dont l'image est la vari{\'e}t{\'e} torique projective  $X_{\cA_Q}$. De fait, ceci est le morphisme de normalisation de $X_{\cA_Q}$~\cite[Cor.~13.6]{Stu96}. 

\smallskip 
Cependant, notons que la notion d'Žventail permet de dŽfinir des variŽtŽs toriques qui, mme compltes, n'admettent pas nŽcessairement de diviseur $\G_m^n$-invariant ample. 
Toutefois, lorsqu'une variŽtŽ torique abstraite possde un tel diviseur, elle
admet une application monomiale régulière vers un espace projectif 
et son image est une variŽtŽ torique projective du type que nous Žtudions ici. 
On notera encore que les variŽtŽs toriques abstraites sont des variŽtŽs normales, {\it voir}~\cite[\S~2.1, p.~29]{Ful93}, ce qui n'est pas toujours le cas pour les variétés toriques projectives.

%

\typeout{Equations et indices d'obstruction successifs}

\section{{\'E}quations et indices d'obstruction successifs} 

\label{obstruction} 

Soit $X \subset \P^N$ une vari{\'e}t{\'e} projective de dimension $n\ge 0$
et $ I(X) \subset \K[x_0,\dots,x_N]$ son id\'eal de d\'efinition,
l'{\it indice d'obstruction} $\omega(X)$ est d{\'e}fini comme le plus petit degr{\'e} d'une {\'e}quation homog{\`e}ne $f \in I(X) \setminus\{0\}$. Plus g{\'e}n{\'e}ralement, pour $i=1,\dots,N-n$ on définit le {\em $i$-{\`e}me
indice d'obstruction de $X$} par 
$$
\omega_i(X):= \min\Big\{ D \in\N \, : \ \dim (Z(f \, : \ f \in I(X)_D)) \le N-i \Big\} 
\enspace,$$
o $I(X)_D$ dŽsigne la partie de degré $D$ de l'idŽal homogène $I(X)$. Alternativement, $\omega_i(X)$ est le plus petit entier $D \ge 0$ tel qu'il existe des polyn{\^o}mes homog{\`e}nes $f_1, \dots, f_{i} \in
I(X)$ de degr{\'e} born{\'e} par $D$ formant une intersection compl{\`e}te. 
{\'E}videmment $\omega(X) = \omega_1(X) $ et 
$$
1\le \omega_1(X) \le \cdots \le \omega_{N-n}(X)\enspace.
$$
L'invariant $\omega(X)$ joue un r{\^o}le important dans les probl{\`e}mes de Lehmer généralisé et de Bogomolov sur les tores~\cite{Dav02},~\cite{AD04},~\cite{Rat04}. D'un autre c™tŽ, les lemmes de z{\'e}ros consistent {\`a} minorer le premier indice $\omega(X)$ pour une vari{\'e}t{\'e} de dimension~0, {\'e}ventuellement munie de multiplicit{\'e}s; 
ces 
r{\'e}sultats sont des outils fondamentaux en th{\'e}orie des nombres transcendants, {\it voir~\cite{Ber87}}. Signalons que pour ces applications, il est souvent important de consid{\'e}rer des indices d'obstruction 
{\it sur des  sous-corps de $\Qbar$} ({\it i.e.} dont les {\'e}quations sont {\`a} coefficients dans un sous-corps 
fixé) et pour des {\it sch{\'e}mas projectifs} quelconques. 

Plus g{\'e}n{\'e}ralement encore, on peut consid{\'e}rer des {\it indices d'obstruction successifs 
relatifs {\`a} un ouvert} donné $U\subset \P^N$:
$$
\omega_i(X; U) := \min\Big\{ D \in\N\, : \ \dim (Z(f \, : \ f \in
I(X)_D) \cap U) \le N-i \Big\} \enspace. 
$$ 

\medskip 

Les idéaux des variétés toriques sont {binomiaux}, \cad engendr{\'e}s par des famille de polyn{\^o}mes de la forme 
$\alpha\, x^a - \beta \, x^b$ avec $a,b \in \N^n$ et $\alpha, \beta \in \K^\times$. 
Dans la suite on explicite la relation entre variétés toriques et idéaux binomiaux, qui est la clé de notre étude 
des indices d'obstruction successifs de ces variétés.

\smallskip

Soit $\cA =(a_0, \dots, a_N) \in (\Z^n)^{N+1} $ tel que $L_\cA=\Z^n$ et $\alpha= (\alpha_0, \dots, \alpha_N) \in (\K^\times)^{N+1}$. Consid{\'e}rons l'application lin{\'e}aire 
$$
\eta_\cA:\Z^{N+1} \to \Z\times\Z^{n} 
\enspace , \quad \quad  
\lambda \mapsto (\lambda_0+\cdots +\lambda_N,\lambda_0\, a_0 + \cdots
+\lambda_N \, a_N)
\enspace,$$ 
et posons $\Gamma_\cA:= \ker(\eta_\cA)$ son noyau. 
C'est un sous-module {\it satur{\'e}} de $\Z^{N+1}$ 
(\cad que le quotient $\Z^{N+1} /\Gamma_\cA$ est sans torsion) 
de rang $N-n$. 
Notons $\Delta\subset \R^{N+1}$ l'hyperplan d'Žquation $\lambda_0+\dots+\lambda_N=0$ et posons 
$\Delta^\Z:=\Delta\cap\Z^{N+1}$; \'evidemment on a 
$\Gamma_\cA\subset \Delta^\Z$.

Pour $b \in \R^{N+1}$ on {\'e}crit de fa{\c c}on unique $b= b_+-b_-$ avec 
$b_+, b_-\in \R_+^{N+1}$ {\`a} supports disjoints, \cad $(b_+)_i = b_i$ si
$b_i>0$ et 0 sinon, 
et  $(b_-)_i = -b_i$ si
$b_i<0$ et 0 sinon. 
Le r{\'e}sultat suivant est une reformulation de~\cite[Cor.~4.3]{Stu96}:

\begin{prop} \label{pasteque} 
$I(X_{\cA,\alpha})= \big( x^{b_+} -
\alpha^b \,x^{b_-} \, : \ b \in \Gamma_\cA\big)$. 
\end{prop} 

\begin{demo} 
Le  {c{\^o}ne} de $X_{\cA,\alpha}$ coïncide avec l'adh{\'e}\-rence de Zariski de l'image de
$$
\G_m^{n} \times \G_m \to \A^{N+1} 
\enspace , \quad \quad 
(s,t)= (s_1,\dots,s_n,t) \mapsto (\alpha_0\,t\,s^{a_0},\dots,
\alpha_N\,t\,s^{a_N})
\enspace,
$$
donc $I(X_{\cA,\alpha})$ co{\"\i}ncide avec le noyau de l'homomorphisme 
$$
\K[x_0, \dots, x_N] \to \K[s_1^{\pm1},\dots,s_n^{\pm1},t^{\pm1}]
\enspace,\quad\quad
x_i \mapsto \alpha_i\, t \,s^{a_i} \enspace; 
$$
le r{\'e}sultat devient une cons{\'e}quence directe
de~\cite[Cor.~4.3]{Stu96}. 
\end{demo}

Ainsi, l'id{\'e}al d'une vari{\'e}t{\'e} torique projective est binomial, 
il en r\'esulte automatiquement que c'est un idéal premier et homog{\`e}ne (puisque $X_\arith$ est une variété projective) 
ne contenant aucune des variables $x_i$, {\`a} cause de l'hypoth{\`e}se $\alpha \in (\P^N)^\circ$. 

\smallskip
Soit maintenant $\Gamma \subset \Z^{N+1}$ un sous-module quelconque et $\rho$ un {\it caract\`ere partiel}, \cad un 
homomorphisme $\rho:\Gamma\to \K^\times$. Ces données définissent un idéal binomial
$$
I(\Gamma,\rho) := \big( x^{b+}-\rho(b) x^{b_-} \, : \ b\in \Gamma\big) \subset \K[x_0,\dots,x_N]\enspace.
$$

\begin{prop}{(\cite[Cor.~2.6]{ES96})} --- \label{citrouille} 
La correspondence
$$
(\Gamma,\rho) \mapsto I(\Gamma,\rho)
$$
est une bijection entre les sous-modules satur{\'e}s $\Gamma \subset \Z^{N+1}$ munis d'un 
caract\`ere partiel $\rho$, et les idéaux de $\K[x_0, \dots,x_N]$ binomiaux, premiers  et  ne 
contenant aucune des variables $x_i$. 
De plus ${\rm rang}_\Z(\Gamma)=\rang(I(\Gamma,\rho))$.
\end{prop}

On vérifie sans peine que $I(\Gamma,\rho)$ est homogène si et seulement si $\Gamma\subset \Delta^\Z$. 
Ainsi, la donnée d'un couple $(\arith)$ définit un sous-module saturé $\Gamma_\cA \subset\Delta^\Z$ 
et un caract\`ere partiel
$\rho_\arith: b\mapsto \alpha^b$; la proposition~\ref{pasteque} peut être reformulée 
sous la forme 
$$
I(X_\arith)=I(\Gamma_\cA,\rho_\arith).
$$

R{\'e}ciproquement, à partir d'un idéal binomial de $\K[x_0, \dots,x_N]$ 
premier, homog{\`e}ne et ne contenant aucune des variables $x_i$, on peut construire $(\arith)$ tel que 
$I=I(X_\arith)$: soient  $\Gamma \subset \Delta^\Z$ le sous-module saturé et $\rho$ le caract\`ere partiel associés à $I$, $n:=N-{\rm rang}_\Z(\Gamma)$ et prenons $v_0, \dots, v_n \in \Z^{N+1}$ une base de l'orthogonal $\Gamma^\perp$ de
$\Gamma$ dans $\Z^{N+1}$ relativement au produit scalaire usuel. 
Puisque $\Gamma \subset \Delta^\Z$ on peut supposer \spdg $v_0=(1,\dots,1)$; on pose alors 
$$
a_i:= (v_{1,i}, \dots,v_{n,i}) \in \Z^{n} \quad \quad \mbox{ pour } i=0,\dots,N\enspace,
$$
et $\cA:= (a_0,\dots,a_N) \in (\Z^{n})^{N+1}$. 
En outre, $\rho$ peut s'étendre (de manière pas forcément unique) en un caract\`ere total $\rho:\Z^{N+1}\to \K^\times$ et 
on prend $\alpha:= (\rho(e_0),\dots, \rho(e_N))\in (\K^\times)^{N+1}$ où les $e_i$ désignent les vecteurs de la base standard de $\Z^{N+1}$.

Par construction $\Gamma_\cA=\Gamma$
et $L_\cA=\Z^n$, parce que $v_0, \dots,v_n$ est une {base} 
de $\Gamma^\perp$ et $\Gamma$ est saturŽ. La proposition~\ref{pasteque} entra{\^\i}ne $I(X_{\cA,\alpha})=I$. 
En particulier, on en déduit que la correspondance $X \mapsto I(X)$ est une bijection entre l'ensemble des vari{\'e}t{\'e}s 
toriques de $\P^N$ et l'ensemble des id{\'e}aux de $\K[x_0, \dots,x_N]$ 
binomiaux, premiers, homog{\`e}nes, ne contenant aucune des variables de $x_i$, $i=0,\dots,N$.
Dans la suite, on notera $\Gamma_X$ et $\rho_X$ le $\Z$-module et le caract\`ere partiel associés à une variété torique donnée  $X$. 
 
\begin{exmpl} 
Soit
$$
\G_m^2 \to \P^3 \enspace, \quad \quad (s,t) \mapsto
(1:s:3\, s^2\,t: s\,t^2) 
$$
et posons $S \subset \P^3$ la surface torique, adh{\'e}rence de Zariski de l'image de cette application. Soit
$$
[\eta_\cA]=\left[ 
\begin{array}{cccc} 
1&1&1&1\\[2mm]
0&1&2&1\\[2mm]
0&0&1&2
\end{array} 
\right]
$$
la matrice de l'application lin{\'e}aire $\eta_\cA:\Z^4\to\Z^3$ dans les bases canoniques, alors $\Gamma_\cA=\ker(\eta_\cA)$ est engendr{\'e} par le seul vecteur
$$
\gamma=(M_0,-M_1,M_2,-M_3) =(2,-3,2,-1) \in \Z^4
$$
o{\`u} $M_i$ d{\'e}signe le $i$-{\`e}me mineur de la matrice $[\eta_\cA]$. Et donc une {\'e}quation de $S$ est 
$$
x^{\gamma_+}-\alpha^\gamma\, x^{\gamma_-} = x_0^2\,x_2^2 -9\,
x_1^3\,x_3 \in \K[x_0,x_1,x_2,x_3]\enspace. 
$$
R{\'e}ciproquement, partons de l'{\'e}quation binomiale $f:=x_0^2\,x_2^2 -9\,x_1^3\,x_3 \in \K[x_0,x_1,x_2,x_3]$. 
Le $\Z$-module saturŽ $\Gamma \subset \Delta^\Z$  correspondant est engendr{\'e} par $\gamma= (2,-3,2,-1)$ et le caract\`ere partiel est
$\rho:\Gamma\to \K^\times, m\cdot \gamma \mapsto 9^m$. 

On v{\'e}rifie que $(1,1,1,1), (0,1,2,1), (0,0,1,2) \in \Z^4$ forment bien une base de $\Gamma^\perp$ et que 
$\Z^4\to \K^\times, b\mapsto (1,1,3,1)^b= 3^{b_2}$, est une extension possible de $\rho$, donc $Z(f) =X_{\cA,\alpha}$ avec $\cA=((0,0),(1,0),(2,1),(1,2)) \in (\Z^2)^4$ et $\alpha= (1,1,3,1)\in (\K^\times)^4$. 

La figure suivante montre les exposants et le polytope associ\'es \`a cette surface: 
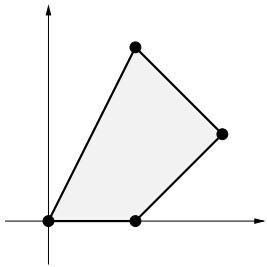
\vspace{-1mm} 

\begin{figure}[htbp]
$$\kern0.4cm\epsfig{file=figures/fano.eps, height= 35mm}$$
\end{figure}
\end{exmpl}

\vspace{10mm}
Considérons la norme $\ell^1$ 
$$
||v||_1= \sum_{i=0}^N |v_i| \quad \quad \mbox{ pour } v\in \R^{N+1} 
$$
et pour un $\Z$-module quelconque $\Gamma\subset \R^{N+1}$ notons $\mu_i := \mu_i(\Gamma; ||\cdot||_1)$ le $i$-{\`e}me minimum successif de $\Gamma$ relativement {\`a} cette norme. Rappelons que $\mu_i$ est  d{\'e}fini comme le plus petit 
$ \nu \in \R_+$ tel qu'il existe $i$ vecteurs ind{\'e}pendants $v_1,\dots, v_{i} \in \Gamma$ de norme born{\'e}e par $\nu$. Alternativement 
$$
\mu_i(\Gamma; ||\cdot||_1):= \min \{  \nu\in\R_+ \, : \ 
{\rm rang}_\Z(\nu\cdot  B_{||\cdot||_1}  \cap \Gamma) \ge i\} \enspace,
$$
o{\`u} $B_{||\cdot||_1}$ d{\'e}signe la boule unit{\'e} pour la norme $\ell^1$. 
Le r{\'e}sultat suivant montre que les indices d'obstruction d'une vari{\'e}t{\'e} torique relatifs 
ˆ l'ouvert $(\P^N)^\circ$ coïncident avec les minimums successifs du $\Z$-module associ\'e:  

\begin{prop} \label{framboise}  
Soit $X\subset \P^N$ une variété torique de dimension $n$, 
alors 
$$ \omega_i(X; (\P^N)^\circ) = \frac12\,
\mu_i(\Gamma_X; ||\cdot||_1) \quad \mbox{ pour } i=1,\dots,N-n \enspace.
$$ 
\end{prop}

\begin{demo}
On peut se restreindre \spdg au cas $\alpha = (1,\dots,1)$. Soient $v_1,\dots, v_{N-n}$ des vecteurs formant une base de $\Gamma_X$. D'aprs~\cite[Thm.~2.1(b)]{ES96}, $X^\circ = X\cap (\P^N)^\circ$ est une intersection complète 
decoupé par les bin™mes 
$$
x^{(v_1)_+} - x^{(v_1)_-} \enspace,  \dots, \enspace 
x^{(v_{N-n})_+} - x^{(v_{N-n})_-} 
\in \K[x_0,\dots,x_N]
$$
et on a $\displaystyle \deg(x^{(v_j)_+} - x^{(v_j)_-}) = \frac12\,||v_j||_1$ (car $\sum_{k=0}^Nv_{j,k}=0$). En prenant $v_1,\dots,v_{N-n}$ r{\'e}alisant les minimums successifs du module  $\Gamma_X$ on obtient 
$$
\omega_i(X; (\P^N)^\circ) \le \max\{ \deg(x^{(v_j)_+}-x^{(v_j)_-})\,: \ 1\le j \le i\} = \frac12\, \mu_i(\Gamma_X; ||\cdot||_1) \enspace.  
$$

Dans l'autre direction, soit $f_1, \dots, f_{N-n} \in I(X)$ 
une suite de polyn{\^o}mes homog{\`e}nes  r{\'e}alisant les indices d'obstruction successifs de $X$ sur $(\P^N)^\circ$. Par la proposition~\ref{pasteque} 
$$
f_j(x) = \sum_{b \in \Gamma_X} g_{j,b}(x) \, (x^{b_+}-x^{b_-}) 
$$
pour certains $g_{j,b} \in \K[x_0,\dots,x_N]$ et, pour tout $b\in\Gamma_\cA$ tel que $g_{j,b}\not=0$, on a 
\begin{equation}\label{minodegfj}
\deg(f_j) = \deg(g_{j,b})+ \deg(x^{b_+}-x^{b_-}) 
= \deg(g_{j,b})+ \frac12 \, ||b ||_1 \ge \frac12\, ||b||_1
\enspace.\end{equation}
Soit $1\le i \le N-n$, posons $L_i$ le $\Z$-module engendr{\'e} par les $b \in \Gamma_X$ tels que $g_{j,b}\ne0$ pour un certain $1\le j \le i$. Alors, $f_j \in J_i:=(x^{b_+}-x^{b_-} \, : \ b \in L_i) \subset \K[x_0,\dots,x_N] $ pour $j=1,\dots,i$ et donc 
$$
i = {\rm rang}(f_1,\dots,f_i) \le {\rm rang}(J_i) = {\rm rang}_\Z(L_i) 
$$
par la proposition~\ref{citrouille} ci-dessus. Ainsi $L_i$ est un sous-module de $\Gamma_X$ de rang au moins $i$ et on peut trouver $i$ vecteurs linŽairement indŽpendants parmi les $b\in\Gamma_\cA$ tels que $g_{j,b}\not=0$ pour un certain $1\le j \le i$. L'un au moins de ces vecteurs est de norme $\Vert b\Vert_1\geq \mu_i(\Gamma_X;\Vert\cdot\Vert_1)$ et on en dŽduit avec~(\ref{minodegfj})
$$
\omega_i(X; (\P^N)^\circ) = \max\{ \deg(f_j) \, : \ j=1,\dots,i\} 
\ge \frac12\, \mu_i(\Gamma_X, ||\cdot||_1)\enspace. 
$$
\end{demo} 

Soit $X\subset \P^N$ une vari{\'e}t{\'e} de dimension $n$, comme cons{\'e}quence directe de sa majoration de la fonction de Hilbert g{\'e}om{\'e}trique, Chardin a montr{\'e} qu'on a toujours 
$$
\omega_1(X) \le  N\,\deg(X)^{\frac{1}{N-n}} \enspace,
$$
et aussi $\omega_2(X)^\frac12 \le  N\,\deg(X)^{\frac{1}{N-n}}$~\cite{Cha89}. 
Ult{\'e}rieurement, Chardin et Philippon ont consi\-d{\'e}\-r{\'e} le probl{\`e}me de l'{\it interpolation alg{\'e}brique}, qui consiste \`a estimer le degr{\'e} minimal de sous-vari{\'e}t{\'e}s 
$$ 
X = Y_{N-n} \subset \cdots \subset Y_1 \subset \P^N
$$
telles que $\codim(Y_j)=j$. Ils ont montrŽ qu'on peut choisir $Y_1,\dots, Y_{N-n}$ comme ci-dessus, satisfaisant~\cite{CP99}
\begin{equation}\label{intCP}
\deg(Y_j)^\frac{1}{j} \le  N \, 4^{N-1} \, \deg(X)^{\frac{1}{N-n}}
\enspace. 
\end{equation}
L'interpolation alg{\'e}brique est tr{\`e}s proche des  probl{\`e}mes d'estimation 
des indices d'obstruction successifs. En fait, une reformulation de leur d{\'e}monstration montre qu'il existe des {\'e}quations $f_1,\dots,f_{N-n} \in I(X)$ formant une intersection compl{\`e}te au voisinage du point g{\'e}n{\'e}rique de $X$ et telles que 
$$
\deg(X) \le \deg(f_1)\cdots \deg(f_{N-n}) \le c(N) \, \deg(X)
$$
pour une constante $c(N) >0$ explicite. Ceci \'equivaut aux in\'egalit\'es
\begin{equation} \label{patate}  
\deg(X) \leq \omega_1(X;U) \cdots \omega_{N-n}(X;U) \le c(N) \, \deg(X)
\end{equation} 
pour un ouvert $U \subset \P^N$ tel que $X \cap U \ne \emptyset$. 
L'interpr{\'e}tation des indices d'obstruction des vari{\'e}t{\'e}s toriques comme 
minimums successifs d'un module permet de traduire le deuxi{\`e}me th{\'e}or{\`e}me de Minkowski en des estimations pour le produit de ces indices relatifs ˆ l'ouvert $(\P^N)^\circ$:

\begin{cor} Soit $X\subset \P^N$ une variété torique de dimension $n$, alors
$$
\deg(X) \le \omega_1(X ; (\P^N)^\circ)
\cdots  
\omega_{N-n}(X; (\P^N)^\circ)
\le c(N,n) \, \deg(X)
$$
avec 
$$
c(N,n):={N+1\choose n+1}^\frac12 \, \left(\frac{N+1}{\pi}\right)^{\frac{N-n}{2}} \Gamma\left(1+\frac{N-n}{2}\right) \le \left(\frac{N+1}{\sqrt{\pi}}\right)^{N-n} \enspace.$$
\end{cor}

Ceci est la proposition~\ref{intro2} dans l'introduction, avec une constante plus précise. 

\begin{demo}
La premi{\`e}re in{\'e}galit{\'e} est cons{\'e}quence directe du th{\'e}or{\`e}me de B{\'e}zout; passons  
{\`a} la deuxi{\`e}me. Posons $\Gamma_X^\R:= \Gamma_X \otimes \R \cong\R^{N-n}$, de sorte que $\Gamma_X$ est un {r{\'e}seau} de $\Gamma_X^\R$. Posons aussi $B:=  B_{||\cdot||_1} \cap \Gamma_X^\R$ la boule unit{\'e} de $\Gamma_X^\R$ par rapport {\`a} la restriction de la norme $\ell^1$. Ceci {\'e}tant un ensemble convexe et sym{\'e}trique par rapport
\`a l'origine, le deuxi{\`e}me th{\'e}or{\`e}me de Minkowski sur les minimums successifs d'un r{\'e}seau entra{\^\i}ne 
$$
\frac{2^{N-n}}{(N-n)!} \le 
\frac{\Vol_{N-n}(B)}{\Vol_{N-n}(\Gamma_X^\R/\Gamma_X)} \,  \prod_{i=1}^{N-n} \mu_i(\Gamma_X;\Vert\cdot\Vert_1)
\le 2^{N-n} 
$$
et donc
\begin{equation}\label{jolieformule}
\prod_{i=1}^{N-n} \omega_i(X; (\P^N)^\circ) \le 
\frac{\Vol_{N-n}(\Gamma_X^\R/\Gamma_X)}{\Vol_{N-n}(B)} 
\end{equation}
par la proposition~\ref{framboise}. 

Pour conclure, il suffit de majorer le quotient
${\Vol_{N-n}(\Gamma_X^\R/\Gamma_X)}/{\Vol_{N-n}(B)}$. Soient $\cA\in (\Z^n)^{N+1}$ et $\alpha\in (\K^\times)^{N+1}$
tels que $X=X_\arith$; en particulier $\Gamma_\cA=\Gamma_X$. 
Soient 
$$
v_i:= (a_{0,i},\dots,a_{N,i}) \in \Z^{N+1} \quad \mbox{ pour } i=0,\dots,n
$$
les lignes de la matrice (dans les bases standard) de l'application $\eta_\cA:\Z^{N+1}\to\Z^{n+1}$.
Le module $\Gamma_X$ est l'orthogonal du sous-module $V \subset \Z^{N+1}$ engendr{\'e} par $v_0, \dots, v_n$. 
C'est aussi un sous-module satur{\'e} {\`a} cause de l'hypoth{\`e}se $L_\cA=\Z^n$; la formule de Brill-Gordan~({\it voir} par exemple~\cite[formule~(II.2)]{PS04}) entra{\^\i}ne alors
$$
\Vol_{N-n}(\Gamma_X^\R/\Gamma_X) = \Vol_{n+1}(V^\R /V) 
$$
et par la formule de Cauchy-Binet
$$
\Vol_{N-n}(\Gamma_X^\R/\Gamma_X) = \Vol_{n+1}(V^\R /V)  = ||v_0\wedge \cdots\wedge v_n||_2
= \bigg(\sum_{J: \card(J) = n+1}  \det([\eta_\cA]_J)^2 \bigg)^{1/2} 
\enspace,
$$
o\`u $[\eta_\cA]_J$ d\'esigne le mineur $(n+1)\times (n+1)$ de la matrice de 
$\eta_\cA$ dont les colonnes sont index\'ees par $J$.  
Chaque terme $\det([\eta_\cA]_J)$ est $\pm n!$ fois le volume d'un simplexe
contenu dans le polytope $Q_\cA$; ainsi on peut majorer cette quantit{\'e} par
$$
\Vol_{N-n}(\Gamma_X^\R/\Gamma_X) \le {N+1\choose n+1}^{1/2} \, n!\, \Vol_n(Q_\cA) = {N+1\choose n+1}^{1/2} \, \deg(X) 
\enspace.$$
En outre, notons que $B$ contient la boule euclidienne de $\Gamma^\R_X$ de rayon $(N+1)^{-1/2}$ 
centrŽe en l'origine, ainsi
$$
\Vol_{N-n}(B) \ge \Vol_{N-n} ((N+1)^{-1/2}B_{||\cdot||_2} \cap \Gamma_X^\R) = 
\left(\frac{\pi}{N+1}\right)^{\frac{N-n}{2}} \cdot\frac{1}{\Gamma\left(1+\frac{N-n}{2}\right)} \enspace.
$$ 
On en dŽduit
\begin{eqnarray*} 
\prod_{i=1}^{N-n} \omega_i(X;(\P^N)^\circ)
&\le &\frac{\Vol_{N-n}(\Gamma_X^\R/\Gamma_X)}{\Vol_{N-n}(B)} \\[0mm]
&\le &{N+1\choose n+1}^\frac12 \cdot\left(\frac{N+1}{\pi}\right)^{\frac{N-n}{2}} \cdot {\Gamma\left(1+\frac{N-n}{2}\right)} \cdot  \deg(X)\\[2mm] 
&\le &\frac{(N+1)^{\frac{N-n}{2}}}{(N-n)!^{\frac{1}{2}}}\cdot \left(\frac{N+1}{\pi}\right)^{\frac{N-n}{2}} \cdot \prod_{k=0}^{\left[\frac{N-n}{2}\right]-1} \left(\frac{N-n}{2}-k\right)\cdot \deg(X)\\[0mm]
&\le &\left(\frac{N+1}{\sqrt{\pi}}\right)^{N-n} \cdot \deg(X) \enspace.
\end{eqnarray*}
\end{demo}

Ce r{\'e}sultat am{\'e}liore la constante $c(N)$ dans~(\ref{patate}) pour le
cas torique et pr{\'e}cise l'ouvert $U$ par rapport auquel les indices peuvent tre considŽrŽs. 

\medskip
En revenant au cas général d'une {vari{\'e}t{\'e}} $X \subset \P^N$ quelconque de dimension $n$, il est naturel de se demander si l'estimation~(\ref{patate}) 
reste valide pour $U=\P^N$. Autrement-dit, s'il existe toujours des polyn{\^o}mes homog{\`e}nes $f_1,\dots, f_{N-n} \in I(X)$ formant une
intersection compl{\`e}te {\it globale} et tels que 
$$
\deg(X) \le \deg(f_1) \cdots \deg(f_{N-n}) \le c(N) \,  \deg(X)\enspace. 
$$
Il serait int{\'e}ressant de d{\'e}cider cette question dŽjˆ sur la classe des vari{\'e}t{\'e}s toriques. Pour ce faire, on voudrait expliciter les indices $\omega_i(X;\P^N)$ de manière analogue à  la proposition~\ref{framboise} 
pour $\omega_i(X;(\P^N)^\circ)$.

%

\typeout{Volumes et hauteurs des espaces tangents}

\section{Volumes, hauteurs d'espaces tangents et degrŽs}\label{voltang}

Dans la dŽfinition du~\S~\ref{vartor}, les sous-groupes algŽbriques connexes de ${\bf G}_m^N$ correspondent aux variŽtŽs toriques $X_\cA^\circ$ {\it via} l'identification $\iota:\G_m^N\hookrightarrow(\P^N)^\circ$. D'un autre c™tŽ, l'ensemble des points complexes d'un sous-groupe algŽbrique connexe $G$ de ${\bf G}_m^{N}$, de dimension $n$, est dŽcrit {\it via} l'application exponentielle par son espace tangent $TG(\C)\subset\C^N$ ˆ l'origine et son rŽseau de pŽriodes 
$$\Lambda:=TG({\bf C})\cap (2{\rm i}\pi{\bf Z})^{N}
\enspace,$$ 
de rang $n$ sur ${\bf Z}$. 

Munissons l'espace tangent ˆ ${\bf G}_m^N$ en l'origine $T{\bf G}_m^N({\bf C}) \simeq {\bf C}^N$ de la structure euclidienne pour laquelle les $2N$ ŽlŽments $e_i$ et $2{\rm i}\pi e_i$ ($i=1,\dots,N$) forment une base orthonormŽe. 

Soit $\cA:=(a_0,\dots,a_N)\in(\Z^n)^{N+1}$ tel que $L_\cA=\Z^n$ et $G:=\iota^{-1}(X_{\cA}^\circ)$. Des gŽnŽrateurs de $\Lambda$ sont donnŽs par les vecteurs $2{\rm i}\pi\lambda_1,\dots, 2{\rm i}\pi\lambda_n\in(2{\rm i}\pi{\bf Z})^N$ 
o $\lambda_{i,j}=a_{j,i}-a_{0,i}$ pour $i=1,\dots,n$ et $j=1,\dots,N$.  
En fait, le $\Z$-module $\frac{1}{2{\rm i}\pi}\Lambda$ est l'orthogonal de 
la projection par $\pi:\R^{N+1}\to \R^N, (x_0,x_1,\dots,x_N)\mapsto (x_1,\dots,x_N)$, 
du $\Z$-module $\Gamma_\cA \subset \Delta^\Z$, introduit au~\S~\ref{obstruction} prŽcŽdent. 

\medskip

ConsidŽrons encore l'adhŽrence de Zariski $\overline G$ de $G$ dans la compactification $({\bf P}^1)^N$ de ${\bf G}_m^{N}$, ses multidegrŽs sont indexŽs par des multi-indices $c=(c_1,\dots,c_N)\in\N_n^N$ o $c_j\in\{0,1\}$. {\`A} toute suite croissante $1\leq j_1<\dots< j_n\leq N$ on associe un multi-indice $c(j_1,\dots,j_n)$ dŽfini par $c(j_1,\dots,j_n)_j=1$ si $j\in\{j_1,\dots,j_n\}$ et $0$ sinon. Indiquons comment on retrouve par ce biais la proposition~4 de~\cite{BP}. 

\begin{prop}\label{voltangent}
Avec les notation ci-dessus, le volume euclidien de la projection de $\Lambda$ sur le produit des $n$ facteurs de ${\bf R}^N$ indexŽs par $j_1,\dots,j_n$ est Žgal ˆ $\MV_n(\overline{a_0\,a_{j_1}}, \dots, \overline{a_0\,a_{j_n}}) = \deg_{c(j_1,\dots,j_n)}(\overline G)$, o $\ov{a_0\,a_j}$ dŽsigne le segment de droite joignant $a_0$ ˆ $a_j$ dans $\R^n$. 
\end{prop}
\begin{demo}
En effet, la projection de $\Lambda$ considŽrŽe est un rŽseau de $(2{\rm i}\pi{\bf R})^n$ engendrŽ par les vecteurs $(2{\rm i}\pi\lambda_{i,j_1}, \dots, 2{\rm i}\pi\lambda_{i,j_n})$ pour $i=1,\dots,n$. Son volume pour la structure euclidienne fixŽe est donc Žgal ˆ $\left|\Lambda_{j_1, \dots, j_n}\right| = \MV_n(\ov{a_0\,a_{j_1}}, \dots, \ov{a_0\,a_{j_n}})$ o 
$$\Lambda_{j_1,\dots,j_n} = \det\left[\begin{array}{ccc}
\lambda_{1,j_1} &\dots &\lambda_{1,j_n}\\
\vdots & &\vdots\\
\lambda_{n,j_1} &\dots &\lambda_{n,j_n}
\end{array}\right] \in {\bf Z}
\enspace.$$ 
Finalement, par la formule~(\ref{multidegvol}) ce multivolume est Žgal
au multidegrŽ $\deg_{c(j_1,\dots,j_n)}(\overline G)$ correspondant.
\end{demo}

On notera que les $\Lambda_{j_1,\dots,j_n}$ introduits dans la dŽmonstration prŽcŽdente sont des coordonnŽes grasmaniennes de l'espace tangent $TG({\bf C})$ en l'origine de $G$ dans ${\bf C}^N$. Comme $\Lambda$ est un rŽseau primitif de $TG({\bf C})$ on a 
$${\rm ppcm}(\Lambda_{j_1,\dots,j_n} \, : \ 1\leq j_1<\dots<j_n\leq N)=1
\enspace.$$
\smallskip

L'image de $G$ par le plongement de Segre $s:({\bf P}^1)^N\rightarrow{\bf P}^{2^N-1}$ est dŽcrite par la somme de Minkowski 
$$\overline{a_0\,a_1} + \dots + \ov{a_0\,a_N}\subset{\bf R}^n
\enspace,$$ 
dont le volume est Žgal ˆ la somme des volumes mixtes pour toutes les projections du type envisagŽ dans la proposition~\ref{voltangent}. En particulier
$$\begin{array}{rcl}
\deg_{{\bf P}_{2^N-1}}(\overline{s(G)}) &= &n!{\rm Vol}_n(\ov{a_0\,a_1} + \dots + \ov{a_0\,a_N})\\[2mm] 
&= &n!\sum_{1 \leq j_1 < \dots < j_n\leq N} \MV_n(\ov{a_0\,a_{j_1}}, \dots, \ov{a_0\,a_{j_n}})\\[2mm]
&= &n!\sum_{1 \leq j_1 < \dots < j_n\leq N} |\Lambda_{j_1,\dots,j_n}|
\enspace.\end{array}$$
Ainsi, le degrŽ de $\overline{s(G)}$ est Žgal ˆ $n!$ fois la hauteur associée à la norme $\ell^1$ des coordonnŽes grasmaniennes de son espace tangent dans ${\bf C}^N$.

Si l'on considère la hauteur de Schmidt $h_S(T)$ d'un sous-espace $T\subset\Qbar^N$ dŽfinie comme la hauteur projective du point qui le reprŽsente dans la variŽtŽ grasmanienne~\cite[page~28]{Sch91}, on a la formule
$$h_S(TG) = \left(\sum_{1\leq j_1<\dots<j_n\leq N}|\Lambda_{j_1,\dots,j_n}|^2\right)^{\frac{1}{2}} = {\rm Vol}_n(TG({\bf C})\cap(2{\rm i}\pi{\bf R})^N/\Lambda) = {\rm Vol}_{N-n}(\Gamma_\cA^\R/\Gamma_\cA)
\enspace,$$
car $TG({\bf C})\subset{\bf C}^N$ est dŽfini sur ${\bf Q}$ et ses coordonnŽes grasmaniennes $\Lambda_{j_1,\dots,j_n}$ dans ${\bf Z}$ sont premires entre elles dans leur ensemble. La seconde ŽgalitŽ n'est autre que la formule de Cauchy-Binet dŽjˆ utilisŽe au~\S~\ref{obstruction}. 
La troisime ŽgalitŽ s'obtient en identifiant l'hyperplan 
$\Delta = \{\lambda_0+\dots+\lambda_N=0\} \subset \R^{N+1}$ avec $\R^N$ par la projection sur les $N$ dernires coordonnŽes, de sorte que l'image de $\Delta^\Z$ est $\Z^N$.

En remarquant que $Q_\cA={\rm Conv}(a_0,\dots,a_N)$ 
contient tous les simplexes de sommets pris parmi les $a_i$ et est contenu dans l'union de ces mmes simplexes ayant $a_0$ comme sommet fixe, on vŽrifie facilement (le volume d'un simplexe ${\rm Conv}(a_0,a_{i_1},\dots,a_{i_n})$ 
est Žgal ˆ $\frac{1}{n!}|\Lambda_{i_1,\dots ,i_n}|$)
$$\frac{1}{n!}\max_{1\leq j_1,\dots,j_n\leq N}|\Lambda_{j_1,\dots,j_n}| \leq {\rm Vol}_n(Q_\cA) \leq \frac{1}{n!}\sum_{1 \leq j_1 < \dots < j_n\leq N} |\Lambda_{j_1,\dots,j_n}|
\enspace.$$
D'aprs la formule~(\ref{degvol}), on sait que $\deg_{{\bf P}^N}(\overline G) = n!{\rm Vol}_n(Q_\cA)$ et 
en particulier

\begin{prop}\label{compdegseg}
Pour tout sous-groupe algŽbrique $G\subset \G_m^N$ de dimension $n$ on a
$$n!\deg_{{\bf P}^N}(\overline{G}) \leq \deg_{{\bf P}^{2^N-1}}(\overline{s(G)}) \leq n!\binom{N}{n} \deg_{{\bf P}^N}(\overline{G})
\enspace.$$
\end{prop}
On pourra comparer avec le rŽsultat analogue pour la hauteur normalisŽe dans~\cite[Prop.~2.2]{DP99}.

\smallskip

Si l'on introduit la fonction $f_{\cal A}$ sur l'enveloppe convexe $Q_{\cal A}$ des points ${\cal A}=(a_0,\dots,a_N)$ dans ${\bf R}^n$, ˆ valeur dans ${\bf N}$, dŽfinie pour $u\in Q_{\cal A}$ par~:
$$f_{\cal A}(u) := {\rm Card}\left\{(j_1,\dots,j_n);1\leq j_1<\dots<j_n\leq N\enspace {\rm et}\enspace u\in{\rm Conv}(a_0,a_{j_1},\dots,a_{j_n})\right\}
\enspace,$$
on vŽrifie facilement l'ŽgalitŽ
$$\frac{\deg_{{\bf P}^{2^N-1}}(\overline{s(G)})}{n!\deg_{{\bf P}^N}(\overline G)} = \frac{1}{{\rm Vol}_n(Q_{\cal A})}\cdot\int_{Q_{\cal A}}f_{\cal A}(u)du
\enspace.$$

\medskip

\begin{exmpl}
En considŽrant des points $a_0,\dots,a_N$ dans ${\bf Z}^n$ tels que $a_{n+1},\dots,a_N$ soient trs proches de $a_0$ comparativement ˆ $a_1,\dots,a_n$ ({\it voir} figure ci-dessous pour $n=2$), on a
$$\sum_{1 \leq j_1 < \dots < j_n\leq N} |\Lambda_{j_1,\dots,j_n}| \approx \max_{1\leq j_1,\dots,j_n\leq N}|\Lambda_{j_1,\dots,j_n}|=|\Lambda_{1,\dots,n}|
\enspace.$$
On obtient ainsi une variŽtŽ torique pour laquelle $n!\deg_{{\bf P}_N}(\overline G) \approx \deg_{{\bf P}_{2^N-1}}(\overline{s(G)})$, montrant que l'inŽgalitŽ de gauche dans la proposition~\ref{compdegseg} est optimale.
\end{exmpl}



\begin{exmpl}
ConsidŽrons maintenant $N=mn$ points $a_1,\dots,a_N$ de ${\bf Z}^n$ Žgalement con\-cen\-trŽs autour de $a_1,\dots,a_n$, par exemple $a_{kn+i}\approx a_i$ ˆ l'intŽrieur de ${\rm Conv}(a_0,a_1,\dots,a_n)$, pour $k=1,\dots,m-1$ et $i=1,\dots,n$ ({\it voir} figure ci-dessous pour $n=2$). On a alors
$$\sum_{1 \leq j_1 < \dots < j_n\leq N} |\Lambda_{j_1,\dots,j_n}|\  \rlap{$\geq$}\raise-7pt\hbox{$\sim$}\ \left(\frac{N}{n}\right)^n n!{\rm Vol}_n(a_0,a_1,\dots,a_n)\enspace.$$
On obtient ainsi une variŽtŽ torique pour laquelle $\deg_{{\bf P}_{2^N-1}}(\overline{s(G)})\ \rlap{$\geq$}\raise-7pt\hbox{$\sim$}\ n!\left(\frac{N}{n}\right)^n\deg_{{\bf P}_N}(\overline G)$, laissant supposer que l'inŽgalitŽ de droite dans la proposition~\ref{compdegseg} pourrait tre amŽliorŽe d'un facteur ${\rm e}^n$. De fait, lorsque $n=2$ on vŽrifie pour toute configuration ${\cal A}$ de points et tout $u\in Q_{\cal A}$ la majoration $f_{\cal A}(u) \leq \left(\frac{N}{2}\right)^2$, d'o $\deg_{{\bf P}_{2^N-1}}(\overline{s(G)}) \leq 2\left(\frac{N}{2}\right)^2\deg_{{\bf P}_N}(\overline G)$, qui est optimal.
\end{exmpl}

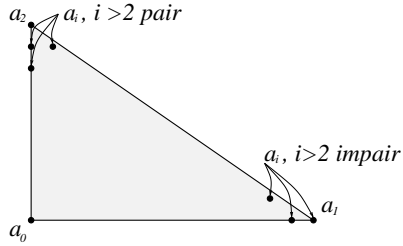
\vspace{-1mm} 

\begin{figure}[htbp]
$$\kern0.4cm\epsfig{file=figures/droite.eps, height= 32mm}$$
\end{figure}


%

\typeout{Un théorème de Bézout pour les poids de Chow}

\section{Un théorème de Bézout pour les poids de Chow}\label{Chow}

Soit $X\subset{\bf P}^N$ et $\tau=(\tau_0,\dots,\tau_N)\in{\bf R}^{N+1}$, 
rappelons que $e_\tau(X)\in \R$ dŽsigne le $\tau$-poids de Chow de $X$ dŽfini 
dans l'introduction. 

\medskip

Pour $\lambda\in{\bf R}^+$ on a $e_{\lambda\tau}(X) = \lambda \, e_\tau(X)$. De mme, on vŽrifie facilement que pour $\tau'\in{\bf R}$ on a $e_{\tau+(\tau',\dots,\tau')}(X) = e_\tau(X) + \tau'(n+1)\deg(X)$. On vŽrifie encore que si $\sigma_\delta: {\bf P}^N\rightarrow {\bf P}^M$, $M+1={{N+\delta}\choose N}$, dŽsigne le plongement de Veronese de degrŽ $\delta$ et $\tau^{(\delta)}$ le vecteur des $\tau$-poids des mon™mes de degrŽs $\delta$, on a $e_{\tau^{(\delta)}}(\sigma_\delta(X)) = \delta^{n+1} \cdot e_\tau(X)$. 
Plus profonde est l'expression de $e_\tau(X)$ comme coefficient dominant d'une {\it fonction poids de Hilbert} dŽmontrŽe par Mumford~\cite[Prop.~2.11]{Mum}~:
$$
s_\tau(X;D) := \max_J\sum_{\lambda\in J}(\tau_0\lambda_0+\dots+\tau_N\lambda_N) = \frac{e_\tau(X)}{(n+1)!}D^{n+1} + O(D^n)
\enspace,$$
o le maximum porte sur tous les ensembles $J$ d'ŽlŽments de ${\bf N}^{N+1}_D$  tels que les mon™mes associŽs induisent une base de la partie graduŽe de degrŽ $D$ de l'anneau de la variŽtŽ $X$.

\begin{exmpl}\label{poidshyper}
Si $H$ est une hypersurface de $\P^N$ d'Žquation $f\in\K[x_0,\dots,x_N]$ et $\tau\in\R^{N+1}$, on a
$$
e_\tau(H) = (\tau_0+\dots+\tau_N)\deg(H) - w_t(\lambda_\tau^*(f))
$$
o $\lambda_\tau^*(f) := f(t^{\tau_0}x_0,\dots,t^{\tau_N}x_N)$ et $w_t$ dŽsigne la valuation $t$-adique.
\end{exmpl}

\medskip

Lorsque $\tau\in\Z^{N+1}$ le poids de Chow peut s'interprŽter en termes de bidegrŽ de la dŽfor\-ma\-tion torique $X_\tau \subset \P^1\times \P^N$ de $X$ relative {\`a} $\tau$. Une forme rŽsultante d'indice $(1,n+1)$ de $X_\tau$ s'Žcrit 
$$v_1^{e_\tau(X)} v_0^{e_{-\tau}(X)} Ch_X(\dots,(-v_0/v_1)^{\tau_j}u_{i,j},\dots)\enspace,$$
son degrŽ en $(v_0,v_1)$ est $e_\tau(X)+e_{-\tau}(X)$ d'aprs la dŽfinition du poids de Chow et c'est par ailleurs le bidegrŽ $\deg_{(0,n+1)}(X_\tau)$, d'o 
$$\deg_{(0,n+1)}(X_\tau) = e_\tau(X) + e_{-\tau}(X)\enspace.$$
Rappelons que le degrŽ $\deg_{(0,n+1)}$ est obtenu en intersectant par $n+1$
formes linŽaires relevŽes du second facteur $\P^{2N+1}$, {\it voir} \S~\ref{vartor} pour plus de dŽtails.

\smallskip 

Lorsque $\tau\in\N^{N+1}$ on peut aussi interprŽter le poids de Chow $e_\tau(X)$ isolŽment comme un bidegrŽ en considŽrant une variation de la dŽformation torique prŽcŽdente. PrŽcisŽment, il s'agit maintenant de l'adhŽrence de Zariski $\wt X_\tau\subset{\bf P}^1\times{\bf P}^{2N+1}$ de
$$\{(1:t)\times(t^{\tau_0}x_0:\dots:t^{\tau_N}x_N:x_0:\dots:x_N);t\in{\bf G}_m,x\in X\}
\enspace.$$
En appliquant ce qui prŽcde on a
$$\deg_{(0,n+1)}(\wt X_\tau) = e_{\wt\tau}(\wt X) + e_{-\wt\tau}(\wt X)\enspace,$$
o $\wt X$ dŽsigne le {\it plongement diagonal} de $X$ dans ${\bf P}^{2N+1}$ et $\wt\tau=(\tau_0,\dots,\tau_N,0,\dots,0)\in{\bf N}^{2N+2}$. On vŽrifie sans difficultŽ $e_{\wt\tau}(\wt X) = e_\tau(X)$ et $e_{-\wt\tau}(\wt X) = 0$, d'o $\deg_{(0,n+1)}(\wt X_\tau) = e_{\tau}(X)$.

\smallskip 

Donnons maintenant la dŽmonstration du thŽorme~\ref{BezoutpoidsChow}, en commenant par une variante de ce thŽorme de BŽzout pour les poids de Chow. Notons 
$$\pi:\P^1\times\P^{2N+1} \to \P^N\enspace,\quad ((t_0:t_1),(x_0:\dots:x_N:y_0:\dots:y_N)) \mapsto (y_0:\dots:y_N)$$ 
la projection de $\P^1\times\P^{2N+1}$ sur $\P^N$ donnŽe par les $N+1$ dernires coordonnŽes de $\P^{2N+1}$ et 
$$\tilde\iota:\P^N \to \P^1\times\P^{2N+1} \enspace,\quad (x_0:\dots:x_N) \mapsto  ((0:1),(x_0:\dots:x_N:0:\dots:0))
\enspace.$$

\begin{thm}\label{Bezoutvar}
Soit $X\subset\P^N$ une variŽtŽ projective et $H$ un diviseur de $\P^N$ ne contenant pas $X$, alors pour tout $\tau\in\Z^{N+1}$ on a
$$e_\tau(X\cdot H) = e_\tau(X)\,\deg(H) - \sum_{Y\in\irr(\initial_\tau(X))} m(\tilde X_\tau\cdot \pi^*(H);\tilde\iota(Y))\,\deg(Y) 
\enspace.$$
De plus, si $H$ est une hypersurface et $f\in\K[x_0,\dots,x_N]$ est une Žquation de $H$ on a, avec les notations du thŽorme~\ref{BezoutpoidsChow},
$$m(\tilde X_\tau\cdot \pi^*(H);\tilde\iota(Y)) = m(X_\tau \cdot \div(\lambda_{\tau}^*(f));\iota(Y))
\enspace,$$
o $\lambda_{\tau}^*(f) := f(t_0^{\tau_0}x_0,\dots,t_0^{\tau_N}x_N) \in \K[t_0,t_1][x_0,\dots,x_N]$.
\end{thm}

\begin{demo}
Il suffit d'\'etablir le rŽsultat pour $\tau\in{\bf N}^{N+1}$, on se ramne ˆ ce cas en ajoutant ˆ $\tau$ un vecteur $(\tau',\dots,\tau')$ o $\tau'\in{\bf N}$ est suffisamment grand et on remarque qu'avec les propriŽtŽs d'homogŽnŽitŽ du poids de Chow l'ŽgalitŽ dŽsirŽe reste invariante.

Dans ce cas, comme $H$ ne contient pas $X$, le cycle intersection $\wt X_\tau\cdot \pi^*(H)$ s'Žcrit comme la somme du cycle $(\widetilde{X\cdot H})_\tau$ et d'un cycle supportŽ par $\tilde\iota(\init_\tau(X))$, le second vivant donc dans $\{(0:1)\}\times\P^{2N+1}$. Comme $\deg_{(0,n+1)}(\wt X_\tau) = e_\tau(X)$ et $\deg_{(0,n)}((\widetilde{X\cdot H})_\tau) = e_\tau(X\cdot H)$ on a, d'aprs le thŽorme de BŽzout multi-projectif \cite[Thm.~3.4]{Rem01b},
$$\begin{array}{rcl}
\deg(H) e_\tau(X) &= &\deg(H) \deg_{(0,n+1)}(\wt X_\tau)\\
&= &\deg_{(0,n)}(\wt X_\tau \cdot \pi^*(H))\\
&= &\displaystyle\deg_{(0,n)}(\widetilde{(X \cdot H})_\tau) + 
\sum_{Y\in \irr(\initial_\tau(X))} m(\wt X_\tau \cdot \pi^*(H); \tilde\iota(Y)).\deg(Y)\\
&= &\displaystyle e_\tau(X\cdot H) + \sum_{Y\in \irr(\initial_\tau(X))} m(\wt X_\tau \cdot \pi^*(H); \tilde\iota(Y)).\deg(Y)
\enspace.\end{array}$$

Notons finalement que le morphisme
$$\begin{array}{rccc}
\pi': &\P^1\times\P^{2N+1} &\to &\P^1\times\P^N\\
&((t_0:t_1),(x_0:\dots:x_N:y_0:\dots:y_N)) &\mapsto  &((t_0:t_1),(x_0:\dots:x_N))
\end{array}$$
est un isomorphisme de $\wt X_\tau$ sur $X_\tau$ dans l'ouvert $t_1\ne0$, tel que $\pi'\circ\tilde\iota=\iota$ et $\pi'_*(\pi^*(H)) = \div(\lambda_{\tau}^*(f))$. Cela entra{\^\i}ne, pour toute composante irrŽductible $Y$ de $\init_\tau(X)$,
$$
m(\wt X_\tau \cdot \pi^*(H);\tilde\iota(Y)) = m({X_\tau}\cdot \div(\lambda_{\tau}^*(f));\iota(Y)) \enspace. 
$$
\end{demo}

\begin{demo}[DŽmonstration du thŽorme~\ref{BezoutpoidsChow}]
On reprend la fin de la dŽmonstration du thŽorme~\ref{Bezoutvar} en se placant dans une carte affine contenant une composante $Y$ de $\init_\tau(X)$ et sur laquelle $H$ est dŽcrit par une Žquation $f$, on a
$$\begin{array}{rcl}
m(\wt X_\tau \cdot \pi^*(H);\tilde\iota(Y))
&\kern-9pt= &\kern-9pt{\rm long}(\K[\tilde X_\tau]/(\pi^*(f)))_{\tilde\iota(Y)}\\[3mm] 
&\kern-9pt= &\kern-9pt{\rm long}(\K[X_\tau]/(\lambda_\tau^*(f)))_{\iota(Y)}\\[3mm]
&\kern-9pt= &\kern-9pt m(X_\tau\cdot \div(\lambda_{\tau}^*(f));\iota(Y))\\[3mm]
&\kern-9pt= &\kern-9pt m(X_\tau\cdot H_\tau;\iota(Y)) + w_{t_0}(\lambda_\tau^*(f))\, m(\init_\tau(X);Y)\\[3mm]
&\kern-9pt= &\kern-9pt m(X_\tau\cdot H_\tau;\iota(Y)) + ((\tau_0+\dots+\tau_N)\deg(H) - e_\tau(H))\,m(\init_\tau(X);Y)
\end{array}$$
d'aprs le calcul de l'exemple~\ref{poidshyper}. Enfin, en sommant sur toutes les composantes $Y$ de $\init_\tau(X)$ on a $\sum_Y m(\init_\tau(X);Y)\deg(Y) = \deg(\init_\tau(X)) = \deg(X)$.
\end{demo}

Dans les notations du thŽorme~\ref{Bezoutvar}, soit $f\in\K[x_0,\dots,x_N]$ une Žquation de $H$ et $\tau\in\R^{N+1}$, posons
\begin{equation}\label{defvaluationrelative}
w_{X,\tau}(f) := \frac{1}{\deg(X)}.\left(e_\tau(X\cdot H) - e_\tau(X)\deg(H)\right)
\enspace,\end{equation}
c'est une fonction continue de $\tau$. Lorsque $\tau\in\N^{N+1}$ on a donc par le thŽorme~\ref{Bezoutvar}
$$w_{X,\tau}(f) = -\sum_{Y\in\irr(\initial_\tau(X))} m(\tilde X_\tau\cdot \pi^*(H);\tilde\iota(Y)).\frac{\deg(Y)}{\deg(X)}
\enspace.$$
Toujours dans ce cas, on vŽrifie $w_{X,k\tau}(f) = kw_{X,\tau}(f)$
pour tout $k\in\N^\times$ et $w_{X,\tau+\tau'(1,\dots,1)}(f) =
w_{X,\tau}(f) - \tau'\deg(H)$ pour $\tau'\in\N$. On peut donc
Žcrire en gŽnŽral (\cad pour $\tau\in\R^{N+1}$) l'ŽgalitŽ
$$e_\tau(X\cdot H) = e_\tau(X)\,\deg(H) + w_{X,\tau}(f)\,\deg(X)$$
en posant
$$w_{X,\tau}(f) := \lim_{k\rightarrow\infty}\frac{1}{k}w_{X,[k(\tau-\tau'(1,\dots,1))]}(f) + \tau'\deg(H)$$
o $\tau':=\min(\tau_0,\dots,\tau_N)$ et $[\cdot]$ dŽsigne le
vecteur des parties entires dans $\N^{N+1}$.

\medskip

Dans le cas d'une variŽtŽ torique $X_{{\cal A},\alpha}$ le poids de Chow $e_\tau(X_{{\cal A},\alpha}) = e_\tau(X_{\cal A})$ s'exprime naturellement en termes d'un polytope construit ˆ l'aide du vecteur $\tau$, au-dessus du polytope $Q_{\cal A}$ dŽjˆ introduit au paragraphe~\ref{vartor}.

\begin{prop}{(\cite[Prop.~III.1]{PS04}) ---} \label{poidsChow}
Avec les notations introduites on suppose ${\bf Z}a_0+\dots+{\bf Z}a_N={\bf Z}^n$. Soit $Q_{\cal A,\tau}$ l'enveloppe convexe des points $(a_0,\tau_0),\dots,(a_N,\tau_N)$ dans ${\bf R}^{n+1}$ et $\vartheta_{{\cal A},\tau}:Q_{\cal A}\rightarrow{\bf R}$ la paramŽtrisation de la toiture de $Q_{{\cal A},\tau}$ au-dessus de $Q_{\cal A}$, alors
$$e_\tau(X_{\cal A}) = (n+1)!\int_{Q_{\cal A}}\vartheta_{{\cal A},\tau}(u)du_1 \dots du_n \enspace.$$
\end{prop}

La d{\'e}monstration qu'on donne de cette proposition dans~\cite[\S~III]{PS04} est basŽe sur un calcul explicite des poids de Hilbert des vari{\'e}t{\'e}s toriques;  nous y montrons {\'e}galement qu'on peut l'obtenir comme cons{\'e}quence des r{\'e}sultats de I.M.~Gelfand, M.M.~Kapranov et A.V.~Zelevinski sur le polytope de Newton du $\cA$-r{\'e}sultant~\cite[Ch.~7 et 8]{GKZ94}. Elle est {\'e}galement implicite dans~\cite[\S~4.2]{Don02}. L'interprŽtation ci-dessus du poids de Chow comme bi-degr{\'e} permet d'en donner encore une autre d{\'e}monstration simple et directe~: 

\begin{demo} 
Tout d'abord on remarque qu'il suffit de d{\'e}montrer l'{\'e}nonc{\'e} pour un
choix g{\'e}n{\'e}rique (au sens de Zariski) du vecteur $\tau$ dans
$\R^{N+1}$, puisque les termes consid{\'e}r{\'e}s sont continus par rapport {\`a}
$\tau$.
L'identit{\'e} {\'e}tant invariante par homoth{\'e}ties et translations, 
on peut se ramener ˆ supposer $\tau \in \N^{N+1}$ {\`a} coordonn{\'e}es
premi{\`e}res entre elles dans leur ensemble. 

On a $L_{\cA,\tau}=\Z^{n+1}$; la proposition~\ref{multidegvol}
entra{\^\i}ne alors
$$e_{\tau}(X) = \deg_{(0,n+1)}(\wt X_\tau) 
= \MV_{n+1}(Q_{\cA,\wt\tau}, \dots,Q_{\cA,\wt\tau}) 
= (n+1)!\int_{Q_{\cal A}}\vartheta_{{\cal A},\tau}(u) du_1 \cdots du_n
\enspace,
$$
o $\wt\tau = (\tau_0,\dots,\tau_N,0,\dots,0) \in \N^{2N+2}$.
\end{demo} 

\medskip

La proposition~\ref{poidsChow} permet de reformuler de faon amusante le thŽorme du sous-espace de Schmidt, dans la version projective qu'en ont donnŽ Evertse et Ferretti~\cite{EF02} mais dans le cas trs particulier des variŽtŽs toriques.

\begin{convention}\label{convention}
Si $K$ est un corps de nombres et $v$ une place de $K$ on notera $|\cdot|_v$ la valeur absolue de $K$ dans $v$ Žtendant la valeur absolue usuelle de $\Q$ si $v$ est archimŽdienne et la valeur absolue $p$-adique standard de $\Q$ ({\it i.e.} $|p|_v=p^{-1}$) si $v$ est ultramŽtrique. On notera encore $M_K$ l'ensemble de ces valeurs absolues reprŽsentant les places de $K$. 
\end{convention}

Soit $K$ un corps de nombres et $S$ un ensemble fini de places de $K$. On considre une variŽtŽ torique $X_{{\cal A},\alpha}\subset{\bf P}_N$ dŽfinie sur $K$, pour toute place $v\in S$ des rŽels $\tau_{v,0},\dots,\tau_{v,N}\geq 0$ et le systme d'inŽquations en $x\in X_{{\cal A},\alpha}(\Qbar)$~:
\begin{equation}\label{sousespace}
\log\left(\frac{|x_i|_w}{\Vert x\Vert_w}\right) \leq -\tau_{v,i}.h(x) \qquad \hbox{\rm pour tout $i=0,\dots,N$ et toute place $w\mid v$ de $K(x)$,}
\end{equation}
o $\Vert x\Vert_w:=\max(|x_i|_w;i=0,\dots,N)$ si $w$ est ultramŽtrique, $\Vert x\Vert_w = \left(\sum_{i=0}^N|x_i|_w^2\right)^{1/2}$ si $w$ est archimŽdienne, et $h(x)$ dŽsigne la hauteur projective. 

Dans cette situation notons $\vartheta_{{\cal A},\tau}(u) := \sum_{v\in S}\frac{[K_v:{\bf Q}_v]}{[K:{\bf Q}]}.\vartheta_{{\cal A},\tau_v}(u)$ pour $u\in Q_{\cal A}$. Le thŽo\-r\-me~3.2 de~\cite{EF02} explicite pour tout $\varepsilon>0$ des rŽels $$c_1(N,n,\deg(X_{{\cal A},\alpha}),\varepsilon)\enspace,\enspace c_2(N,n,\deg(X_{{\cal A},\alpha}),\varepsilon)\geq 1$$ 
tels que l'ŽnoncŽ suivant soit vrai ($\tau_v:=(\tau_{v,0},\dots,\tau_{v,N})$, $v\in S$)~:

\begin{thm}\label{schmidttorique} (\cite{EF02}, Thm.~3.2) --
Si la valeur moyenne de $\vartheta_{{\cal A},\tau}$ sur $Q_{\cal A}$ est $\geq 1+\varepsilon$, alors les points $x\in X_{{\cal A},\alpha}(\overline{\bf Q})$ solutions du systme d'inŽquations~(\ref{sousespace}) et satisfaisant 
$$h(x) \geq c_1(N,n,\deg(X_{{\cal A},\alpha}),\varepsilon).(1+h(X_{{\cal A},\alpha}))$$
appartiennent ˆ un sous-ensemble algŽbrique propre de $X_{{\cal A},\alpha}$, dŽfini sur $K$ et de degrŽ $\leq c_2(N,n,\deg(X_{{\cal A},\alpha}),\varepsilon)$.
\end{thm}
\begin{demo}
Notons $X=X_{{\cal A},\alpha}$, la condition donnŽe dans le thŽorme~3.2 de~\cite{EF02} s'Žcrit
$$\frac{1}{(n+1)\deg(X)}.\sum_{v\in S}\frac{[K_v:{\bf Q}_v]}{[K:{\bf Q}]}.e_{\tau_v}(X) \geq 1+\varepsilon
\enspace.$$
Or $\deg(X)=n!{\rm Vol}_n(Q_{\cal A})$ d'aprs la formule~(\ref{degvol}) et, vue la proposition~\ref{poidsChow}, on a 
$$\sum_{v\in S}\frac{[K_v:{\bf Q}_v]}{[K:{\bf Q}]}.\frac{e_{\tau_v}(X)}{(n+1)\deg(X)} = \frac{1}{{\rm Vol}_n(Q_{\cal A})}. \int_{Q_{\cal A}} \vartheta_{{\cal A},\tau}(u)du_1\dots du_n$$ 
qui est bien la valeur moyenne de $\vartheta_{{\cal A},\tau}$ au-dessus $Q_{\cal A}$.
\end{demo}
\medskip

Le rŽsultat trivial, provenant de la formule du produit, permet d'affirmer que le sytme (\ref{sousespace}) n'a pas de solution dans ${\bf P}_N^\circ(\overline{\bf Q})$ ds qu'il existe $i\in\{0,\dots,N\}$ tel que 
$$\sum_{v\in S}\frac{[K_v:{\bf Q}_v]}{[K:{\bf Q}]}.\tau_{v,i} > 1
\enspace.$$ 
On donne un exemple o ce rŽsultat trivial s'applique bien que la condition du thŽorme~\ref{schmidttorique} ne soit pas remplie et un autre exemple o le thŽorme~\ref{schmidttorique} donne un rŽsultat non trivial.

\begin{exmpl}
Soient $0<\varepsilon<n+1$, ${\cal A}\subset{\bf Z}^n$ et $\tau_v$ tels que $\tau_{v,i}=0$ pour tout $v\in S$ et $i=1,\dots,N$ et $\sum_{v\in S}\frac{[K_v:{\bf Q}_v]}{[K:{\bf Q}]} \tau_{v,0} = 1+\varepsilon>1$. On calcule ˆ l'aide de la proposition~\ref{poidsChow} $e_{\tau_v}(X_{\cal A}) \leq \deg(X_{\cal A}).\tau_{v,0}$, d'o $\sum_{v\in S}\frac{[K_v:{\bf Q}_v]}{[K:{\bf Q}]}.\frac{e_{\tau_v}(X_{\cal A})}{(n+1)\deg(X_{\cal A})} \leq \frac{1+\varepsilon}{n+1} <1$. On sait donc que dans ce cas le systme d'inŽquations~(\ref{sousespace}) n'a pas de solution telle que $x_0\not=0$, bien que l'hypothse principale du thŽorme~\ref{schmidttorique} ne soit pas satisfaite.
\end{exmpl}

\begin{exmpl}
Soit maintenant $N'>3$ un entier, $N:=2N'$ et $D>N'(N'-1)$ un autre entier. On considre $\cA:=(0,\dots,N'-1,D-N'+1,\dots,D)\in{\bf Z}$ et $S=\{v_1,\dots, v_{N'}\}$ un ensemble de $N'$ places de ${\bf Q}$. \smash{\`A} chaque place $v_i$ de $S$ on associe le vecteur $\tau_{v_i}\in{\bf R}^{N+1}$ dont toutes les coordonnŽes sont nulles sauf celles d'indices $i$ et $D-N'+i$ qui valent $\frac{1}{N'-1}$. On a ainsi $e_{\tau_{v_i}}(X_\cA) \geq 2\frac{D-N'+1}{N'-1}$ et $\sum_{v\in S}\frac{e_{\tau_v}(X_\cA)}{2\deg(X_\cA)} \geq \frac{D-N'+1}{D}.\frac{N'}{N'-1} >1$, tandis que $\sum_{v\in S}\tau_{v,i}=\frac{1}{N'-1}<1$ pour tout $i=0,\dots,N$. Le thŽorme~\ref{schmidttorique} s'applique donc et donne un rŽsultat {\em a priori} non trivial dans ce cas.
\end{exmpl}


\typeout{Hauteur normalisee}
\section{Hauteur normalisŽe} \label{hauteurs}

L'espace projectif peut {\^e}tre vu comme 
une compactification {\'e}quivariante du groupe multiplicatif $\G_m^N \simeq (\P^N)^\circ $, cette structure 
de groupe permettant de d{\'e}finir une 
notion de  hauteur pour les sous-vari{\'e}t{\'e}s de $\P^N$ plus canonique que
les autres, appel{\'e}e {\it hauteur normalis{\'e}e}. Cette notion joue un
r{\^o}le central dans l'approximation diophantienne sur les tores, et tout
particuli{\`e}rement dans les probl{\`e}mes de Lehmer g{\'e}n{\'e}ralis{\'e} et de Bogomolov sur les tores, 
{\it voir}~\cite{DP99}, \cite{AD} et leurs r{\'e}f{\'e}rences,
{\it voir} {\'e}galement~\cite{Dav02} pour un aper{\c c}u historique.

Suivant~\cite{DP99}, la hauteur normalis{\'e}e peut se d{\'e}finir par un
proc{\'e}d{\'e} \ovg {\`a} la Tate{}\fmg. 
De fa{\c c}on pr{\'e}cise, pour $k \in \N$ 
on pose $ [k]: \P^N \to \P^N$, $(x_0: \cdots : x_N) \mapsto (x_0^k: \cdots : x_N^k)$ l'application puissance $k$-i{\`e}me; restreinte au tore $(\P^N)^\circ$ 
c'est l'application de multiplication par $k$. La {hauteur normalis{\'e}e}  
d'une vari{\'e}t{\'e}  projective $X$
est par d{\'e}finition
$$
\wh{h}(X) := \deg(X) \cdot \lim_{k\to \infty} \, 
\frac{h({[k] \, X})}{k\, \deg([k]\, X)} \ \in \R_+
\enspace, 
$$ 
o{\`u} $\deg$ et $h$ d{\'e}signent  le degr{\'e} et la  hauteur projective, {\it
  voir}~\cite[\S \ 2]{DP99} ou~\cite[\S~I.2]{PS04}. 
Lorsque $X$ est de dimension $0$, il s'agit de la hauteur de Gauss-Weil des points.

\smallskip 

On sait que le comportement de cette hauteur normalisŽe est lié à la nature gŽomŽtrique de la variŽtŽ $X$. En particulier, $\wh{h}(X)=0$ si et seulement si $X$ est une variŽtŽ de torsion, \cad un translatŽ de 
sous-tore un point de torsion.
Dans~\cite{PS04} on a {\'e}tabli un analogue arithm{\'e}tique de l'expression du degré d'une variété torique comme volume du polytope associé, que nous rappelons maintenant. 

Soit $M_K$ l'ensemble des places du corps $K$ et rappelons la convention~\ref{convention}. Pour chaque $v \in M_K$ on consid{\`e}re le vecteur $\tau_{\alpha v} := (\log|\alpha_0|_v,\dots,\log|\alpha_N|_v)\in\R^{N+1}$ et le polytope 
$$
Q_{\cA,\tau_{\alpha v}}:= \Conv \Big( (a_0,\log|\alpha_0|_v), 
\dots, (a_N,\log|\alpha_N|_v)  \Big) \ \subset \R^{n+1} 
\enspace, 
$$
dont la {\it toiture} au-dessus de $Q_\cA$ (\cad l'enveloppe sup{\'e}rieure) s'envoie bijectivement sur $Q_\cA$ par la projection standard $ \R^{n+1} \to \R^n$. 
On pose alors 
$$
\vartheta_{\cA,\tau_{\alpha  v}}  : Q_\cA \to \R
\enspace , \quad \quad 
x \mapsto \max \Big\{ y\in\R\, : \ (x,y) \in Q_{\cA,\tau_{\alpha v}} \Big\} 
$$ 
la param{\'e}trisation de cette toiture; c'est une fonction {\it concave} 
et {\it affine par morceaux}. 

\begin{thm} \label{thmprincipal} (\cite[Thm.~0.1]{PS04}) ---
Soit $\geom \in (\Z^n)^{N+1}$ tel que $L_\geom= \Z^n$
et $\alpha \in (K^\times)^{N+1}$, alors
$$
\wh{h}(X_{\cA,\alpha}) = \sum_{v\in M_K}{{[K_v:{\bf Q}_v]}\over{[K:{\bf
      Q}]}} \, e_{\tau_{\alpha\,v}}(X_{\cA}) =  (n+1)! \, \sum_{v \in M_K} 
\frac{[K_v:\Q_v]}{[K:\Q]} \, \int_{Q_\cA} 
\vartheta_{\cA,\tau_{\alpha v}}(u)\, \ du_1 \cdots du_n 
\enspace .
$$
\end{thm} 

Notons que $\tau_{\alpha v} = 0$ pour presque tout $v$, donc cette somme ne contient qu'un nombre fini de termes non nuls. Dans les cas des points ($n=0$) la formule se r\'eduit \`a la d\'efinition usuelle de la hauteur normalis\'ee (hauteur de Gauss-Weil) d'un point de $\P^N$.

\medskip

Comme on a $\wh{h}(X_\cA)=0$, on peut s'interroger si cette formule n'est pas la manifestation d'une propriŽtŽ gŽnŽrale de la hauteur normalisŽe par translation. Dans ce sens, B.~Sturmfels nous a demandé si pour toute vari\'et\'e projective $X$ on a
$$\wh{h}(\alpha X) = \wh{h}(X) + \sum_{v\in M_K}{{[K_v:{\bf Q}_v]}\over{[K:{\bf Q}]}} e_{\tau_{\alpha\,v}}(X)
\enspace?$$
Soit $\prod_{i=0}^n\prod_{j=0}^NU_{i,j}^{b_{i,j}}$ un mon™me apparaissant dans la forme de Chow de $X$, on vŽrifie facilement 
$$e_{\tau_{\alpha\, v}}(X) 
\geq \left(\sum_{i=0}^nb_{i,0}\right)\tau_{\alpha\,v,0} + \dots + \left(\sum_{i=0}^nb_{i,N}\right)\tau_{\alpha\,v,N}
\enspace.$$
Comme $\displaystyle \sum_{v\in M_K}{{[K_v:{\bf Q}_v]}\over{[K:{\bf Q}]}} \tau_{\alpha\,v,i} = 0$ pour tout $i=0,\dots,N$ par la formule du produit, on en dŽduit 
$$\sum_{v\in M_K}{{[K_v:{\bf Q}_v]}\over{[K:{\bf Q}]}} 
e_{\tau_{\alpha\, v}}(X)  \geq 0 
\quad \mbox{ pour tout } \alpha\in(\Qbar^\times)^{N+1} \enspace.$$ 
Une formule du type prŽcŽdent ne peut donc tre valable en gŽnŽral.

\medskip
Illustrons cette formule par un exemple en dimension 1. 
Soient $\cA_N:=(0,1,2,\dots,N) \in \Z^{N+1}$ et $\alpha_N := (1,2,3,\dots,N+1) \in (\Qbar^\times)^{N+1}$; ainsi 
$$
\varphi_{\cA_N,\alpha_N} : \G_m \to \P^N \quad , \quad s \mapsto (1:2s:\dots:(N+1)s^N).$$ 
Les figures suivantes montrent pour $N=3$ les polytopes associ{\'e}s et leurs toitures, pour chaque place $v \in M_\Q$:




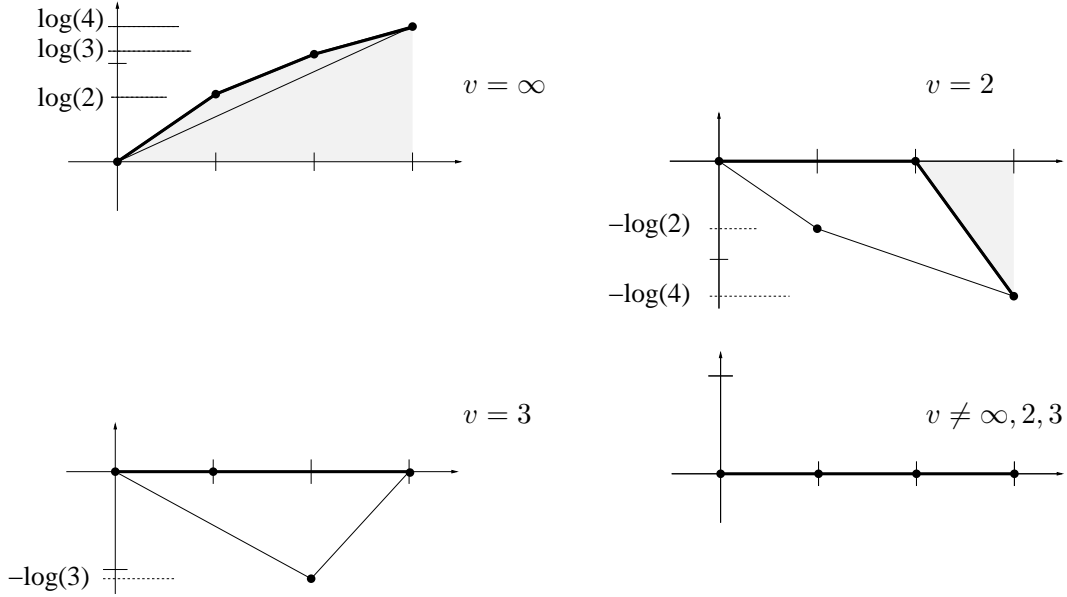
$$\raise-39mm\vtop to 39mm{

\begin{picture}(0,0) 

\put(30,50){
\epsfig{file=figures/infini.eps, height= 28 mm} 
} 
\put(195,95){$v=\infty$} 

\put(246,12.5){
\epsfig{file=figures/2.eps, height= 26.5 mm} 
} 
\put(370,95){$v=2$}

\put(19,-95){
\epsfig{file=figures/3.eps, height= 23 mm} 
} 
\put(195,-30){$v= 3$} 

\put(270,-68){
\epsfig{file=figures/autre.eps, height= 23 mm} 
} 
\put(370,-30){$v\ne \infty, 2,3$}

\end{picture} 

\vfil}$$



\goodbreak 


\noindent Ainsi $\vartheta_v \equiv 0$ pour $v \ne \infty, 2$, d'o{\`u} 
$$
\hnorm(X_{\cA_3,\alpha_3}) = 2!\,\bigg(\int_0^3\,\vartheta_\infty(u)\, du\
+\int_0^3\,\vartheta_2(u)\, du\bigg) = 2\, \log(2) + 2\, \log(3) = \log(36)
\enspace.$$
En g{\'e}n{\'e}ral
\begin{eqnarray*}
\wh{h}(X_{\cA_N,\alpha_N}) &= &2!\,\bigg(\int_0^N\,\vartheta_\infty(u)\, du\
+ \sum_{p \mid N+1} \int_0^N\,\vartheta_p(u)\, du \bigg) \\[0mm] 
&= &2\, \bigg( \log(2)+\dots +\log(N) \ + \frac12 \, \log(N+1) + \sum_{p\mid N+1} \frac12 \, \log|N+1|_p\bigg) \\[0mm]  
&= &2\, \log(N!)\enspace.
\end{eqnarray*}

\smallskip 

{\`A} l'instar des multidegrŽs, les multihauteurs du tore $\G_m^n$ plongŽ dans un produit d'espaces projectifs {\it via} plusieurs applications monomiales peuvent aussi s'expliciter ˆ l'aide {\it d'intŽgrales mixtes} des fonctions concaves apparaissant dans le thŽorme~\ref{thmprincipal}. 
 
Soit $Z\subset\P^{N_0}\times\dots\times\P^{N_m}$ de dimension $n$ et $c=(c_0,\dots,c_m)\in\N^{m+1}_{n+1}$ avec $0\le c_i \le N_i$. 
Soit $\res_{d(c)}(I(Z)))$ la forme résultante associée, dans les notations du \S~\ref{vartor}. 
La {\it multihauteur projective de $Z$ d'indice $c$} est d{\'e}finie par 
$$
h_{c}(Z) := h(\res_{d(c)}(I(Z))) \enspace, $$
où $h$ désigne la hauteur des polynômes multihomogènes définie à l'aide de
la {\it $S_{N_0+1} \times \cdots \times S_{N_n+1}$-mesure} pour les places archimediennes; on 
renvoie {\`a}~\cite{Rem01a} ou encore~\cite[\S~I.2]{PS04} pour les détails. 
Notons  $s:
\P^{N_0} \times \cdots \times \P^{N_m} \to
\P^{(N_0+1) \cdots (N_m+1) - 1}$
le plongement de Segre, 
la suite 
$$  k \mapsto  \deg(s(Z)) \cdot\frac{h_c([k] \, Z)}{k\, \deg([k]\, s(Z))}$$
converge lorsque $k$ tend vers l'infini~\cite[Prop.~I.2]{PS04}). Sa limite est par définition  
la {\it multihauteur normalis{\'e}e de $Z$ d'indice $c$}, notŽe $\hnorm_c(Z)$. 

Pour
 $D=(D_0,\dots,D_m)\in(\N^\times)^{m+1}$ soit  
$
\Psi_D :  \P^{N_0} \times \cdots \times \P^{N_m}  \to  
\P^{ {D_0+N_0 \choose N_0} \cdots {D_m+N_m \choose N_m} -1} $
le {plongement mixte} introduit au paragraphe~\ref{vartor}, alors on a la formule~\cite[(I.4)]{PS04}
\begin{equation}\label{pinson}
\hnorm(\Psi_D(Z)) = \sum_{c\in\N^{m+1}_{n+1}}{n+1\choose c}\hnorm_c(Z)D^c
\enspace.\end{equation}

\medskip
Les multihauteurs des variétés toriques s'explicitent en termes de ce qu'on appelle 
des {\it intŽgrales mixtes}  d'une famille de fonctions concaves~\cite[\S~IV.3]{PS04}, 
notion analogue au volume mixte. 
Soient $f: Q \to \R$ et $g: R \to \R$ des fonctions concaves d{\'e}finies sur des ensembles convexes $Q, R \subset \R^n$ respectivement. On pose 
$$
f \boxplus g : Q+R \to \R
\enspace , \quad \quad
x \mapsto \max \{ f(y) + g(z) \, : \ y \in Q, \ z \in R, \ y+z =x\}
\enspace,  
$$
qui est une fonction concave d{\'e}finie sur la somme de Minkowski $Q+R$; on obtient ainsi une structure de semi-groupe commutatif sur l'ensemble des  fonctions concaves (dŽfinies sur des ensembles convexes).
Pour une famille de fonctions concaves 
$
f_0 : Q_0 \to \R, \dots,  f_n: Q_n\to \R$ l'{int{\'e}grale mixte} 
(ou {\em multi-int{\'e}grale}) est d{\'e}finie {\em via} la formule
\begin{equation*} 
\MI(f_0, \dots, f_n) := 
\sum_{j=0}^n (-1)^{n - j} 
\sum_{0 \le i_0 < \cdots < i_j \le n} 
\int_{Q_{i_0} + \cdots + Q_{i_j}} (f_{i_0} \boxplus 
\cdots \boxplus f_{i_j})(u) \, du_1 \cdots du_n \enspace.
\end{equation*} 
Comme pour le volume mixte, l'intŽgrale mixte est une fonctionnelle positive,  sym{\'e}trique et lin{\'e}aire en chaque
variable $f_i$, {\it voir}~\cite[\S~IV.3]{PS04}. 

\smallskip

Soit $\cA_0\in(\Z^n)^{N_0+1}, \dots, \cA_m\in(\Z^n)^{N_m+1}$ tels 
que $L_{\cA_0}+\dots+L_{\cA_m}=\Z^n$ et $\underline{\cA}:=(\cA_0,\dots,\cA_m)$. Soit  $\alpha_0 \in (K^\times)^{N_0+1}, \dots , \alpha_m \in (K^\times)^{N_m+1}$ et posons $\un \alpha:=(\alpha_0,\dots,\alpha_m)$. Con\-si\-dŽ\-rons alors l'action monomiale $*_{\underline{\cA}}$ de $\G_m^n$ sur le produit d'espaces projectifs $\P^{N_0}\times\dots\times\P^{N_m}$ 
associŽe; on note $X_{\underline{\cA},\un \alpha}$ l'adhŽrence de Zariski de l'orbite du point $\un \alpha\in\P^{N_0}\times\dots\times\P^{N_m}$.

Pour chaque $v \in M_K$ on note aussi $\vartheta_{\cA_i,\tau_{\alpha_i\,v}}: Q_{\cA_i} \to \R$  la fonction param{\'e}trant la toiture du polytope $Q_{\cA_i,\tau_{\alpha_i\,v}} \subset \R^{n+1}$ associ{\'e} au vecteur $\cA_i$ et au poids $\tau_{\alpha_i \, v}$.

\begin{thm}{\cite[Thm.~0.3 et Rem.~IV.7]{PS04}}--- \label{multihauteurint}
Soit $c\in\N^{m+1}_{n+1}$, dans la situation ci-dessus on a 
$$\hnorm_{c}(X_{\underline{\cA},\un \alpha}) = \sum_{v\in M_K} \frac{[K_v:\Q_v]}{[K:\Q]} \, \MI_c( \vartheta_{\underline{\cA}, \tau_{\underline{\alpha}\, v}}) 
$$
avec $\MI_c( \vartheta_{\underline{\cA}, \tau_{\underline{\alpha}\, v}}) := \MI\Big( \underbrace{\vartheta_{\cA_0, \tau_{\alpha_0\, v}}, \dots, \vartheta_{\cA_0, \tau_{\alpha_0\, v}}}_{c_0 \mbox{ \scriptsize \rm fois}}
\enspace, \dots, \enspace  
\underbrace{\vartheta_{\cA_m, \tau_{\alpha_m\, v}}, \dots, \vartheta_{\cA_m, \tau_{\alpha_m\, v}}}_{c_m \mbox{ \scriptsize \rm fois}} \Big)$.
\end{thm}

\medskip

\begin{exmpl}
Soient $\xi_1,\dots,\xi_N\in K^\times$ et considŽrons l'application monomiale
$$\begin{array}{rcccl}
\G_m &\stackrel{\varphi}{\longrightarrow}&(\P^1)^N &
\stackrel{\mbox{\rm \small Segre}}{\longrightarrow} &\P^{2^N-1}\\
s &\longmapsto &((1:\xi_1s),\dots,(1:\xi_N s)) &\longmapsto & 
\big( (\prod_{j\in J}\xi_j) \cdot s^{{\scriptsize \rm Card}(J)} \, : J\subset \{1,\dots, N\}\big) \enspace. 
\end{array}$$
Soit $X$ l'image de $\varphi$ dans $(\P^1)^N$ et pour $1\leq i,j\leq
N$ notons $c(i,j)\in\N^N$ le vecteur dont les coordonnŽes d'indices
$i$ et $j$ valent $1$ et les autres $0$, on vŽrifie ˆ l'aide du
thŽorme~\ref{multihauteurint} et de la définition des multi-intégrales
$$\begin{array}{rcl}
\wh{h}_{c(i,j)}(X) &= &\sum_{v\in M_K} \frac{[K_v:\Q_v]}{[K:\Q]} \,
\MI_c( \vartheta_{\underline{\cA}, \tau_{\underline{\alpha}\,
    v}})\\[2mm] 
&= &\sum_{v\in M_K} \frac{[K_v:\Q_v]}{[K:\Q]} \,
\max(\log|\xi_i|_v;\log|\xi_j|_v)\\[2mm] 
&= &h(\xi_i:\xi_j)
\enspace.\end{array}$$
Dans la figure~\label{dessin4} suivante les parties grisŽes
montrent pour $N=2$ et en supposant $\log|\xi_i| \le \log|\xi_j|$,
pour chaque cas le calcul de la multi-intŽgrale $\MI_c( \vartheta_{\underline{\cA}, \tau_{\underline{\alpha}\,v}})$.



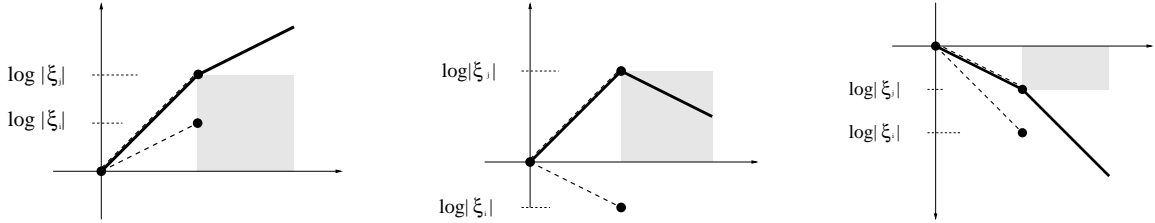
\begin{figure}[htbp]

\vtop to29mm{$$\raise-43mm\hbox{

\begin{picture}(0,0) 

\put(-218,55){\epsfig{file=figures/multi1.eps, height=29mm}} 

\put(-55,55){\epsfig{file=figures/multi2.eps, height=29mm}} 

\put(100,55){\epsfig{file=figures/multi3.eps, height=29mm}} 

\end{picture}

}$$}

\end{figure}



Par la formule~(\ref{pinson}), on en dŽduit que la hauteur de l'image de  
$\mbox{\rm Segre}\circ \varphi$  dans $\P^{2^N-1}$ est Žgale ˆ
$$\wh{h}(\mbox{\rm Segre}(X)) = 2\sum_{1\leq i<j\leq N}\wh{h}_{c(i,j)}(X) = 2\sum_{1\leq i<j\leq N} h((\xi_i:\xi_j)) = \sum_{1\leq i,j\leq N} h((\xi_i:\xi_j))\enspace.$$
\end{exmpl}

\medskip

Une consŽquence du thŽorme~\ref{thmprincipal} est que 
$\wh{h}(X_{{\cal A},\alpha})$ est le logarithme d'un nombre alg\'ebrique. 
En effet, $e_{\tau_{\alpha \, v}}(X_{\cal A})=0$ pour toutes les places $v$ sauf un nombre fini et, comme le poids de Chow est combinaison ˆ coefficients entiers des composantes du vecteur poids $\tau_{\alpha\,v}$, c'est ˆ dire des $\log|\alpha|_v$, notre assertion est claire. Un tel logarithme de nombre alg\'ebrique \'etant nul ou transcendant, 
on peut \'enoncer: 

\begin{prop}\label{transcendance}
Soit $X$ une variété torique qui n'est pas de torsion, alors $\wh{h}(X)\notin \Qbar$. 
\end{prop} 

Les hauteurs normalisŽes des variŽtŽs toriques qui ne sont pas de torsion, sont donc des nombres transcendants parmi les plus simples que l'on connaisse, ˆ savoir les logarithmes de nombres algŽbriques. On sait qu'il devrait rŽsulter des conjectures les plus gŽnŽrales sur les variŽtŽs et les motifs que les hauteurs des variŽtŽs projectives s'expriment rationnellement en termes de valeurs de fonctions $L$ et de leurs dŽrivŽes, ou encore en termes de polylogarithmes. 

Ds qu'on considre des variŽtŽs projectives qui ne sont plus toriques on doit s'attendre ˆ ce que leur hauteur normalisŽe fasse effectivement intervenir des polylogarithmes supŽrieurs. De mme, si l'on s'intŽresse ˆ la hauteur {\it projective} d'une variŽtŽ torique, une quantitŽ sup\-plŽ\-men\-taire s'ajoute aux 
intégrales des fonctions $\vartheta_{{\cal A},\tau_{\alpha\,v}}$, qui s'Žcrit bien sur des exemples ˆ l'aide de polylogarithmes mais est en gŽnŽral difficile ˆ Žvaluer exactement (calculs non publi\'es).

%

\typeout{Optimalité du th{\'e}or{\`e}me des minimums algŽbriques successifs} 

\section{Optimalité du th{\'e}or{\`e}me des minimums algŽbriques successifs} 

\label{minimums}

Dans ce paragraphe on étudie les minimums alg\'ebriques successifs des 
vari\'et\'es toriques, d\'efinis dans l'introduction.   
Comme résultat principal, on démontre le théorème~\ref{densite} \'enonc\'e 
dans l'introduction qui entra\^ine l'optimalité des estimations 
dans le théorème des minimums alg\'ebriques successifs: on construit des 
exemples montrant
qu'à des $\varepsilon$-près, toute configuration 
possible des minimums successifs 
se réalise et que le quotient ${\wh{h}(X) }/{\deg (X)}$
peut atteindre n'importe quel valeur dans l'intervalle autorisé par
les in\'egalit\'es~(\ref{zhang}). 

\medskip
Notons que l'analogue abélien du théorème~\ref{densite} est faux: 
soit $A$ une variété abélienne 
définie sur $\Qbar$ 
munie d'un fibré en droites $L$  ample et symétrique, permettant de définir une notion de hauteur normalisé
$\hnorm = \hnorm_L$ pour les sous-variétés de $A$ \cite{Phi91}.
En particulier, on a une notion de minimums successifs
et le théorème de Zhang y est toujours valable, {\it voir}~\cite[Thm.~5.2]{Zha95}. 

Soit  $\alpha + B \subset A$ le translat{\'e} d'une sous-vari{\'e}t{\'e} 
ab{\'e}lienne $B$ par un point $\alpha$
et $\Tors(B)$ le sous-groupe des points de torsion de $B$. 
Soit $\beta$ un point quelconque dans $\alpha+B$, alors 
$\beta+\Tors(B)$ est un sous-ensemble de points de hauteur 
$\hnorm(\beta)$ dense dans $\alpha+B$
(puisque $\Tors(B)$ est Zariski dense dans $B$), 
donc 
$\wh{\mu}^\ess(\alpha+B) \le \hnorm(\beta)$. 
On en déduit 
$ \wh{\mu}^\ess (\alpha+B) = \wh{\mu}^\abs (\alpha+B)$ et donc
$$
\wh{\mu}_i(\alpha +B)= \wh{\mu}^\ess (\alpha +B)
 \quad \mbox{ pour } i=1,\dots,n+1 \enspace. 
$$
Ainsi, dans cette situation l'intervalle du théorème des minimums 
successifs se réduit à un point, et on a les égalités
$$
\frac{\wh{h}(\alpha +B)}{\deg(\alpha +B)} 
= (n+1) \, \wh{\mu}^\ess (\alpha +B) 
= \wh{\mu}_1 (\alpha +B) + \dots + \wh{\mu}_{n+1} (\alpha +B)\enspace. 
$$
La situation est plus riche dans le cas torique, 
la diff{\'e}rence tenant au fait 
que les translat{\'e}s de sous-tores de $(\P^N)^\circ$ 
ne sont pas des ensembles ferm{\'e}s. 

\medskip
Soit $\cA= \Big(a_0, \dots, a_N \Big)\in (\Z^n)^{N+1}$ et 
$\alpha= (\alpha_0, \dots, \alpha_N) \in (\Qbar^\times)^{N+1}$.
L'ouvert principal
$X_{\cA,\alpha}^\circ$ 
est le translat{\'e} de sous-tore $\alpha\cdot X_\cA$  
et avec le même raisonnement que pour le cas abélien 
\begin{equation} \label{ess=abs} 
\wh{\mu}^\ess (X_{\cA,\alpha}) =
\wh{\mu}^\ess (X_{\cA,\alpha}^\circ) 
= \wh{\mu}_i (X_{\cA,\alpha}^\circ) 
\quad \mbox{ pour } i=1,\dots,n+1 \enspace.
\end{equation} 
Cependant, les autres minimums successifs  de 
$X_{\cA,\alpha}$ dépendent des orbites de dimension inférieure, 
et peuvent donc différer du minimum essentiel.

\begin{lem} \label{mus toriques}
Avec les notations ci-dessus, 
pour $i=1, \dots, n+1$
$$
\wh{\mu}_i (X_{\cA,\alpha}) = \min \Big\{ \wh{\mu}^\abs(X_{\cA,\alpha,P}^\circ) \, : \, P
\in \Faces(Q_\cA), \ \dim(P) = n-i+1 \Big\}\enspace, 
$$
où $P$ parcours l'ensemble des faces $\Faces(Q_\cA)$ du polytope $Q_\cA$ de dimension $n-i+1$.
\end{lem}

\begin{demo}
Considérons la décomposition en orbites: 
$ X_\arith  = \bigcup_{P \in \faces(Q_\cA) } X_{\arith, P}^\circ $
(formule~(\ref{decomposition})).
C'est un recouvrement de $X_\arith$ et donc par~\cite[Lem.~2.2]{Som02}
on a 
\begin{eqnarray*}
\munorm_i (X_\cA) &=& 
\min \Big\{ \, \munorm_{\dim(P) - n +i} \Big(X^\circ_{\cA,\alpha,P}\Big) \, : \ 
\dim (P) \ge n-i+1 \,\Big\}  \\[1mm]
&=& \min \Big\{ \, \munorm^\abs \Big(X^\circ_{\cA,\alpha, P}\Big) \, : \ 
\dim (P) = n-i+1 \,\Big\}
\end{eqnarray*}
car ce minimum est atteint sur les faces de dimension $n-i+1$.
\end{demo}

Ainsi, le  calcul des minimums successifs 
de $X_{\cA,\alpha}$ 
se r{\'e}duit ˆ celui du minimum essentiel (ou absolu) 
d'un sous-tore. 
Le lemme suivant donne le minimum essentiel pour 
certaines variétés toriques particulières. 

\begin{lem}\label{tous-sur-a}
Soit $\cA= \Big(a_0, \dots, a_N \Big)\in (\Z^n)^{N+1}$ et 
$\alpha= (\alpha_0, \dots, \alpha_N) \in (K^\times)^{N+1}$, et 
supposons qu'il existe $a  \in \Expo(\cA)$ 
tel que pour toute place $v \in M_{{K}}$ le maximum des 
$|\alpha_i|_v$ pour $i=0,\dots,N$ est atteint au-dessus de $a$, autrement-dit
$$\max\{|\alpha_i|_v\, :\ 0 \le i \le N\} 
= \max\{ |\alpha_\ell|_v \, : \ 0 \le \ell \le N, a_\ell=a\}\enspace.
$$
Alors $\wh{\mu}^{\ess}(X_{\cA,\alpha}) = \wh{\mu}^\abs(X_{\cA,\alpha}^\circ) = \wh{h}(\alpha)$.  
\end{lem} 
\begin{demo} 
Notons $0 \le \ell_0, \dots, \ell_M \le N$ les indices pour lesquels $a_\ell=a$ 
et soit  $ \varpi: \P^N \to \P^M$ la projection $x \mapsto (x_{\ell_0} : 
\cdots : x_{\ell_M})$. 
Alors pour un point quelconque 
$
\xi = (\alpha_0\, s^{a_0} :\cdots : \alpha_N \, s^{a_N} )
\in X_{\cA,\alpha}^\circ $ on a 
$\varpi(\xi) = (\alpha_{\ell_0} \, s^a : \cdots : \alpha_{\ell_M} \,
s^a) = (\alpha_{\ell_0} : \cdots : \alpha_{\ell_M})
$ et donc 
$$
\wh h (\xi) \ge \wh h (\varpi(\xi)) = \wh h (\alpha_{\ell_0} : \cdots : \alpha_{\ell_M})\enspace,  
$$
 et $\wh h (\alpha_{\ell_0} : \cdots : \alpha_{\ell_M}) = \wh h (\alpha)$ gr{\^a}ce {\`a} l'hypoth{\`e}se du lemme, donc 
$\wh{\mu}^{\abs}(X_{\cA,\alpha}^\circ) \geq \wh{h}(\alpha)$.
En outre 
$\wh{h}(\alpha) \geq \wh{\mu}^{\ess}(X_{\cA,\alpha}) $ d'o\`u la conclusion. 
\end{demo}

\begin{demo}[D{\'e}monstration du th{\'e}or{\`e}me~\ref{densite}]
On peut supposer \spdg $N=3n+1$, puisque le cas général 
se déduit de celui-ci par immersion de $\P^{3n+1}$
comme un sous-espace standard de $\P^N$. 

\smallskip
Soient $d$ un nombre premier, $1\le k\le n$, $1\le f\le d-1$ 
des paramètres 
entiers à fixer ultérieurement et 
encore $q_0\geq\dots\geq q_n\geq 0$ des paramètres 
rationnels. 
Rappelons que $e_1,\dots,e_n$ désigne la base standard de $\Z^n$ et  
$S$ le simplexe standard $\Conv(0,e_1,\dots,e_n) \subset \R^n$, on pose 
$$
\begin{array}{l}
a_i:=d\,e_i \hspace{24mm} \mbox{ pour } i=1,\dots,n \enspace   \\[4mm]
b_i:=\left\{\begin{array}{ll}
            (d-1)\,e_i \enspace  \quad 
&\mbox{ pour } 1\le i \le k-1 \enspace \\[2mm] 
            f\,e_i \enspace  \quad & \mbox{ pour } i=k \enspace\\[2mm] 
            e_i \enspace  \quad & \mbox{ pour }  k+1 \le i \le n \enspace
          \end{array} \right. 
\end{array} 
$$
puis 
$$
\begin{array}{l}
\cA:= (0, a_1 \dots, a_n, 0, a_1, \dots, a_n, b_1, \dots,b_n)  
\in (\Z^n)^{3n+2} \enspace, \\[2mm] 
\alpha := (\underbrace{1, \dots\dots\dots\dots, 1}_{n+1\ \mbox{\footnotesize fois}}, 
2^{q_0}, 2^{q_1}, . .\dots, 2^{q_n}, 
\underbrace{2^{q_0}, \dots, 2^{q_0}}_{k\ \mbox{\footnotesize fois}}, 
\underbrace{1, \dots, 1}_{n-k\ \mbox{\footnotesize fois}}) \in (\Qbar^\times)^{3n+2}
\end{array}
$$
et on considère $X \subset \P^{3n+1}$ la variété torique 
associée au couple $(\cA,\alpha)$ ainsi défini. 
On a $L_\cA = \Z^n$ et 
$Q_\cA = d\,S$, donc la dimension et le degré de cette variété sont égaux à $n$ et $n!\,\Vol_n(Q_\cA) = d^n$
respectivement, et 
si $\ell$ est un dŽnominateur commun de $q_0,\dots,q_n$,
cette variété est définie sur l'extension kummerienne 
$K:=\Q(2^{\frac{1}{\ell}})$. 
De plus, il rŽsulte du lemme~\ref{mus toriques} et du 
lemme~\ref{tous-sur-a} appliquŽ ˆ toutes les faces de~$Q_\cA$
\begin{equation}\label{minessexemple}
\wh{\mu}_i(X) = q_{i-1}\,\log(2)\enspace,\quad i=1,\dots,n+1
\enspace,
\end{equation}
le $i$-ème minimum se réalisant 
sur la face $\Conv(a_{i-1},\dots,a_n)$. En outre, on estime la hauteur en utilisant la formule dans le th{\'e}or{\`e}me \ref{thmprincipal}. On décompose le polytope de base en 
$$
Q_\cA= Q_1 \cup Q_2 \cup Q_3
$$
avec $Q_1:= {\rm Conv}(0,b_1,\dots,b_k,a_{k+1},\dots,a_n)$, 
$Q_2:= {\rm Conv}(b_1,\dots,b_k,a_k,\dots,a_n)$ et 
$Q_3:= Q_\cA\setminus (Q_1 \cup Q_2)$.  
Pour $n=3$ on a une figure du type montrŽ ci-dessous~:

\begin{figure}[htbp]

\vtop to72mm{$$\raise-85mm\hbox{

\begin{picture}(0,0) 

\put(-170,55){\epsfig{file=figures/tetra.eps, height=60mm}} 
\put(-40,180){\raise9mm\hbox{$Q_1$}\epsfig{file=figures/tetra1.eps, height=29mm}} 
\put(60,150){\raise9mm\hbox{$Q_2$}\epsfig{file=figures/tetra2.eps, height=29mm}} 
\put(50,60){\raise9mm\hbox{$Q_3$}\epsfig{file=figures/tetra3.eps, height=29mm}} 

\end{picture}

}$$}

\end{figure}

\vspace{15mm}

Pour toute place finie $v$ de $K$ 
on vérifie $\vartheta_{\cA,\tau_{\alpha v}} \equiv 0$; 
la formule pour $\hnorm(X)$ se rŽduit alors
aux contributions des places archimŽdiennes. 
Pour  $v\in M_K^\infty$, la décomposition du domaine 
considérée sépare l'intégrale en trois morceaux
$I_i:=\int_{Q_i} \vartheta_{\cA,\tau_{\alpha\, v}}(u)du$ 
$(i=1,2,3)$.
On vérifie d'abord  que la fonction $\vartheta_{\cA,\tau_{\alpha v}}$
est {\it linéaire} sur chacun des simplexes $Q_1$ et $Q_2$. 
L'intŽgrale d'une fonction linŽaire sur un simplexe 
étant  Žgale au volume du simplexe multipliŽ par la 
valeur moyenne de la fonction, on trouve 
\begin{eqnarray*}
I_1 &= &
{\Vol_n(Q_1)} \, 
\frac{((k+1)\, q_0 + q_{k+1} + \dots + q_n)\,\log(2)}{n+1} 
\\[0mm]
&= &
\frac{f\, (d-1)^{k-1} \, d^{n-k}}{(n+1)!} \, 
((k+1)\, q_0 + q_{k+1} + \dots + q_n)\,\log(2) \enspace,
\\[2mm]
I_2 &= &
{\Vol_n(Q_2)} \, 
\frac{(k\, q_0 + q_{k} + \dots + q_n)\,\log(2)}{n+1} \\[0mm]
&= & \frac{(d-f)\, (d-1)^{k-1} \, d^{n-k}}{(n+1)!} \, 
(k\, q_0 + q_{k} + \dots + q_n)\,\log(2) \enspace.  
\end{eqnarray*}
En outre on estime brutalement la troisième intégrale
$$
0 \leq I_3
\leq \Vol_n(Q_3) \cdot \max(\vartheta_{\cA,\tau_{\alpha\, v}}) 
= \frac{d^n - (d-1)^{k-1} \, d^{n-k+1}}{n!}\, q_0 \, \log(2)
\enspace, 
$$
puisque $ \Vol_n(Q_3) = \Vol_n(Q_\cA) - \Vol_n(Q_1 \cup Q_2)= 
d^n- (d-1)^{k-1}\,d^{n-k+1}$.
Il s'ensuit
$$
(n+1)!\, \int_{Q_\cA} \vartheta_{\cA,\tau_{\alpha\, v}}(u)du 
= (n+1)!\, (I_1+I_2+I_3) =  
d^n \, (\theta+E(\arith, v)) 
$$
avec
$$\begin{array}{rcl}
\theta &= &\bigg(1-\frac{1}{d}\bigg)^{k-1}
\bigg(\frac{d-f}{d}\, (k\, q_0 + q_{k} + \dots + q_n) +
\frac{f}{d}\, ((k+1)\, q_0 + q_{k+1} + \dots + q_n) \bigg)\,\log(2)\\
&= &\bigg(1-\frac{1}{d}\bigg)^{k-1}\left((kq_0+q_k+\cdots+q_n) + \frac{f}{d}(q_0-q_k)\right)\,\log(2)
\end{array}$$
et 
$$ 0\le E(\arith,v)= \frac{(n+1)!}{d^n}\, I_3 \le 
  (n+1)\, \bigg(1-\bigg(1-\frac{1}{d}\bigg)^{k-1}\bigg) \, q_0 \, \log(2)
 \le \frac{(n+1)\, (k-1)}{d}\, q_0 \, \log(2) \enspace.
$$
On déduit alors du théorème~\ref{thmprincipal}
\begin{equation} \label{radisnoir}
0 \le \theta - \frac{\wh{h}(X)}{\deg(X)} \le 
\frac{(n+1)\, (k-1)}{d}\, q_0 \, \log(2) \enspace.
\end{equation}
Maintenant on fixe les paramètres~: d'abord on prend  
$\ell:=\lceil\log(2)\,\varepsilon_1^{-1}\rceil +1$, 
et pour $ 0 \le i \le n$ on pose 
$
q_i := \frac{1}{\ell}{\left\lfloor \frac{\ell\cdot \mu_{i+1}}{\log(2)} \right\rfloor}$
de sorte que  $\wh\mu_{i+1}(X) = q_i \,\log(2)$ satisfait
\begin{equation} \label{citronnelle}
0 \le \mu_{i+1} - q_i \,\log(2) 
<\varepsilon_1
\end{equation}
comme voulu. 
On vŽrifie
$$
(q_0+\dots+q_{n})\log(2) \leq \mu_1+\dots+\mu_{n+1} \leq \nu \leq (n+1)\mu_1 -\varepsilon_1 \leq (n+1)q_0\log(2) 
$$
et on fixe $ 1\le k \le n $ tel que 
$$
(k \, q_0+ q_{k}+\cdots+q_{n})\log(2) \leq 
\nu \leq ((k+1)\, q_0+q_{k+1}
+ \cdots+q_{n})\log(2) \enspace.
$$
Soit $\lambda\in[0,1]$ tel que 
$$\nu =(1-\lambda) \, 
(k \, q_0+ q_{k}+\cdots+q_{n})\log(2) 
+ \lambda \, ((k+1)\, q_0+q_{k+1}
+ \cdots+q_{n})\log(2) \enspace.
$$ 
On prend alors $1\le f \le d-1$ 
tel que 
\begin{equation} \label{escargot}
|\lambda-f/d|\le 1/d
\end{equation}
 et on considère le $\theta=\theta(d)$ associé à ces paramètres. 
On vŽrifie facilement
$$\begin{array}{rcl}
|\nu - \theta| &\kern-3pt\leq &\kern-3pt\left|\lambda - \frac{f}{d}\right|(q_0-q_k)\log(2)  + \left(1-\left(1-\frac{1}{d}\right)^{k-1}\right) \cdot ((k+1)q_0+q_{k+1}+\cdots+q_n)\log(2)\\[2mm]
&\kern-3pt\leq &\kern-3pt\displaystyle\frac{\mu_1}{d} + \frac{(k-1)(n+1)\mu_1}{d} = \frac{n^2\mu_1}{d}
\enspace.\end{array}$$
Finalement on obtient 
le rŽsultat cherchŽ en sommant avec l'inégalité~(\ref{radisnoir})
$$
\left|\frac{\wh{h}(X)}{\deg(X)} - \nu\right| \le 
\frac{(n+1)(k-1)}{d}\, \mu_1 + \frac{(n+1)(k-1)+1}{d}\, \mu_1 \leq
\frac{2n^2}{d}\, \mu_1
\enspace.
$$
En prenant, gr‰ce au postulat de Bertrand, $d$ un premier entre $2n^2 \, \mu_1 \, \varepsilon_2^{-1}
$ et $4n^2\mu_1\varepsilon_2^{-1}$ on arrive au résultat annoncé. 
\end{demo}

Ce résultat montre que dŽjˆ dans le cadre torique,
toute configuration possible des minimums $\wh{\mu}_1(X), \dots, \wh{\mu}_{n+1}(X)$ se réalise,  et 
l'encadrement~(\ref{zhang}) est optimal en toute dimension, 
lorsque le degrŽ de $X$ et celui du corps de définition
sont assez grands. 
Toutefois, notre exemple prŽsente une codimension minimale $N-n=2n+1$ 
de l'ordre de la dimension de la variŽtŽ produite. 
La question se pose donc de savoir ce qu'il en est pour les 
variŽtŽs de petite codimension. 
Dans le cas de codimension $1$ on a le rŽsultat suivant qui laisse ouverte la possibilitŽ de raffinements de~(\ref{zhang}) en petit codimension.

\begin{prop}\label{cashyper}
Soit $X\subset\P^{N}$ une hypersurface {\em torique}
d'Žqua\-tion homogne 
minimale $f_X=x^{b}- \lambda \, x^{c}\in \Qbar[x_0, \dots, x_N]$, 
alors
$$
 \wh{\mu}^{\ess}(X) 
= \frac{h(1:\lambda)}{\deg(f_X)}
= \frac{\wh{h}(X)}{\deg(X)} 
\enspace,\quad
\wh{\mu}_2(X)=\cdots=\wh{\mu}_{N}(X)=0\enspace.$$ 
\end{prop}

\begin{demo} 
Les contributions des 
places archimŽdiennes ˆ la hauteur normalisŽe de $X$
sont Žgales aux mesures de Mahler des conjuguŽes de $f_X$
et donc $\wh{h}(X) = h(1:\lambda)$, {\it voir}~\cite{PS04}, 
suite de l'exemple III.7. 
Par la proposition~\ref{pasteque} on a
$\lambda=\alpha^{b-c}$ pour tout point $\alpha \in X^\circ$.
En choisissant $\alpha$ de hauteur normalisŽe aussi proche
du minimum essentiel $\wh{\mu}^{\ess}(X)$ que l'on veut,
 on a 
$$
\frac{\wh{h}(X)}{\deg(X)} 
= \frac{h(1:\alpha^{b-c})}{\deg(f_X)}
= \frac{h(\alpha^{b}:\alpha^{c})}{\deg(f_X)}
\leq  \wh{h}(\alpha)
\leq \wh{\mu}^{\ess}(X)+\varepsilon
$$ 
pour tout $\varepsilon>0$. 
Avec la formule~(\ref{zhang}) de l'introduction on obtient $
\wh{\mu}_1 (X)+ \cdots + \wh{\mu}_{n+1} (X) 
\le  {\wh{h}(X)}/{\deg(X)} \le 
\wh{\mu}_1 (X)$, ce qui entraîne l'énoncé. 
\end{demo}

Il serait très intéressant d'expliciter les minimums successifs d'une
variété torique quelconque, en généralisant à la fois le lemme~\ref{tous-sur-a} 
et la proposition~\ref{cashyper}. 
Grâce au lemme~\ref{mus toriques} on sait qu'il suffit de le faire 
pour le minimum essentiel.

Une question liée est celle de construire explicitement des 
points de $X_{\cA,\alpha}^\circ$ de hauteur comparable au
 minimum essentiel. 
De plus, on peut se demander s'il existe un point 
dans la variété réalisant ce minimum essentiel, \cad
de savoir s'il existe un point de hauteur minimale. 
Dans le mme ordre d'idŽe, existe-t-il un point de $X_{\cA,\alpha}$ dont la hauteur normalisŽe 
soit Žgale au quotient $\displaystyle {\wh{h}(X_{\cA,\alpha})}/{\deg(X_{\cA,\alpha})}$?

\smallskip

Pour conclure ce paragraphe, explicitons l'encadrement 
pour le quotient hauteur-sur-degré qui découle de la formule pour la hauteur d'une variété torique: 

\begin{prop}\label{casa} 
Soient $\cA \in (\Z^n)^{N+1}$ et $\alpha\in(\overline{\Q}^\times)^{N+1}$, alors 
$$
\wh{h}(\alpha) - n \, \hnorm \Big( \alpha_j^{-1} \, : \, a_j \in \Faces_0(Q_\cA) \Big) \le \frac{\wh{h}(X_{\cA, \alpha})}{\deg(X_{\cA, \alpha})} \le (n+1) \, \wh{h}(\alpha)
\enspace,$$
où $a_j$ parcours l'ensemble $\Faces_0(Q_\cA)$ des sommets de $Q_\cA$. 
\end{prop} 

\begin{demo}
Soit $K$ le corps de définition de $\alpha$ et $v\in M_K$. 
On a $\max( \vartheta_{\cA, \tau_{\alpha \, v}}) = \log\max\{|\alpha_0|_v,\dots, |\alpha_N|_v\}=: \log(||\alpha||_v)$, ce
qui entraîne 
$$
 \int_{Q_\cA} \vartheta_{\cA, \tau_{\alpha \, v}} (u) \, du
\le \Vol_n(Q_\cA) \,\log (\Vert\alpha\Vert_v) 
= \frac{1}{n!}\, \deg(X_{\cA,\alpha}) \,\log (\Vert\alpha\Vert_v) 
$$ 
et donc $\wh{h}(X_{\cA,\alpha}) \le (n+1) \, \deg(X_{\cA,\alpha}) \, \wh{h}(\alpha)$, ce qui établit la majoration. 
Pour la minoration, soit $v\in M_K$ et posons 
$$
m_v: = \min_{u\in Q_\cA} ( \vartheta_{\cA, \tau_{\alpha \, v}}(u))=\min \Big\{ \max \{ \log|\alpha_i|_v  \, : \ 0\le i \le N, 
a_i=a \}  \, :\ 
a\in \Faces_0(Q_\cA) \Big\}
$$
et considérons le polytope 
$Q_v:= \Conv\Big((a_i,\log|\alpha_i|_v), (a_i,m_v) \, : \ i=0,\dots,N\Big)\subset \R^{N+1}$. 
Soit $0\le \ell\le N$ tel que $\log|\alpha_{\ell}|_v$ soit maximal, 
alors 
$Q_v \supset \Conv \Big( (a_{\ell}, \log|\alpha_{\ell}|_v), 
Q_\cA \times \{m_v\} \Big) $
et donc 
$ 
\Vol_{n+1} (Q_v)  \ge   \frac{1}{n+1} \, 
\Big( \log \Vert\alpha\Vert_v - m_{v} \Big) \, \Vol_{n} (Q_\cA)$.
On en déduit
$$
\int_{Q_\cA} \vartheta_{\cA, \tau_{\alpha \, v}} (u) \, du 
=   \Vol_{n+1} (Q_v) + m_{v} \, \Vol_{n} (Q_\cA) 
\ge  \bigg( \frac{1}{n+1} \, \log( \Vert\alpha\Vert_v)
+ \frac{n}{n+1} \, m_{v} \bigg) \, \Vol_{n} (Q_\cA) \enspace, 
$$
d'où
$\wh{h}(X_{\cA,\alpha}) \ge \deg(X_{\cA,\alpha}) \, \Big(\wh{h}(\alpha)  - n\, \wh{h}(\alpha_j^{-1}:a_j\in \Faces_0(Q_\cA))\Big)$
car $
\sum_{v\in M_K}\frac{[K_v:\Q_v]}{[K:\Q]}(-m_v) \leq \wh{h}(\alpha_j^{-1}:a_j\in \Faces_0(Q_\cA))$. 
\end{demo} 

Notons qu'en remplaant le point $\alpha$ par un point 
de $X_{\cA,\alpha}^\circ$ de hauteur aussi proche  que l'on veut
du minimum essentiel, on retrouve 
simplement la majoration de~(\ref{zhang}); par contre on obtient une minoration différente.

%

\typeout{Poids de Chow et hauteur des diviseurs monomiaux} 

\section{Poids de Chow et hauteur des diviseurs monomiaux} \label{diviseurs}

Dans ce paragraphe on consid{\`e}re l'intersection d'une vari{\'e}t{\'e} torique
avec un diviseur monomial de $\P^N$. 
On montrera comment dans cette situation, le th{\'e}or{\`e}me de B{\'e}zout 
pour les poids de 
Chow (Th{\'e}or{\`e}me~\ref{BezoutpoidsChow} de l'introduction et 
ThŽorme~\ref{Bezoutvar}) 
s'explicite comme la
d{\'e}composition poly\'edrale d'un certain volume, et dans cette situation peut se 
démontrer de façon indépendante.  
En combinant ceci avec la formule pour la hauteur d'une variété torique 
(thŽorme~\ref{thmprincipal}), on obtient 
un th{\'e}or{\`e}me de B{\'e}zout arithm{\'e}tique pour 
la hauteur normalisŽe du cycle intersection
d'une vari{\'e}t{\'e} torique avec un diviseur monomial. 

\medskip 

Soient $\geom \in (\Z^n)^{N+1}$ tel que $L_{\geom}=\Z^n$  
et $b \in \Z^{N+1}$, notons $X_\geom \subset \P^N$ et $\div(x^b) \in 
\Div(\P^N)$ la vari{\'e}t{\'e} torique et le diviseur monomial associ{\'e}s. 
On s'intŽressera au cycle intersection 
d{\'e}coup{\'e} sur $X_\geom$ par le mon{\^o}me $x^b$; 
le lemme ci-dessous explicite ce cycle.

Pour chaque hyperface $F\in F_{n-1}(Q_\cA)$ on considère 
la variŽtŽ 
$X_{\cA,F}  \subset \P^N$, adhérence de Zariski de l'orbite 
associ{\'e}e $X_{\cA, F}^\circ$.
On considère aussi l'hyperplan
d'appui $H_F \subset \R^n$ de cette face et 
$H_{F}^\Z:=H_F \cap \Z^n$; fixons également un point quelconque $a_F \in F$.  
Notons $L_{\geom, F}$ 
le $\Z$-module engendr{\'e} par les diff{\'e}rences
des {\'e}l{\'e}ments de $\geom\cap F$, qui est un sous-r{\'e}seau de $H_F^\Z -a_F $ d'indice  
$$
i(\cA;F) := [ H_{F}^\Z-a_F : L_{\geom,F} ] \enspace.
$$
Rappelons que $v_F \in \Z^n$ d{\'e}signe le plus petit vecteur entier, orthogonal {\`a} $H_F$ et dirig{\'e} vers l'int{\'e}rieur de $Q_{\geom}$. On pose $M_\cA(b) = b_0\,a_0 + \cdots + b_N\,a_N \in \Z^n$ et $D:= \deg(x^b) =\sum_{j=0}^N b_j$. 

\begin{lem} \label{oignon} 
Avec les notations ci-dessus, on a 
$$
X_\geom \cdot \div(x^b)
= \sum_{F \in \faces_{n-1}(Q_\geom)} 
 \langle M_\cA(b)  -D\, a_F,v_F\rangle \, i(\cA;F) \,
[X_{\cA, F} ] 
\ \in Z_{n-1}(\P^N)
$$
o{\`u} $\langle\cdot,\cdot\rangle$ d{\'e}signe le produit scalaire 
ordinaire de $\R^n$ et $[X_{\cA,F}]$ le cycle défini par $X_{\cA,F}$. 
\end{lem} 

En particulier, un cycle  $[X_{\cA,F}]$ intervient comme composante de 
$X_\cA\cdot \div(x^b)$ si et seulement si $M_\cA(b)$ 
n'appartient pas à l'hyperplan  $H_F+(D-1)a_F$.

\begin{demo}
Gr{\^a}ce {\`a} la d{\'e}composition en orbites~(\ref{decomposition}) 
on vérifie que le cycle considéré est supporté par la r{\'e}union
des orbites de codimension 1, celles-ci correspondant aux hyperfaces de $Q_\cA$, soit   
$$
\bigcup_{F \in \faces_{n-1}(Q_\geom)} X_{\cA, F} \enspace.  
$$
Explicitons les multiplicités correspondantes.  
Par lin{\'e}arit{\'e} on peut se ramener \spdg
au cas o{\`u} $b=e_i$ est un des vecteurs
de la base standard de $\R^{N+1}$ pour un certain $0\le i \le N$, \cad
$x^b=x_i$. Fixons une hyperface $F$ 
et soit $a_{j} \in F$ un sommet quelconque. 
Considérons le c{\^o}ne 
dual de l'angle de $Q_\geom$ en $a_j$: 
$$
\sigma:=  \Big\{ u \in \R^n \, : \, \langle u, a_k -a_j \rangle \ge 0 
\mbox{ pour }  k=0, \dots, N \Big\} \subset \R^n \enspace,  
$$
et soit $U_\sigma$ la vari{\'e}t{\'e} torique affine 
correspondante,
d{\'e}finie par 
$$
U_\sigma := \Spec(\Qbar[S_\sigma])
$$
o{\`u} $S_\sigma:= \sigma^\vee \cap \Z^n$
est le semi-groupe des points entiers dans l'angle
de $Q_\geom$ en $a_j$~\cite[\S~1.3]{Ful93}.
Consid{\'e}rons encore l'application naturelle 
$ N: U_\sigma \to (X_\geom)_{x_j}  \subset (\P^N)_{x_j}$
donn{\'e}e par l'inclusion d'alg{\`e}bres
$$
\Qbar[(X_\geom)_{x_j}] =  \Qbar[s^{a_0 - a_j}, \dots, s^{a_N-a_j}] 
\hookrightarrow \Qbar[S_\sigma] = \Qbar[U_\sigma] \enspace. 
$$
La vari{\'e}t{\'e} $U_\sigma$ est {\it normale} car $S_\sigma$ est un
semi-groupe satur{\'e}, et de ce fait  $N$ est le  morphisme de 
normalisation de la carte affine $(X_\geom)_{x_j}$~\cite[Cor.~13.6]{Stu96}. 

\smallskip 

Soit $\rho$ l'ar{\^e}te du c{\^o}ne $\sigma$ duale de la face $F$, notons
$V(\rho)$ la cl\^oture dans $U_\sigma$ de l'orbite correspondante~\cite[\S~3.1]{Ful93}. On a un diagramme commutatif
$$
\begin{array}{ccc} 
\Qbar[(X_\geom)_{x_j}] & {N_*}\atop{\displaystyle \hookrightarrow} & \Qbar[U_\sigma] \\[2mm]
\downarrow & & \downarrow \\[2mm] 
\Qbar[(X_{\cA, F})_{x_j}] & \hookrightarrow & \Qbar[V(\rho)]
\end{array} 
$$
qui implique $N^{-1}(X_{\cA, F}) = V(\rho)$ et $\deg(N|_{V(\rho)}) = [ H_{F}^\Z : L_{\geom,F} ]=i(\geom; F)$.  
Le mon{\^o}me $\chi:= N^*(x_i)= s^{a_i-a_j} \in \Qbar[s_1^{\pm 1}, \dots, s_n^{\pm 1}]$ d{\'e}finit une fonction rationnelle $U_\sigma \dashrightarrow  \Qbar$. On a $\Ord_{X_{\cA, F}} (x_i) = \deg(N|_{V(\rho)}) \,
\Ord_{V(\rho)}(\chi)$ puisque $U_\sigma$ est normal~\cite[Exerc.~1.2.3.]{Ful84}. 
On en dŽduit 
$$\begin{array}{rcl}
\longueur_{\Qbar[(X_{\cA, F})_{x_j}]} \Big(\Qbar[(X_\cA)_{x_j}] / (x_i) \Big)
&= &\Ord_{X_{\cA, F}} (x_i)\\[0mm]
&= &\deg(N|_{V(\rho)}) \, \Ord_{V(\rho)}(\chi)
= i(\geom;F)  \, \Ord_{V(\rho)}(\chi)
\enspace.\end{array}$$
Finalement, le lemme de \cite[\S~3.3, p. 61]{Ful93} entra{\^\i}ne 
$
\Ord_{V(\rho)}(\chi) = \langle a_i-a_j , v_F \rangle
$
car $v_F$ est le g{\'e}n{\'e}rateur du semi-groupe $\rho \cap \Z^n \cong \N$, 
d'o{\`u} 
$$m(X_\cA\cdot\div(x_i);X_{\cA,F}) =  \longueur_{\Qbar[(X_{\cA, F})_{x_j}]} \Big(\Qbar[(X_\cA)_{x_j}] / (x_i) \Big) = \langle a_i-a_j , v_F \rangle \, i(\geom;F)
\enspace.$$ 
\end{demo}

Le th{\'e}or{\`e}me de B{\'e}zout g{\'e}om{\'e}trique 
\begin{equation}\label{Bezoutgeomtorique}
\deg \Big(X_\geom \cdot \div(x^b) \Big)= 
D \, \deg(X_\cA) 
\end{equation} 
s'interprte en termes de d{\'e}composition du volume du polytope $Q_\geom$~: 
on a $\deg(X_{\geom}) = n! \, {\rm Vol}_n (Q_{\geom})$ et
\begin{eqnarray*}
\deg(X_{\cA, F}) &=&(n-1)! \mu_{\cA(F)}(F) \\[2mm]
&=& \frac{(n-1)!}{\Vol_{n-1}((H_F-a_F)/L_{\cA,F})}
\, {\Vol_{n-1}(F)} 
= \frac{(n-1)!}{{i(\cA;F) \, ||v_F||_2}}
\, {\Vol_{n-1}(F)} \enspace, 
\end{eqnarray*}
{\`a} cause de la normalisation de la forme volume $\mu_{\cA(F)}$ sur $H_F$ et du fait que
$$
\Vol_{n-1}(H_F/L_{\cA,F}) = 
i(\cA; F) \cdot \Vol_{n-1}(H_F / H_{F}^\Z) = i(\cA;F)\cdot ||v_F||_2 \enspace,
$$ 
cons{\'e}quence de la formule de Brill-Gordan. On a encore 
$$ 
\langle M_\cA(b)- Da_F,v_F\rangle = 
\varepsilon(b,F)\, \Vert v_F\Vert_2 \, 
{\rm dist}(M_\cA(b),H_{F}+(D-1)a_F) 
\enspace,
$$ 
o $\varepsilon(b,F)=+1$ si $M_\cA(b)$ et $v_F$ sont d'un mme c™tŽ de l'hyperplan $H_F+(D-1)a_F$ et $\varepsilon(b,F)=-1$ sinon. 
Combiné avec le lemme~\ref{oignon} cela donne
\begin{eqnarray*}
\deg \Big(X_\geom \cdot \div(x^b) \Big)&
= &\sum_{F\in \faces_{n-1}(Q_\cA)} \langle M_\cA(b)-Da_F,v_F\rangle\,i(\cA,F)\, \deg(X_{\cA,F}) \\[0mm]
&= &\sum_{F\in \faces_{n-1}(Q_\cA)} 
\varepsilon(b,F)\, (n-1)! \, \Vol_{n-1}(F) \, 
{\rm dist}(M_\cA(b),H_{F}+(D-1)a_F)  \\[0mm]
&=& \sum_{F\in \faces_{n-1}(Q_\cA)} 
n! \, \varepsilon(b,F)\,\Vol_n(\Conv(F+(D-1)a_F,M_\cA(b)))\enspace.
\end{eqnarray*}
L'identit{\'e}~(\ref{Bezoutgeomtorique}) ci-dessus (multipli{\'e}e par $n!^{-1}$) 
se traduit ainsi en la d{\'e}composition de volume 
\begin{equation}\label{radis}
\displaystyle\sum_{F \in \faces_{n-1}(Q_\geom)}  \varepsilon(b,F)\,
\Vol_n(\Conv(F+(D-1)a_F,M_\cA(b)))
=D\, \Vol_n(Q_{\geom})
\enspace.
\end{equation}
La figure suivante illustre cette d{\'e}composition, pour 
$\geom= \Big( (1,1),(1,0),(3,1),(2,2),(0,2)\Big) \in (\Z^2)^6$ et 
$b= (1,0,0,0,0)$ (donc $M_\cA(b)=(1,1)=a_0$ et $D=1$)~:   

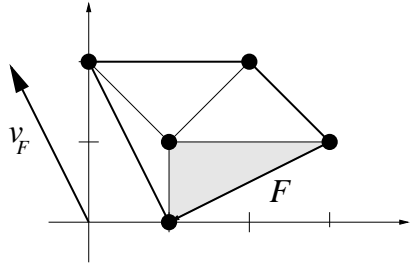
\vspace{-3mm} 

\begin{figure}[htbp]
$$\epsfig{file=figures/bezout-geom.eps, height= 35mm}$$
\end{figure}

\vspace{-3mm}

Soit $\tau = (\tau_0,\dots,\tau_N)\in\Z^{N+1}$. 
Pour le cas $\tau \in (\N^\times)^{N+1}$, le théorème
de Bézout pour les poids de Chow (théorème~\ref{Bezoutvar})
s'écrit 
dans les notations du \S~\ref{Chow}
\begin{equation}\label{litchi} 
e_\tau(X_\geom \cdot  \div(x^b)) = D \, e_\tau(X_{\geom}) - \sum_{Y \in \irr( \initial_\tau(X_\cA))}  m((X_{\cA})_\tau\cdot \div(\lambda_{\tau}^*(x^b)) ;\iota(Y)) \cdot \deg(Y) 
\enspace. 
\end{equation}
On va expliciter les termes intervenant dans cet enoncé. 
Soit 
$$Q_{\geom, \tau} = \Conv((a_0,\tau_0), \dots, (a_N,\tau_N)) \subset \R^{n+1}$$ le polytope associ{\'e} au couple  ($\geom, \tau)$, 
dont la {\it toiture} $E_{\geom, \tau}$ (\cad l'enveloppe sup{\'e}rieure) 
s'envoie bijectivement sur $Q_\cA$ par la projection $\R^n \times \R \to \R^n$. 
On appelle {\it pan} de la toiture toute face de dimension $n$ de
$E_{\cA,\tau}$. 
De m{\^e}me, on appelle {\it mur} tout polytope de la forme
$\Conv( Q_{\cA(F)}, F\times\{0\})$ pour une hyperface $F$ de $ Q_\cA$; 
dans le cas
$\tau \ge 0$ ceci est un des murs de la 
{\it maison}   $\cM_{\geom,\tau }:= \Conv(Q_{\cA,\tau} \cup (Q_{\cA} \times \{ 0\}))$, \cad une face de dimension $n$ se projetant 
sur une face de dimension $n-1$ de $Q_{\geom}$.

La d\'eformation torique $(X_\cA)_\tau$ est l'adhérence de Zariski de l'application 
$$
\G_m\times \G_m^n \to \P^1\times \P^N \enspace, \quad  
(t, s) \to \Big((1:t), (t^{\tau_0}s^{a_0}: \cdots: t^{\tau_N}s^{a_N})\Big) \enspace.
$$
C'est donc la variété torique bi-projective associée aux vecteurs
$(0,1)\in (\Z^1)^2$ et $(\tau, \cA)=\Big((\tau_0,a_0), \dots, (\tau_N,a_N)\Big)
\in (\Z\times \Z^n)^{N+1}$, {\it voir} \S~\ref{vartor}. 
Le couple de polytopes associé est donc 
${\mathbf 0}\times [0,1]$ et $Q_{\cA,\tau}$; 
comme cons{\'e}quence de la d{\'e}composition en orbites décrite au~\S~\ref{vartor}, on v{\'e}rifie que
les points de $(X_\cA)_\tau$ contenus dans l'hyperplan 
$\{(0:1)\}\times\P^N$ correspondent nŽcessairement aux couples de 
faces de la forme $(\{ ({\mathbf 0},1)\},F) $ avec 
$F \in \Faces_n(E_{\cA,\tau})$. 
On voit ainsi qu'il y a bijection entre les pans de $E_{\cA,\tau}$ et les orbites de $X_{\cA,\tau}$ contenues dans l'hyperplan $\{(0:1)\}\times\P^N$, en particulier le support de $X_{\cA,\tau} \cdot (\{(0:1)\}\times\P^N)$ est contenu dans
$$
\bigcup_{P \in \faces_{n}(E_{\geom,\tau})} X_{\cA, \tau, P} \enspace.
$$
L'identification $\iota: \P^N \to \{(0:1)\}\times\P^N \subset \P^1\times\P^N$ met en correspondance le cycle $X_{\cA,\tau} \cdot (\{(0:1)\}\times\P^N)$ et la variŽtŽ initiale $\init_\tau(X_\cA)$, on en d{\'e}duit qu'il y a bijection entre les composantes de $\init_\tau(X_\geom)$ et les pans de $E_{\cA,\tau}$. 

Pour chaque pan $P$ 
on considère son hyperplan
d'appui $H_P \subset \R^{n+1}$, 
posons $L_{\geom, \tau, P} \subset \Z^{n+1}$ 
le $\Z$-module engendr{\'e} par les diff{\'e}rences
des {\'e}l{\'e}ments de $(\tau,\geom) \cap P$. 
Modulo une translation, ce dernier    
est un sous-r{\'e}seau de $H_P^\Z:=H_P \cap \Z^{n+1}$ d'indice  
$$
i(\cA, \tau;P) := [ H_{P}^\Z : L_{\geom, \tau,P} ] \enspace.
$$
Avec ces notations, d'apr{\`e}s~\cite[formule (27), page~222]{Stu94} 
({\it voir} aussi~\cite[Thm. 5.3.]{KSZ92}),  on a
$$
X_{\cA,\tau} \cdot (\{(0:1)\}\times\P^N)=\iota(\init_\tau(X_\geom))
= \sum_{P \in \faces_n(E_{\geom,\tau})} 
i(\cA, \tau;P) 
\, [X_{\geom,\tau,P}] \enspace.     
$$
Pour chaque pan $P$ on note $(v_P,w_P) \in \Z^n\times\Z $ le plus
petit vecteur entier orthogonal au plan d'appui $H_P \subset \R^{n+1}$
tel que $w_P <0$. 
On note aussi  $(a_P,\tau_P)$ un point quelconque de $P$. 

\medskip

Le vecteur $\tau$ induit une {\it d{\'e}composition poly\'edrale  coh{\'e}rente}  
$\DPC_\tau(Q_\cA) $ du polytope de base $Q_\geom$. Les faces $S$ de dimension $n$ de cette d{\'e}composition sont en correspondance avec les pans de la toiture; pour $S \in\DPC_\tau(Q_\cA) $ on {\'e}crit $\Pan(S) \in \Faces_{n}(E_{\cA,\tau})$ pour le pan correspondant.

\begin{lem} \label{ratatouille} 
Soient $\tau\in \Z^{N+1}$ et $P \in \Faces_n(E_{\cA,\tau})$ un pan de la toiture de $Q_{\cA,\tau}$
alors, 
$$
m(X_{\cA,\tau} \cdot \div(\lambda_\tau^*(x^b)); X_{\cA,\tau,P}) = ( \langle M_\cA(b) - D\,a_P, v_P \rangle - D\, \tau_P\, w_P ) \, i(\cA,\tau; P) 
\enspace,$$
avec $\lambda_{\tau}^*(x^b) = (t^{\tau_0} x_0, \dots, t^{\tau_N}x_N)^b=t_0^{b_0\tau_0+\dots+b_N\tau_N}x^b$. 
\end{lem} 

\begin{demo} 
Cette d{\'e}monstration étant tout {\`a} fait analogue {\`a} celle du
lemme~\ref{oignon}, on n'indiquera que les pas principaux. 

\smallskip 
 
Par lin{\'e}arit{\'e} on peut se ramener au cas o{\`u} $b=e_i$ est un des vecteurs de la base standard de $\R^{N+1}$, 
donc $x^b=x_i$. Soit $(a_j,\tau_j)$ un {sommet} quelconque du pan $P$, on se place dans la carte affine $\A^1\times \A^N \subset \P^1\times\P^N$ correspondant {\`a} $t_1 \ne 0$ et $x_j \ne 0$ (puisqu'on veut calculer une multiplicit{\'e} le long d'une sous-vari{\'e}t{\'e} de $Z(t_0)$). Dans cette carte, l'application monomiale s'{\'e}crit
$$
\G_m\times \G_m^n \to \A^1\times \A^N \enspace, \quad \quad 
(t,s) \mapsto \Big(t; t^{\tau_j-\tau_0} \, s^{a_0-a_j}, \dots,
t^{\tau_j-\tau_N} \, s^{a_N-a_j} \Big) 
$$
et donc l'alg{\`e}bre de cette carte affine de   $X_{\cA,\tau}$ est 
$$
\Qbar[(X_{\cA,\tau})_{t_1,x_j}]=
\Qbar\Big[t, t^{\tau_j-\tau_0} \, s^{a_0-a_j}, \dots,
t^{\tau_j-\tau_N} \, s^{a_N-a_j} \Big] \enspace.
$$ 
La normalisation de cette alg{\`e}bre correspond au c{\^o}ne 
$$
\sigma:= \{(u,v) \in \R^n \times \R \, : \  
v \ge 0, \langle u,a_k-a_j\rangle + v \, (\tau_j-\tau_k) \ge 0, \ k=0, \dots, N\} \enspace. 
$$
Le reste de la d{\'e}monstration suit les lignes de celle du lemme~\ref{oignon}, en consid{\'e}rant la normalisation de 
$(X_{\cA,\tau})_{t_1,x_j}$ donn{\'e}e par le semi-groupe des points entiers $S_\sigma:= \sigma^\vee \cap \Z^n$. 

L'hyperplan $H_P \subset \R^{n+1}$ est un des plans d'appui du c{\^o}ne $\sigma^\vee$, car $P$ est un pan. 
Il d{\'e}finit donc une ar{\^e}te $\rho$ du c{\^o}ne dual $\sigma$, dont le semi-groupe $\rho \cap \Z^n$ est engendr{\'e} par $(v_P,-w_P) $. Comme $\lambda_{\tau}^*(x_i) =t^{\tau_i}x_i = t^{\tau_j}s^{a_i-a_j}$ sur la carte considŽrŽe, on en conclut
\begin{eqnarray*} 
m(X_{\cA,\tau} \cdot \div(\lambda_{\tau}^*(x_i)); X_{\cA,\tau,P}) &=& 
\langle (a_i,0)-(a_j, -\tau_j) , (v_P,-w_P)  \rangle \,
i(\cA,\tau; P) \\[2mm] 
&=& ( \langle a_i-a_j, v_P \rangle - \tau_j\, w_P ) \, i(\cA,\tau; P) \enspace. 
\end{eqnarray*}
\end{demo}

\begin{prop} \label{aneth}
Soit $\cA\subset (\Z^n)^{N+1}$ tel que $L_\cA=\Z^n$, $b\in\Z^{N+1}$ et $\tau \in \N^{N+1}$, alors l'\'egalit\'e~(\ref{litchi}) correspond terme \`a terme \`a la suivante, multipli\'ee par $(n+1)!$,
$$\begin{array}{rcl} 
\displaystyle D\,{\rm Vol}_{n+1}(\cM_{\cA,\tau}) &\kern-1mm= &\kern-2mm\displaystyle\sum_{F \in \faces_{n-1}(Q_\geom)} \varepsilon(b,F)\Vol_{n+1} (\Conv(M(F)+(D-1)(a_F,0),(M_\cA(b),0)))\\[8mm] 
&&\displaystyle\kern3mm+ \sum_{P \in \faces_n(E_{\cA,\tau})} \varepsilon(b,P)\Vol_{n+1}(\Conv(P+(D-1)(a_P,\tau_P), (M_\cA(b),0))) 
\end{array}$$ 
o $\varepsilon(b,F)=+1$ ({\it resp.} $\varepsilon(b,P)=+1$) si $M_\cA(b)$ et $v_F$ ({\it resp.} $(M_\cA(b),0)$ et $(v_P,w_P)$) sont d'un mme c™tŽ de l'hyperplan $H_F+(D-1)a_F$ ({\it resp.} $H_P+(D-1)(a_P,\tau_P)$) et $\varepsilon(b,F)=-1$ ({\it resp.} $\varepsilon(b,P)=-1$) sinon. 
\end{prop}

La figure suivante (pour $b=e_i$, $D=1$ et $M_\cA(b)=a_i$) illustre cette décomposition:  

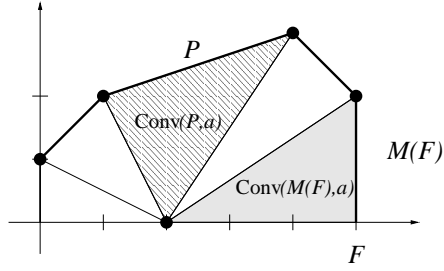
\begin{figure}[htbp]
\vbox to30mm{$$
\epsfig{file=figures/bezout-poids.eps, height= 35 mm}
$$}
\vspace*{-0mm} 
\end{figure}

\begin{demo}

Montrons comment~(\ref{litchi}) se traduit en la d{\'e}composition de l'int{\'e}grale de la fonction $\vartheta_{\cA,\tau}: Q_\cA \to \R$ param{\'e}trant la toiture $E_{\cA,\tau}$ de la proposition~\ref{aneth}. On a d'abord 
$$
e_\tau(X_\cA) = (n+1)!\, \int_{Q_\cA} \vartheta_{\cA, \tau}(u) \, du_1
\cdots du_n = (n+1)!\, \Vol_{n+1} (\cM_{\cA,\tau}) 
$$
puis, par le lemme~\ref{oignon} il vient
$$
e_\tau(X_\geom \cdot \div(x^b) ) = \sum_{F \in \faces_{n-1}(Q_\geom)} \langle M_\cA(b)  - D\, a_F,v_F\rangle \, i(\cA;F) \, e_\tau(X_{\cA, F})
\enspace,$$
et pour chaque face $F$ de $Q_\cA$  
$$
e_\tau(X_{\cA, F}) = \frac{n!}{i(\cA;F) \, ||v_F||_2}  \, 
\int_F \vartheta_{\cA, \tau} \, d\mu_{n-1} 
= \frac{n!}{i(\cA;F) \, ||v_F||_2}  \, \Vol_{n} (M(F))
$$
o{\`u} $M(F) \subset \R^{n+1}$ d{\'e}signe le mur de la maison $\cM_{\cA,\tau}$ au-dessus de $F$. Ainsi
$$\begin{array}{rcl} 
\langle M_\cA(b)-D\, a_F,v_F\rangle\,i(\cA;F)\,e_\tau(X_{\cA, F}) 
&\kern-2mm=&\kern-2mm n!\,\varepsilon(b,F)\dist(M_\cA(b),F+(D-1)a_F)\,\Vol_{n} (M(F))\\[3mm] 
&\kern-5.8cm=&\kern-3cm (n+1)!\,\varepsilon(b,F)\Vol_{n+1}(\Conv(M(F)+(D-1)(a_F,0),(M_\cA(b),0))) 
\enspace.\end{array}$$ 

Finalement, soit $Y$ une composante irr{\'e}ductible de
$\init_\tau(X_\cA)$ puis $S \in \DPC_\tau(Q_\geom)$ et $P:= \Pan(S) \in
\Faces_n(E_{\cA,\tau})$ les faces correspondantes dans la subdivision
et dans la toiture respectivement, on a 
$$
\deg(Y) =  \frac{n!}{i(\cA; S)} \Vol_n(S) = \frac{n!}{ i(\cA,\tau;P) 
 \, ||(v_P,w_P)||_2} \, \Vol_{n} (P) \enspace. 
$$ 
Le lemme~\ref{ratatouille} entra{\^\i}ne alors
$$\begin{array}{rcl} 
m(X_{\cA,\tau}\cdot \div(\lambda_{\tau}^*(x^b));\iota(Y)) \, \deg(Y)
&\kern-3mm= &\kern-3mm\frac{n!}{||(v_P,w_P)||_2} \,  \langle (M_\cA(b),0)\kern-1mm -\kern-1mm D\, (a_P,\tau_P), (v_P,w_P) \rangle  \, \Vol_{n} (P) \\[3mm] 
&\kern-4cm= &\kern-2.1cm(n+1)! \,\varepsilon(b,P) \Vol_{n+1}(\Conv(P+(D-1)(a_P,\tau_P), (M_\cA(b),0))) \enspace. 
\end{array}$$ 
En regroupant ces calculs on voit que l'identit{\'e}~(\ref{litchi}) (multipli{\'e}e par $\frac{1}{(n+1)!}$)  
se traduit dans la d{\'e}compo\-si\-tion cherchée.  
\end{demo}

Pour $S \in \DPC_\tau(Q_\cA)$ on définit  $\theta_{\tau, S} (b) \in \R$ l'unique r{\'e}el tel que $\big( M_\cA(b), \theta_{\tau, S} (b) \big) \in H_{\pan(S)}+(D-1)(a_{\pan(S)},\tau_{\pan(S)})$. Le lemme suivant explicite cette quantitŽ.

\begin{lem} \label{theta}
Soit $S \in \DPC_\tau(Q_\geom) $ et  $a_{j_0}, \dots , a_{j_n} \in S $ des vecteurs affinement ind\'e\-pen\-dants. Soit $b \in \Z^{N+1}$ et $D:= \sum_{j=0}^N b_j$, alors
$$
\begin{array}{rll}
\theta_{\tau,S}(b)
\cdot \det  
\left[ 
\begin{array}{ccc} 
1 & \dots & 1 \\[1mm] 
a_{j_0,1} & \cdots & a_{j_n,1} \\[1mm] 
\vdots & \ddots & \vdots \\[1mm] 
a_{j_0,n} & \cdots & a_{j_n,n} 
\end{array}
\right] 
&= &
-\det \left[
\begin{array}{ccccccc}
1       & \dots  &1           & D \\[1mm] 
a_{j_0,1} & \dots  &a_{j_n,1} & M_\cA(b)_1 \\[1mm] 
\vdots  & \ddots &\vdots      & \vdots\\ [1mm] 
a_{j_0,n} & \dots  &a_{j_{n},n} & M_\cA(b)_n \\[1mm]
\tau_{j_0}     & \dots  & \tau_{j_n}      & 0 
\end{array}
\right]\enspace. 
\end{array}
$$
\end{lem} 

En triangulant $S$ par des simplexes et en utilisant la relation entre d\'eterminant et volume, on peut r{\'e}Žcrire ceci sous la forme (on pose $P=\Pan(S)$)~:
\begin{equation}\label{parmentier}
\begin{array}{rcl}
\theta_{\tau,S}(b) \, \Vol_n(S) &\kern-3mm= &\kern-3mm(n+1) \, \varepsilon(b,P) \Vol_{n+1}\Big( \Conv(P+(D-1)(a_{P},\tau_{P}), (M_\cA(b),0)) \Big)\\
&\kern-3mm= &\kern-3mm m(X_{\cA,\tau}\cdot\div(\lambda_{\tau}^*(x^b));\iota(Y)) \, \deg(Y)
\end{array}
\end{equation}
qui s'interpr\`ete comme l'{\'e}galit{\'e} des volumes montrŽs dans la figure suivante (lorsque $D=1$ et en posant $a=M_\cA(b)$)~: 

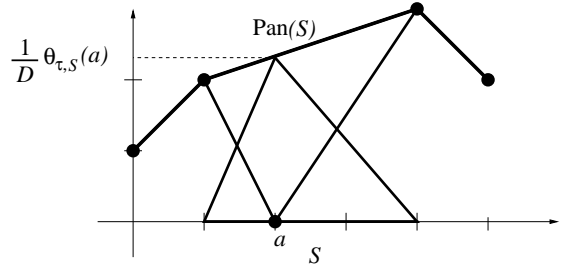
\vspace{-3mm}
\begin{figure}[htbp]
$$
\epsfig{file=figures/triangles.eps, height= 35 mm}
$$
\end{figure}

\vspace*{-7mm} 

\begin{demo} 
C'est un calcul direct de l'intersection des espaces lin{\'e}aires $H_{\pan(S)}+(D-1)(a_{\pan(S)},\tau_{\pan(S)})$ et $\{M_\cA(b)\}\times\R$~: on a 
$\theta_{\tau,S}(b) = \sum_{i=0}^n v_i \, \tau_{j_i}$ o{\`u} $v= (v_0, \dots, v_n) \in \R^{n+1} $ est l'unique 
solution du syst{\`e}me lin{\'e}aire 
$$
 \sum_{i=0}^n v_i = D\quad, \quad \quad 
 \sum_{i=0}^n v_i \, a_{j_i}=M_\cA(b) \enspace. 
$$
On r{\'e}sout ce syst{\`e}me par les formules de Cramer et on trouve ainsi 
$$
\theta_{\tau,S}(b)\cdot 
\det 
\begin{bmatrix}
1 & \dots & 1 \\[1mm] 
a_{j_0,1} & \cdots & a_{j_n,1} \\[1mm] 
\vdots & \ddots & \vdots \\[1mm] 
a_{j_0,n} & \cdots & a_{j_n,n} 
\end{bmatrix}
\kern-1mm = \displaystyle \sum_{i=0}^{n} 
\tau_{j_i} 
\,  \det
\begin{bmatrix}
1 &\dots &1 &\kern-1mm D &\kern-1mm 1 &\dots &1 \\[1mm] 
a_{j_0,1} &\dots &a_{j_{i-1},1} &\kern-1mm M_\cA(b)_1 &\kern-1mm a_{j_{i+1} ,1} &\dots &a_{j_{n},1}\\[1mm] 
\vdots & \ddots &\vdots &\vdots &\vdots & \ddots &\vdots\\[1mm] 
a_{j_0,n} &\dots &a_{j_{i-1},n} &\kern-1mm M_\cA(b)_n &\kern-1mm a_{j_{i+1},n} &\dots &a_{j_{n},n}
\end{bmatrix}
$$
qui est le d\'eveloppement du d\'eterminant dans le membre droite de l'ŽnoncŽ, par rapport ˆ la derni{\`e}re ligne. 
\end{demo}

Ceci permet d'{\'e}crire l'identit{\'e}~(\ref{litchi}) de la fa{\c c}on
suivante: 

\begin{prop}\label{Bezoutnorm}
Soit $\geom \in (\Z^n)^{N+1}$, $b \in \Z^{N+1}$ et $\tau = (\tau_0, \dots, \tau_N) \in \R^{N+1}$, posons $D:= \sum_{j=0}^N b_j$, alors
$$
e_\tau \Big( X_\geom \cdot \div(x^b) \Big) = D \, e_\tau(X_{\geom}) - n! \, \sum_{S \in \dpc_\tau(Q_\geom)} \theta_{\tau, S}(b) \, \Vol_n (S)
$$
o{\`u} la seconde somme porte sur les faces $S$ de dimension $n$ de la d{\'e}composition poly\'edrale $\DPC_\tau(Q_{\geom})$. 
\end{prop}

\begin{demo}
Le cas $\tau\in \N^{N+1}$ r\'esulte de la proposition~\ref{aneth} et de~(\ref{parmentier}). Ceci s'étend successivement à $\tau \in \Z^{N+1}$ à cause de l'invariance de la formule par remplacement de $\tau$ par $\tau+c\cdot (1,\dots, 1)$ avec $c\in \N$. 
Puis, la formule s'étend à $\tau\in \Q^{N+1}$ par homog\'en\'eit\'e et finalement à $\tau\in \R^{N+1}$ par continuité. 
\end{demo}

Avec la notation~(\ref{defvaluationrelative}) on peut encore Žcrire cet ŽnoncŽ sous la forme~:
$$w_{X_{\cA,\tau}}(x^b) = - \sum_{S\in\dpc_\tau(Q_\geom)} \theta_{\tau, S}(b).\frac{\Vol_n (S)}{\Vol_n (Q_\cA)}\enspace,$$
dont on vŽrifie, par continuitŽ et homogŽnŽitŽ, la validitŽ pour tout $\tau\in\R^{N+1}$.
Dans le cas o{\`u} $\div(x^b)$ est effectif, \cad quand $b\in\N^{N+1}$, on
a 
$
\theta_{\tau,S} (b) \ge \tau_0 \, b_0 + \cdots + \tau_N \, b_N
$
pour tout $S \in \DPC_\tau(Q_\cA)$ {\`a} cause de la concavit{\'e} de la
toiture du polytope $Q_{\cA,\tau}$, et donc
$$
(\tau_0 \, b_0 + \cdots + \tau_N \, b_N)
\, \deg(X_\cA) \le n!\, \sum_{S \in \dpc_\tau(Q_\cA)} \theta_{\tau,S} (b)
\, \Vol_n(S) \enspace, 
$$
ainsi
\begin{equation} \label{sauterelle} 
e_\tau \Big( X_\geom \cdot \div(x^b) \Big) 
\le D \, e_\tau(X_{\geom}) - (\tau_0 \, b_0 + \cdots + \tau_N \, b_N)
\, \deg(X_\cA)  \enspace. 
\end{equation} 
Alternativement, on peut d{\'e}montrer cette in{\'e}galit{\'e} par 
application directe du th{\'e}or{\`e}me~\ref{BezoutpoidsChow} et de l'exemple~\ref{poidshyper}. 

\medskip 
On en d\'eduit 
un th{\'e}or{\`e}me de B{\'e}zout arithm{\'e}tique {\it exact} pour la hauteur normalisŽe de  l'intersection d'une vari{\'e}t{\'e} torique avec un diviseur monomial:

\begin{cor} \label{capre} 
Soit $K$ un corps de nombres, $\cA \in (\Z^n)^{N+1} $ tel que $L_\cA=\Z^n$, $\alpha \in (K^\times)^{N+1}$ et $b \in \Z^{N+1}$. Posons 
$\tau_{\alpha\,v}= (\log|\alpha_0|_v, \dots,\log|\alpha_N|_v) $ pour toute place $v\in M_K$
et  $D:= \sum_{j=1}^N b_j$, alors
\begin{eqnarray*}
\wh{h} (X_{\arith}\cdot \div(x^b)) 
&=& D \wh{h}(X_\arith)  - n! \, \sum_{v\in M_K} \frac{[K_v : \Q_v]}{[K:\Q]} \sum_{S\in\dpc_{\tau_{\alpha\,v}}(Q_\geom)} \theta_{\tau_{\alpha\,v},S}(b)\Vol_{n} (S) \\
&=& D \wh{h}(X_\arith)  + \bigg( \sum_{v\in M_K} \frac{[K_v : \Q_v]}{[K:\Q]}
w_{X_{\cA,\tau_{\alpha\, v}}}(x^b) \bigg) \deg(X_\arith)
\enspace.
\end{eqnarray*}
En particulier, si $\div(x^b)$ est effectif (\cad $b\in\N^{N+1}$) on a $\wh{h} ( X_{\arith}\cdot \div(x^b) ) \le  D \, \wh{h}(X_\arith)$. 
\end{cor}

\begin{demo} 
Pour l'identit\'e on remarque que $X_\arith\cdot \div(x^b) = \alpha(X_\cA\cdot \div(x^b))$. En sommant sur $v\in M_K$ l'\'egalit\'e de la proposition~\ref{Bezoutnorm} avec les coefficients $\frac{[K_v : \Q_v]}{[K:\Q]}$, on conclut gr\^ace au th\'eor\`eme~\ref{thmprincipal}. 

Pour Žtablir l'in{\'e}galit{\'e}, on a par~(\ref{sauterelle})
\begin{eqnarray*}
\wh{h} (X_{\arith}\cdot \div(x^b)) 
&\le &D \wh{h}(X_\arith)  - \sum_{v\in M_K} 
\frac{[K_v : \Q_v]}{[K:\Q]}
 (\tau_{\alpha\,v,0} \, b_0 + \cdots + \tau_{\alpha\,v,N} \, b_N) \cdot
\, \deg(X_\cA) \\[2mm] 
&\leq &D \wh{h}(X_\arith)
\end{eqnarray*} 
car $\displaystyle \sum_{v\in M_K} \frac{[K_v : \Q_v]}{[K:\Q]} \tau_{\alpha\,v,i}=0$ pour $i=0,\dots,N$, gr{\^a}ce {\`a} la formule du produit. 
\end{demo} 

Si l'on définit la hauteur de $x^b$ relative à la variété $X_\arith$ par la formule
$$
\hnorm_{X_\arith}(x^b)
:= \sum_{v\in M_K} \frac{[K_v : \Q_v]}{[K:\Q]}\, w_{X_{\cA,\tau_{\alpha\,v}}}(x^b) 
$$
le résultat préc\'edent se r\'eécrit 
$$
\wh{h} (X_{\arith}\cdot \div(x^b)) 
=D \wh{h}(X_\arith)  + \hnorm_{X_\arith}(x^b)\cdot \deg(X_{\cA,\alpha})
\enspace.$$

%

\let\livrefont=\sl

\typeout{References}


\begin{thebibliography}{Gor1873}

{\small

\bibitem[AD03]{AD}
F.~Amoroso, S.~David, 
{\it Minoration de la hauteur normalis{\'e}e dans un tore}, 
J. Inst. Math. Jussieu {\bf 2} (2003) 335-381. 

\bibitem[AD04]{AD04}
F.~Amoroso, S.~David, 
{\it Le probl{\`e}me de Lehmer en dimension sup{\'e}rieure. II,} 
pr{\'e}publication de l'université de Caen, 2004. 

\bibitem[Aud91]{Audin}
M.~Audin,
{\livrefont The topology of torus actions on symplectic manifolds}, 
Progress in Math. {\bf 93}, Birkh{\"a}user, 1991. 

\bibitem[Ber87]{Ber87} 
D. Bertrand, 
{\it Lemmes de z{\'e}ros et nombres transcendants}, 
S{\'e}minaire Bourbaki 1985/86,   Ast{\'e}risque  {\bf 145 \& 146} (1987) 21-44. 

\bibitem[BP88]{BP}
D.~Bertrand, P.~Philippon, 
{\it Sous-groupes algŽbriques de groupes algŽbriques commutatifs}, 
Illinois J. Math. {\bf 32} (1988) 263-280.

\bibitem[Cha89]{Cha89} 
M.~Chardin, 
{\it Une majoration de la fonction de Hilbert et 
ses cons{\'e}quences pour l'interpolation alg{\'e}brique},
Bull. Soc. Math. France {\bf 117} (1989) 305-318. 

\bibitem[CP99]{CP99} 
M.~Chardin, P.~Philippon, 
{\it RŽgularitŽ et interpolation}, 
J. Algebraic Geom. {\bf 8} (1999) 471-481. 

\bibitem[Cox01]{Cox01} 
D.~Cox, 
{\it Minicourse on toric varieties}, notes d'un cours donn{\'e} {\`a} l'université de Buenos Aires en 
Juillet 2001. T{\'e}l{\'e}chargeable {\`a}~{\tt http\hspace{-2mm}://www.amherst.edu/\~{ }dacox/}. 

\bibitem[CLO98]{CLO98} 
{D.~Cox, J.~Little, D.~O'Shea,} 
{\livrefont Using algebraic geometry,} 
Graduate Texts in Math.~{\bf 185}, Springer, 1998. 

\bibitem[Dav03]{Dav02}
{S. David,}
{\it On the height of subvarieties of groups varieties,}  
{\`a} para{\^\i}tre dans le J. Ramanujan Math. Soc.. 

\bibitem[DP99]{DP99} 
{S.~David, P.~Philippon,}
{\it Minorations des hauteurs normalis{\'e}es des sous-vari{\'e}t{\'e}s des
tores,} Ann. Scuola Norm. Sup. Pisa {\bf 28} (1999) 489-543. 

\bibitem[Don02]{Don02} 
{S.K.~Donaldson,}  
{\it Scalar curvature and stability of toric varieties},
J. Differential Geom. {\bf 62} (2002) 289-349.

\bibitem[ES96]{ES96} 
D. Eisenbud, B.~Sturmfels,
{\it Binomial ideals}, 
Duke Math. J. {\bf 84} (1996) 1-45. 

\bibitem[EF02]{EF02} 
{J.-H. Evertse, R.G. Ferretti,}  
{\it Diophantine inequalities on projective varieties},
Internat. Math. Res. Not. {\bf 25} (2002) 1295-1330.

\bibitem[Ewa96]{Ewa96} G. Ewald, 
{\livrefont Combinatorial convexity and algebraic geometry}, 
Graduate Texts in Math. {\bf 168}, Springer, 1996.

\bibitem[Fer03]{Fer03} 
{R.G. Ferretti,}  
{\it Diophantine approximation and toric
 deformations}, Duke Math. J. {\bf 118} (2003) 493-522. 

\bibitem[Ful84]{Ful84}
W. Fulton, 
{\livrefont Intersection theory,} 
Ergeb. Math. Grenzgeb. (3) {\bf 2}, Springer, 1984. 

\bibitem[Ful93]{Ful93} 
{W. Fulton,}
{\livrefont Introduction to toric varieties,} 
Ann. Math. Studies {\bf 131}, Princeton Univ. Press, 1993.


\bibitem[GKZ94]{GKZ94}
{I.M. Gelfand, M.M. Kapranov, A.V. Zelevinsky,}
{\livrefont Discriminants, resultants and multidimensional determinants,}
Birkh{\"a}user, 1994.


\bibitem[KSZ92]{KSZ92} 
M.M. Kapranov, B. Sturmfels, A.V. Zelevinsky,
{\it Chow polytopes and general resultants}, 
Duke Math. J. {\bf 67} (1992)  189-218. 


\bibitem[Mum77]{Mum}
D.~Mumford, 
{\it Stability of projective varieties}, 
Enseign. Math. {\bf 23} (1977) 39-110.

\bibitem[Phi91]{Phi91} 
{P.~Philippon,} 
{\it 
Sur des hauteurs alternatives, I},
   Math. Ann. {\bf 289}  (1991) 255-283. 

\bibitem[PS04]{PS04}
{P. Philippon, M. Sombra,} 
{\it Hauteur normalis{\'e}e des vari{\'e}t{\'e}s toriques projectives,} 
t{\'e}l{\'e}chargeable ˆ~{\tt http\hspace{-2mm}://fr.arxiv.org/abs/math.NT/0406476}, 38 pp.. 

\bibitem[PS05]{PS05}
{P. Philippon, M. Sombra,} 
{\it Géométrie diophantienne et variétés toriques,}
C. R. Math. Acad. Sci. Paris  {\bf 340}  (2005) 507-512.

\bibitem[PS06]{PS06}
{P. Philippon, M. Sombra,} 
{\it À propos du minimum essentiel des translatés de sous-tores,}
tapuscrit, 16 pp..

\bibitem[Rat04]{Rat04} 
N.~Ratazzi, 
{\livrefont Minoration de la hauteur de N{\'e}ron-Tate pour les points et les
  sous-vari{\'e}t{\'e}s: variations sur le probl{\`e}me de Lehmer},
th{\`e}se de Doctorat, Université de~Paris~VI, 2004. 

\bibitem[Rem01a]{Rem01a}
{G. R{\'e}mond,}
{\it \smash{\'E}limination multihomogne}, 
chapitre 5 de {\livrefont Introduction to algebraic independence theory}, 
Lecture Notes in Math. {\bf 1752} (2001) 53-81.

\bibitem[Rem01b]{Rem01b}
{G. R{\'e}mond,}
{\it G{\'e}om{\'e}trie diophantienne multiprojective}, 
chapitre 7 de {\livrefont Introduction to algebraic independence theory}, 
Lecture Notes in Math. {\bf 1752} (2001) 95-131.

\bibitem[Sch91]{Sch91} 
W.M. Schmidt,
{\livrefont Diophantine approximation and diophantine equations}, 
Lecture Notes in Math. {\bf 1467}, Springer, 1991. 

\bibitem[Som05]{Som02}
{M. Sombra,} {\it Minimums successifs des variétés toriques projectives,}
J. Reine Angew. Math.~{\bf 586}  (2005) 207-233.

\bibitem[Stu94]{Stu94} 
B. Sturmfels,
{\it On the Newton polytope of the resultant}, 
J. Algebraic Combin. {\bf 3} (1994), 207-236. 

\bibitem[Stu96]{Stu96} 
B.~Sturmfels, 
{\livrefont Gr{\"o}bner bases and convex polytopes}, 
Amer. Math. Soc., 1996. 


\bibitem[Zha95]{Zha95}
{S.-W. Zhang,} 
{\it Positive line bundles on arithmetic varieties}, J. Amer. Math. Soc. {\bf 8}  (1995)
187-221. 

}

\end{thebibliography}
\end{document}